\documentclass[final,3p,times]{elsarticle}
\usepackage{amssymb}
 \usepackage{amsthm}
\usepackage{graphicx}
\usepackage{mathrsfs}
\usepackage[colorinlistoftodos]{todonotes}
\usepackage[colorlinks=true, allcolors=blue]{hyperref}
\usepackage{amsmath}
\usepackage{multirow}
\usepackage{amsfonts, geometry}
\usepackage{bbold}
\newtheorem{theorem}{Theorem}[section]
\newtheorem{corollary}{Corollary}[section]

\newtheorem{lemma}{Lemma}[section]

\theoremstyle{dgdef}
\newtheorem{definition}{Definition}[section]
\newtheorem{example}{Example}[section]

\begin{document}
\begin{frontmatter}

\title{Global and local  existence  of solutions for a novel type of parabolic  Kirchhoff system with singular term \tnoteref{label1}}

\author[1]{Ahmed Aberqi }
\ead{ahmed.aberqi@usmba.ac.ma}
\author[2]{ Abdesslam Ouaziz }
\ead{abdesslam.ouaziz@usmba.ac.ma}
\author[3,4]{ Maria Alessandra Ragusa  \corref{cor1}}\ead{ maragusa@dmi.unict.it }

\cortext[cor1]{Corresponding author}
%\address[]{
%\textbf{"Dedicated to professor Lars-Erik Persson"}
%}
\address[1]{LSATE, National School of Applied Sciences, Sidi Mohamed Ben Abdellah University,   BP 1796, 30000 Fez, Morocco.}

\address[2]{LSATE, Faculty of Sciences Dhar EL Mahraz, Sidi Mohamed Ben Abdellah University,  BP 1796, 30000 Fez, Morocco.}
\address[3]{Dipartimento di Matematica e Informatica, Universita di Catania, Viale Andrea Doria, 6, 95125 Catania, Italy}
\address[4] {Faculty of Fundamental Science, Industrial University of Ho Chi Minh City, Vietnam}
\begin{abstract}
%%%%
In this paper, we investigate solutions for a fractional system involving a novel class of Kirchhoff functions and logarithmic nonlinearity:
\begin{equation*}
 \left\{\begin{array}{lll}
\displaystyle
\mathfrak{u}_{t}+\mathcal{K}\left([\mathfrak{u}]_p^s\right) \mathscr{L}_p^s u=\vert  \mathfrak{v} \vert^{\sigma }\vert  \mathfrak{u} \vert^{\sigma-2} u \log | \mathfrak{u} \mathfrak{v}|, \, \, & \mbox{in}\quad &\mathcal{U} \times[0, T),\\
\mathfrak{v}_t+\mathcal{K}\left([\mathfrak{v}]_q^s\right) \mathscr{L}_q^s \mathfrak{v}=\vert  \mathfrak{u} \vert^{\sigma }|\mathfrak{v}|^{\sigma-2} \mathfrak{v} \log | \mathfrak{u} \mathfrak{v}|, & \text { in } & \mathcal{U} \times[0, T), \\
\mathfrak{u}(\mathrm{x}, t)=\mathfrak{v}(\mathrm{x}, t)=0, & \text { in } & \partial \mathcal{U} \times[0, T), \\
\mathfrak{u}(\mathrm{x}, 0)=\mathfrak{u}_0(\mathrm{x}), \mathfrak{v}(\mathrm{x}, 0)=\mathfrak{v}_0(\mathrm{x}), & \text { in } & \mathcal{U},
\end{array}%
\right.
\end{equation*}
where  $\mathcal{K}$ is Kirchhoff function, and $\mathscr{L}_{p}^{s}$ is the fractional $p-$ Laplacian operator. We prove the existence of a weak solution using the Faedo-Galerkin method under suitable assumptions on the Kirchhoff function. We investigate the   finite-time blow-up  and global existence of solutions based on  critical, subcritical, and supercritical initial energy levels.
Subsequently, we establish the stabilization of the solution with positive initial energy by applying Komornik's integral inequality.
\end{abstract}
\date{\today}
\end{frontmatter}
\textit{Math. Subj. Classif. (2020):} Primary 58J10, Secondary 58J20, 35J66	
%% main text

\section{Introduction}
In this work, we are concerned with study the local and global existence solution for parabolic system with Kirchhoff function and logarithmic nonlinearity. To be more precise, we study the following system:
\begin{equation}\label{problem}
 \left\{\begin{array}{lll}
\displaystyle
\mathfrak{u}_{t}+\mathcal{K}\left([\mathfrak{u}]_p^s\right) \mathscr{L}_p^s u=\vert  \mathfrak{v} \vert^{\sigma }\vert  \mathfrak{u} \vert^{\sigma-2} u \log | \mathfrak{u} \mathfrak{v}|, \, \, & \mbox{in}\quad &\mathcal{U} \times[0, T),\\
\mathfrak{v}_t+\mathcal{K}\left([\mathfrak{v}]_q^s\right) \mathscr{L}_q^s \mathfrak{v}=\vert  \mathfrak{u} \vert^{\sigma }|\mathfrak{v}|^{\sigma-2} \mathfrak{v} \log | \mathfrak{u} \mathfrak{v}|, & \text { in } & \mathcal{U} \times[0, T), \\
\mathfrak{u}(\mathrm{x}, t)=\mathfrak{v}(\mathrm{x}, t)=0, & \text { in } & \partial \mathcal{U} \times[0, T), \\
\mathfrak{u}(\mathrm{x}, 0)=\mathfrak{u}_0(\mathrm{x}), \mathfrak{v}(\mathrm{x}, 0)=\mathfrak{v}_0(\mathrm{x}), & \text { in } & \mathcal{U},
\end{array}%
\right.
\end{equation}
where $\mathcal{U} \subset \mathbb{R}^{N}$ being a bounded domain with Lipschitz boundary $\partial \mathcal{U},$  $s \in(0,1),$    $p, $  $q,  $ and $\sigma$ are three numbers positives verify some conditions that will be mentioned later.
The main  operator $\mathscr{L}_{p}^{s}$ is
 the so-called the fractional p- Laplacian operator given by
\begin{equation*}
\mathscr{L}_{p}^{s} \mathfrak{u}= 2 \lim _{\varepsilon \rightarrow 0^{+}} \int_{\mathcal{U} \backslash \mathfrak{B}_{\varepsilon}(\mathrm{x})}   \frac{ \vert \mathfrak{u}(\mathrm{x})-\mathfrak{u}(\mathrm{y})\vert^{p-2}  (\mathfrak{u}(\mathrm{x})-\mathfrak{u}(\mathrm{y}))}{ \vert \mathrm{x}-\mathrm{y}\vert^{N+s p}} \textrm{d}\mathrm{y}.
\end{equation*}
Furthermore, we considere  the continuous  function
$\mathcal{K}:[0, \infty)\to  (0, \infty)$ which  satisfy the following assumptions:\begin{enumerate}
\item[ $(\mathcal{H} _{1})$] The function $\mathcal{K}$ is non- decreasing and  continuous  on $(0, \infty). $
\item[ $(\mathcal{H} _{2})$] For  a given  a constant $\beta \geq 0,$  $z \mapsto  \frac{\mathcal{K}(z)}{z^{\beta}}$ is decreasing on $(0, \infty).$
\end{enumerate}
Below, we present several examples of models that satisfy our assumptions.
\begin{example}
\begin{enumerate}
\item [$\bullet$] $\mathcal{K}(t)= a+bt^{c},$  where  $a, $ $b,$ $c>0.$\\
\item [$\bullet$] $\mathcal{K}(t)= \log(1+t).$
\end{enumerate}
\end{example}

Furthermore, we suppose   numbers   $p, $  $q,$ $\beta$ and $\sigma$ satisfy the following inequality:
\begin{equation}\label{equation1}
1<p<q<\sigma< \sigma+1<p_s^* \text { and } \max (2, p(\beta+1), q(\beta+1))<\sigma<2 \sigma \leq \min \left\{p\left(1+\frac{2}{N}\right), q\left(1+\frac{2}{N}\right), p_s^*, q_s^*\right\},
\end{equation}
with $p_s^*$ and $q_s^*$ are the criticals Sobolev exponent given by
$$
p_s^*= \begin{cases}\frac{N p}{N-s p}, & \text { if } N>s p, \\ +\infty, & \text { if } N \leq s p,\end{cases}
$$
and
$$
q_s^*= \begin{cases}\frac{N q}{N-s q}, & \text { if } N>s q, \\ +\infty, & \text { if } N \leq s q.\end{cases}
$$
For the sake of simplicite, we note
\begin{align*}
[\mathfrak{u}]_p^s= \frac{1}{p}
\int_{\mathcal{U} \times \mathcal{U} } \frac{|\mathfrak{u}(\mathrm{x})-\mathfrak{u}(y)|^p}{\vert \mathrm{x}- \mathrm{y}\vert^{N+s p}} d \mathrm{x} d \mathrm{y},
 \text {  and  } [\mathfrak{v}]_p^s= \frac{1}{q}
\int_{\mathcal{U} \times \mathcal{U} } \frac{|\mathfrak{v}(\mathrm{x})-\mathfrak{v}(y)|^p}{\vert \mathrm{x}- \mathrm{y}\vert^{N+s q}} d \mathrm{x} d \mathrm{y}.
\end{align*}

As of late, great attention has been paid on the study of equations involving double phase operator and  variational problems with growth conditions with variable exponent. This kind of problem arises in  numerous applications, such as,  population dynamics, phase transition phenomena,  continuum mechanics, typical outcome of stochastically stabilization
of L{\'e}vy  processes, image processing, electro-rheological fluids and  thermo-rheological fluids \cite{Aboulaich, Applebaum, Antontsev,  Caffarelli, Chen,  Rajagopal}.
Kirchhoff in \cite{kirchhoff} was first introduced the following model of elastic string \begin{equation}\label{kirchhof}
\rho h \mathfrak{u}_{tt} +\delta \mathfrak{u}_t= \left( p_0+\frac{Eh}{L} \mathcal{K}(\mathfrak{u})\right)u_{\mathrm{x}\mathrm{x}} +f,
\end{equation}
where $\displaystyle \mathcal{K}(\mathfrak{u})=\int_{0}^{L} \mathfrak{u}_{\mathrm{x}}^2 d\mathrm{x},$ accounts for the effect of transversal oscillations. Additionally, Alembert  in \cite{Ghisi} developped the model \eqref{kirchhof}, and  he proposed the below equation for the free vibrations:
\begin{align}\label{Alembert}
\varepsilon \mathfrak{u}_{t t}^{\varepsilon}+\mathfrak{u}_t^{\varepsilon}-\mathcal{K}\left(\int_{\mathcal{U}}\left|\nabla \mathfrak{u}^{\varepsilon}\right|^2 \mathrm{~d} \mathrm{x}\right) \Delta \mathfrak{u}^{\varepsilon}=f(\mathrm{x}, t).
\end{align}
 Setting $\epsilon=0$ in equation \eqref{Alembert} yields the Kirchhoff-type problem:
\begin{align}\label{Alembert1}
\mathfrak{v}_t-\mathcal{K}\left(\int_{\mathcal{U}}|\nabla \mathfrak{v}|^2 \mathrm{~d} \mathrm{y}\right) \Delta \mathfrak{v}=f(\mathrm{y}, t).
\end{align}
These models  have attracted significant research interest in recent years, as they encompass various physical phenomena.  Before highlighting the novelty of this work, we provide a review of the existing literature.
In \cite{Han}, the authors studied the $p-$Kirchhoff type models describing motion of a non-stationary fluid given by
$$\mathfrak{u}_t-\mathcal{K}(\Vert \mathfrak{u}\Vert_{p}^p)\left(-\Delta_p\right)(\mathfrak{u})=\vert \mathfrak{u}\vert^{\alpha-2}\mathfrak{u}.$$  For the degenerate case, authors in   \cite{Pan}, have proved the solvability of  the fractional version
$$\mathfrak{u}_{t}-\mathcal{K}(\left[ \mathfrak{u} \right]_{s,p}^p)\left(-\Delta_p^s\right)(\mathfrak{u})=\vert \mathfrak{u}\vert^{\alpha-2}\mathfrak{u}.$$

A third fascinating aspect of our problem is the presence of a logarithmic nonlinearity term. Indeed, considering the following parabolic equation:
\begin{align*}
\mathfrak{u}_{t}= \nabla \mathfrak{u}+ \vert \mathfrak{u} \vert ^{q-2}\log \vert \mathfrak{u} \vert, \quad  \mathfrak{u}: \mathbb{R}^{n}\times (0, \infty) \to \mathbb{R}, \quad q, n\geq 2,
\end{align*}
which it shows up in a lot of physical applications, such
as  theory of superfluidity, nuclear physics,  diffusion phenomena, and transport. See \cite{Zloshchastiev} for more details. Due to these reasons, the following model

\begin{equation*}
 \left\{\begin{array}{lll}
\displaystyle
\mathfrak{u}_{t}+\left(- \triangle \right)^{s} \mathfrak{u}
=\vert  \mathfrak{u} \vert^{r-2} u \log \vert \mathfrak{u}\vert , \, \, & \mbox{in}\quad &\mathcal{U} \times[0, T),\\
\mathfrak{u}(\mathrm{x}, t)=0, & \text { in } & \partial \mathcal{U} \times[0, T), \\
\mathfrak{u}(\mathrm{x}, 0)=\mathfrak{u}_0(\mathrm{x}),  & \text { in } & \mathcal{U},
\end{array}%
\right.
\end{equation*}
has been studied quite extensively; for example, the authors in previous papers \cite{Aberqi, aberqi2022, dAvenia, Le, Tian}
studied the  asymptotic behavior, global existence,  and finite blow up of solutions. In this context, many authors has been studied the blow-up and global existence of solutions to the following problem:
\begin{equation*}
 \left\{\begin{array}{lll}
\displaystyle
\mathfrak{u}_{t}+\left(- \triangle \right) \mathfrak{u}
=g(  \mathfrak{u} ), \, \, & \mbox{in}\quad &\mathcal{U} \times[0, T),\\
\mathfrak{u}(\mathrm{x}, t)=0, & \text { in } & \partial \mathcal{U} \times[0, T), \\
\mathfrak{u}(\mathrm{x}, 0)=\mathfrak{u}_0(\mathrm{x}),  & \text { in } & \mathcal{U},
\end{array}%
\right.
\end{equation*}
where $g$ is a function satisfies some conditions. See the references \cite{aberq, Ball, Payne, Tan, Liu1} for more details.
Inspired by previpously works, in this present paper, we study the effect of logarithmic nonlinearity term end a new class of Kirchhoff function that satsfies  $(\mathcal{H} _{1})$ and $(\mathcal{H} _{2}).$ For this, we use potentiel well method what introduced by Sattinger and Payne  \cite{Sattinger, Payne} and Galerkin's method, we show the following parabolic system:

\begin{equation*}\label{problem}
 \left\{\begin{array}{lll}
\displaystyle
\mathfrak{u}_{t}+\mathcal{K}\left([\mathfrak{u}]_p^s\right) \mathscr{L}_p^s u=\vert  \mathfrak{v} \vert^{\sigma }\vert  \mathfrak{u} \vert^{\sigma-2} u \log | \mathfrak{u} \mathfrak{v}|, \, \, & \mbox{in}\quad &\mathcal{U} \times[0, T),\\
\mathfrak{v}_t+\mathcal{K}\left([\mathfrak{v}]_q^s\right) \mathscr{L}_q^s \mathfrak{v}=\vert  \mathfrak{u} \vert^{\sigma }|\mathfrak{v}|^{\sigma-2} \mathfrak{v} \log | \mathfrak{u} \mathfrak{v}|, & \text { in } & \mathcal{U} \times[0, T), \\
\mathfrak{u}(\mathrm{x}, t)=\mathfrak{v}(\mathrm{x}, t)=0, & \text { in } & \partial \mathcal{U} \times[0, T), \\
\mathfrak{u}(\mathrm{x}, 0)=\mathfrak{u}_0(\mathrm{x}), \mathfrak{v}(\mathrm{x}, 0)=\mathfrak{v}_0(\mathrm{x}), & \text { in } & \mathcal{U},
\end{array}%
\right.
\end{equation*}
has a weak solution.
 However, under some conditions on the intial data, we  show the existence of global solution for system \eqref{problem}.
An one of the main tools, Komormik's integral inequality. In addition, we show the global solution decays algebraically or exponentially depending on $q,$  $p,$ $s,$   $\mathfrak{u}_0,$
 and $\mathfrak{v}_0.$ Otherwise, we study blow-up phenomenna in subscritical initial energy case, and we  give an upper bounded for the maximal existence time.
The difficulty here is apparently the coupling term $\log(\mathfrak{u} \mathfrak{v}),$ in our problem,  and
  the absence of logarithmic Sobolev inequality which seems not to exist logarithmic Sobolev inequality concerning the fractional p-Laplacian again. \\
The remainder of this paper is structured as follows.  In Section \ref{Fractional Sobolev spaces}, we give some preliminary results that will be used throughout the paper. In sections  \ref{Potential well theory},  \ref{Global Existence and Asymptotic Behavior},  and \ref{Blow-up Phenomena},   we will prove the main results.
\section{Preliminaries}\label{Fractional Sobolev spaces}
\subsection*{ Fractional Sobolev spaces}
In this subsection, we review key properties of fractional Sobolev spaces. For more details, see \cite{bucur, iezza}. Let $\mathcal{U}\subset \mathbb{R}^{N}$ be a bounded domain, and $\partial \mathcal{U}$ its boundary. For $s\in (0, 1),$ and
$p \in[1 ;+\infty),$ we define Sobolev space fractional $W^{s, p}(\mathcal{U})$ as follows:
$$
W^{s, p}(\mathcal{U})=\left\lbrace \mathfrak{v} \in L^p(\mathcal{U}): \frac{|\mathfrak{v}(\mathrm{x})-\mathfrak{v}(\mathrm{y})|}{ \vert \mathrm{x}- \mathrm{x}\vert^{\frac{n}{p}+s}} \in L^p(\mathcal{U} \times \mathcal{U})\right\rbrace,
$$
i.e, an intermediary Banach space between $L^p(\mathcal{U})$ and $W^{1, p}(\mathcal{U}),$ endowed with the natural norm
$$
\Vert\mathfrak{v}\Vert_{W^{*, p}(\mathcal{U})}=\left(\int_{\mathcal{U}}\vert \mathfrak{v}(\mathrm{x})\vert ^p  d\mathrm{x} +[\mathfrak{v}]_{W^{s, p}(\mathcal{U})}^p\right)^{\frac{1}{p} },
$$
where  $[\mathfrak{v}]_{W^{s, p}(\mathcal{U})}$ is the Gagliardo semi-norm  given by
$$
[\mathfrak{v}]_{W^{s, p}(\mathcal{U})}=\left(\iint_{\mathcal{U} \times \mathcal{U}}     \frac{\vert \mathfrak{v}(\mathrm{x})-\mathfrak{v}(\mathrm{y})\vert^{p} }{\vert \mathrm{x}- \mathrm{y}\vert^{N+sp}}d\mathrm{x} d\mathrm{y}\right) ^{\frac{1}{p}},
$$
 and the Banach space $$
\Vert \mathfrak{v}\Vert_{W_0^{s, p}(\mathcal{U})}:=\left(\int_{\mathbb{R}^N} \int_{\mathbb{R}^N} \frac{|\mathfrak{v}(\mathrm{x})-\mathfrak{v}(y)|^p}{\vert \mathrm{x}- \mathrm{y}\vert^{N+s p}} d \mathrm{x} d \mathrm{y}\right)^{\frac{1}{p}},
$$
endowed with  the Gagliardo semi-norm:
$$
\Vert \mathfrak{v}\Vert_{W_0^{s, p}(\mathcal{U})}:=\left(\int_{\mathbb{R}^N} \int_{\mathbb{R}^N} \frac{|\mathfrak{v}(\mathrm{x})-\mathfrak{v}(y)|^p}{\vert \mathrm{x}- \mathrm{y}\vert^{N+s p}} d \mathrm{x} d \mathrm{y}\right)^{\frac{1}{p}}.
$$

\begin{theorem}( See \cite{iezza, Guo})
For $1\leq p <\infty$ and $0<s<1$ with $sp<N.$
We have the embedding $$W^{s, p}\left(\mathcal{U}\right) \hookrightarrow  L^q\left(\mathcal{U}\right) \text {  and }   W_{0}^{s, p}\left(\mathcal{U}\right) \hookrightarrow  L^q\left(\mathcal{U}\right)$$ which are respectively continuous for $q\in[1, p^{\star}_{s}),$  $q\in[1, p^{\star}_{s}]$ and both compacts for
$q\in[1, p^{\star}_{s}).$
\end{theorem}

$$
\mathrm{W}_0^{s, p}(\mathcal{U}):=\left\lbrace \mathfrak{v} \in \mathrm{W}^{s, p}\left(\mathbb{R}^{\mathrm{N}}\right),\,\  \mathfrak{v}=0 \text { a.e. in } \mathbb{R}^{\mathrm{N}} \backslash \mathcal{U}\right\rbrace
$$
and the Banach norm in the space $\mathrm{W}_0^{s, p}(\mathcal{U})$ is the Gagliardo semi-norm:
$$
\Vert \mathfrak{v}\Vert_{W_0^{s, p}(\mathcal{U})}:=\left(\int_{\mathbb{R}^N} \int_{\mathbb{R}^N} \frac{|\mathfrak{v}(\mathrm{x})-\mathfrak{v}(y)|^p}{\vert \mathrm{x}- \mathrm{y}\vert^{N+s p}} d \mathrm{x} d \mathrm{y}\right)^{\frac{1}{p}}.
$$
\begin{theorem}( See \cite{iezza}:  Poincar{\'e} inequality )
Let $\mathfrak{v} \in \mathrm{W}_0^{s, p}(\mathcal{U}).$ Then, there exists a positive constant $c>0,$ such that
$$
c^{-1}\Vert \mathfrak{v}\Vert_{\mathrm{W}^{s, p}\left(\mathbb{R}^N\right)} \leq\Vert \mathfrak{v}\Vert_{\mathrm{W}_0^{s, p}(\mathcal{U})} \leq c\Vert \mathfrak{u}\Vert_{\mathrm{W}^{s, p}\left(\mathbb{R}^{\mathrm{N}}\right)},
$$
\end{theorem}
  \begin{lemma}\label{lemma1} (See \cite{adams}
  If $1<p_{0}<p \mu<p_{1}<\infty.$ Then, we have
  $$\Vert \mathfrak{v}\Vert_{p\mu}\leq \Vert \mathfrak{v}\Vert_{p_{0}}^{1-\mu} \Vert \mathfrak{v}\Vert_{p_{1}}^{\mu},   $$
  for any $\mathfrak{v}\in L^{p_{1}}(\mathcal{U})$ with $\mu\in (0,1)$ and  $\frac{1}{\mu p}=\frac{1-\mu}{ p_{0}}+ \frac{\mu}{p_{1}}.$
\end{lemma}
\begin{definition} (See \cite{Han}) \textbf{(Maximal existence time)}\\
The maximal existence time of a weak solution $ \mathfrak{v}$ to the system  \eqref{problem} is defined as:

 \begin{enumerate}
 \item If $\mathfrak{v}$ exists for all $0 \leq t<+\infty,$ then $T_{\max }=+\infty.$
 \item If there exits a $t_0 \in(0,+\infty)$ such that $\mathfrak{v}$ exists for $0 \leq t<t_0,$ but does not exist at $t=t_0,$ then $T_{\max }=t_0.$
 \end{enumerate}
\subsection*{Inequalities algebrics}
\begin{lemma} \label{lemma2.2}(See \cite{ levine}\label{lemma10}) Let    $K$  be a positive function $K(t) \in C^2\left( [0, T), \mathbb{R}\right) $ which  satisfies
$$
K^{\prime \prime}(t) K(t)-\alpha\left(K^{\prime}(t)\right)^2 \geq 0,
\text{   for all } \alpha>1, \text{  and } 0<T \leq \infty.$$
If $K(0)>0, K^{\prime}(0)>0,$ then
\begin{align*}
T \leq \frac{K(0)}{(\alpha-1) K^{\prime}(0)}<\infty
\text {
and } K(t) \rightarrow \infty \text {
as } t \rightarrow T.
\end{align*}
\end{lemma}
\begin{lemma}( See \cite{Komornik})\label{lem1}

Let us consider the non-increasing function $R:[0,+\infty) \rightarrow[0,+\infty).$ Suppose there exist constants $\eta \geq 0$ and $C>0$ such that
$$
\int_t^{\infty} R^{1+\eta}(\tau) \mathrm{d} \tau \leq \frac{1}{C} R^\eta(0) R(t), \quad \text { for all } \quad t \geq 0.
$$

Then, $R$ satisfies the following decay estimates:
\begin{enumerate}
\item[1)]   \textbf{ Exponential decay ( $\eta=0$):} If $\eta=0,$ then $R(t) \leq R(0) e^{1-C t}$ for all $t \geq 0.$
\item[2)]
\textbf{ Polynomial decay ( $\eta>0$):} If $\eta>0,$ then $R(t) \leq R(0)\left(\frac{1+\eta}{1+\eta C t}\right)^{1 / \eta},$ for all $t \geq 0.$
\end{enumerate}
\end{lemma}
\end{definition}
\section{Potential well theory}\label{Potential well theory}
\begin{lemma} \label{lemma2}
For every $\eta>0,$ we get that
\begin{enumerate}
\item[(i)] $\log (t) \leq \frac{t^{\eta}}{\eta \exp (1)},$ for all $t \in(1, \infty).$
\item[(ii)] $t^\eta|\log (t)| \leq \frac{1}{\eta} \frac{1}{\exp (1)},$ for all $t \in(0,1].$
\end{enumerate}
\end{lemma}
\begin{proof}
 Setting the function
$$
g(t)=\left(\log (t)-\frac{t^\eta}{\eta \exp (1)}\right)
  \mathbf {1}_{(1,+\infty)},
$$

where $\mathbf {1}_{(1,+\infty)}  $ is an indicatrice function. It is evident that we prove that the function $g$ achieves its maximum at $t_*=$ $\exp \left(\frac{1}{\eta}\right).$ Thus, $g(t) \leq g\left(t_*\right),$ for all $t \geq 1.$ Now, we prove (ii). Considering that $t \rightarrow t^\eta|\log (t)|$ is a continuous on $(0,1]$ with $\lim _{t \rightarrow 0} t^\eta|\log (t)|=0,$ which achieves the maximum at $t_*=\exp \left(-\frac{1}{\eta}\right),$ we can easily conclude. This completes the  proof.
\end{proof}

Henceforward, $Y$ denotes the Cartesian product of two fractional Sobolev spcaes  $W^{s, p}\left(\mathcal{U}\right)$ and $W^{s, q}\left(\mathcal{U}\right),$ that is $Y= W^{s, p}\left(\mathcal{U}\right)\times  W^{s, q}\left(\mathcal{U}\right)$ equipped with the norm $\Vert \left(  \mathfrak{u}, \mathfrak{v}\right) \Vert_{Y}= \Vert \mathfrak{u} \Vert_{W^{s, p}\left(\mathcal{U}\right)} +
\Vert \mathfrak{v} \Vert_{W^{s, q}\left(\mathcal{U}\right)}.
 $
 \subsection{Potential well}
In this section, we introduce several key lemmas that play a crucial role in establishing our main results.
\begin{lemma} \label{lemma3}  Under assumptions  $\left(\mathcal{H}_1\right)- \left(\mathcal{H}_2\right),$ the following assertions hold:
\begin{enumerate}
\item[(1)] $\mu^\beta \mathcal{K}(t) \leq \mathcal{K}( \mu t),$ for all $0 \leq \mu \leq 1, z \geq 0.$
\item[(2)] $\mathcal{K}( \mu z) \leq \mu^\beta \mathcal{K}(z),$ for all $\mu \geq 1, z \geq 0.$
\item[(3)] $\frac{\mathcal{K}(\mu)}{\mu^\beta} \min \left\{\mu^\beta, z^\beta\right\} \leq \mathcal{K}(z) \leq \frac{\mathcal{K}(\mu)}{\mu^\beta} \max \left\{\mu^\beta, z^\beta\right\}.$
\item[(4)] $\mathcal{K}(z)>0,$ for all $z>0.$
\item[(5)] $\frac{1}{\beta+1} z \mathcal{K}(z) \leq \widehat{\mathcal{K}}( \mu z)\leq t \mathcal{K}(z),$ for all $z \geq 0.$
\item[(6)] $\widehat{\mathcal{K}}(\mu z) \geq \mu^{\beta+1} \widehat{\mathcal{K}}(z),$ for all $0 \leq \mu \leq 1, z \geq 0.$
\item[(7)] $\widehat{\mathcal{K}}(\mu z) \leq \mu^{\beta+1} \widehat{\mathcal{K}}(z),$ for all $\mu \geq 1, z \geq 0.$

\end{enumerate}
\end{lemma}
\begin{proof}
 Let $\mu \in(0,1)$ and $z>0.$ Applying the condition $\left(\mathcal{H}_2\right).$ Then, we have
$$
\frac{\mathcal{K}(z)}{z^\beta} \leq \frac{\mathcal{K}( \mu z)}{(\mu z)^\beta}.
$$

Clearly, we have (1). Identically, we establish (2). Combining (1) and (2) with condition  $(\left.\mathcal{H}_1\right),$ we obtain (3). From (3) and the fact that $\mathcal{K} \neq 0,$  we get (4). Since $\mathcal{K}$ is non-decreasing, it follows that

$$
\widehat{\mathcal{K}}(z)=\int_0^z \mathcal{K}(t) d t \leq \int_0^z \mathcal{K}(z) d t=z \mathcal{K}(z).
$$

Applying the condition $\left(\mathcal{H}_2\right),$ we get that
$$
\widehat{\mathcal{K}}(z)=\int_0^z \frac{\mathcal{K}(t)}{t^\beta} t^\beta d t \geq \int_0^z \frac{\mathcal{K}(z)}{z^\beta} t^\beta d t=\frac{1}{\beta+1} z \mathcal{K}(z).
$$

Let $\mu \in[0,1],$ we get that
$$
\widehat{\mathcal{K}}(\mu z)=\int_0^{\mu z} \mathcal{K}(t) d t=\int_0^z \mu \mathcal{K}(\mu t) d t \geq \mu^{\beta+1} \int_0^z \mathcal{K}(t) d t=\mu^{\beta+1} \widehat{\mathcal{K}}(z).
$$
\end{proof}
Let us define the  energy functional $\varphi: W^{s, p}(\mathcal{U}) \rightarrow \mathbb{R}$ and Nehari functional $\psi: W^{s, p}(\mathcal{U}) \rightarrow \mathbb{R}$ by:
$$
\varphi(\mathfrak{u}, \mathfrak{v}):= \frac{1}{p}\widehat{\mathcal{K}}\left([\mathfrak{u}]_p^s\right)+ \frac{1}{q}\widehat{\mathcal{K}}\left([\mathfrak{v}]_q^s\right)+\frac{1}{\sigma^2} \int_{\mathcal{U}}\vert  \mathfrak{u} \vert^\sigma\vert  \mathfrak{v} \vert^\sigma \mathrm{dx}-\frac{1}{\sigma} \int_{\mathcal{U}}\vert  \mathfrak{u} \vert^{ \sigma}\vert  \mathfrak{v} \vert^\sigma \log | \mathfrak{u} \mathfrak{v}| \mathrm{dx},
$$
and
$$
\psi(\mathfrak{u}, \mathfrak{v}):=\mathcal{K}\left([\mathfrak{u}]_p^s\right) \int_{\mathcal{U} \times \mathcal{U}} \frac{|\mathfrak{u}(\mathrm{x})-\mathfrak{u}(\mathrm{y})|^p}{|\mathrm{x}-\mathrm{y}|^{N+s p}} \mathrm{dxdy}+\mathcal{K}\left([\mathfrak{v}]_q^s\right) \int_{\mathcal{U} \times \mathcal{U}} \frac{|\mathfrak{v}(\mathrm{x})-\mathfrak{v}(\mathrm{y})|^q}{|\mathrm{x}-\mathrm{y}|^{N+s q}} d \mathrm{x} d \mathrm{y}-2 \int_{\mathcal{U}}\vert  \mathfrak{u} \vert^{\sigma+1}\vert  \mathfrak{v} \vert^{\sigma+1} \log | \mathfrak{u} \mathfrak{v}| \mathrm{dx}.
$$

\begin{lemma}\label{injection}  Let $(\mathfrak{u}, \mathfrak{v}) \in W^{s, p}(\mathcal{U}) \times W^{s, q}(\mathcal{U}).$ Then, we have the following inequality:
$$
\int_{\mathcal{U}}|\mathfrak{u}|^\sigma|\mathfrak{v}|^\sigma \log |\mathrm{uv}| d\mathrm{x} \leq \log [\mathfrak{u}]_{s, p} \int_{\mathcal{U}}|\mathfrak{u}(\mathrm{x})|^\sigma d \mathrm{x}+\log [\mathfrak{v}]_{s, q} \int_{\mathcal{U}}|\mathfrak{v}(\mathrm{x})|^\sigma d \mathrm{x}+
S\left( [\mathfrak{u}]_{s, p}^{p_s^*}+[\mathfrak{v}]_{s, p}^{p_s^*}+[\mathfrak{u}]_{s, q}^{q_s^*}+
[\mathfrak{v}]_{s, q}^{q_s^*}+[\mathfrak{u}]_{s, p}^\sigma+[\mathfrak{v}]_{s, q}^{\sigma }
\right)
.
$$
\end{lemma}

\begin{proof} Let $(\mathfrak{u}, \mathfrak{v}) \in W^{s, p}(\mathcal{U}) \times W^{s, q}(\mathcal{U}) \backslash\{(0,0)\}.$

\begin{align*}
\int_{\mathcal{U}}\vert  \mathfrak{u} \vert^{\sigma}\vert  \mathfrak{v} \vert^{ \sigma} \log | \mathfrak{u} \mathfrak{v}| \mathrm{dx} & \leq \int_{\mathcal{U}}\vert  \mathfrak{u} \vert^\sigma\vert  \mathfrak{v} \vert^\sigma|\log | u v|| \mathrm{dx} \\
& \leq \int_{\mathcal{U}} \frac{\vert  \mathfrak{u} \vert^{2 \sigma}+\vert  \mathfrak{v} \vert^{2 \sigma}}{2}(|\log | u||+|\log | v||) \mathrm{dx} \\
& \leq 2^{\sigma-2}\left(\int_{\mathcal{U}}\vert  \mathfrak{u} \vert^{2 \sigma}|\log | u||+\int_{\mathcal{U}}\vert  \mathfrak{u} \vert^{2 \sigma}|\log | v||+\left.\int_{\mathcal{U}}\vert  \mathfrak{v} \vert^{2 \sigma}|| \log \vert  \mathfrak{u} \vert\left|+\int_{\mathcal{U}}\right| v\right|^{2 \sigma}|\log | \mathfrak{v}||\right) \mathrm{dx} \\
& =2^{\sigma-2}\left(I_1+I_2+I_3+I_4\right).
\end{align*}

Now, we calculate the integrals $I_i,$ for all $i \in\{1,2,3,4\}.$ Set $\mathcal{U}_1:=\left\{\mathrm{x} \in \mathcal{U}:\left|\mathfrak{u}(\mathrm{x}) \right|\leq[\mathfrak{u}]_{s, p}\right\}$ and $\mathcal{U}_2:=\left\{\mathrm{x} \in \mathcal{U}:\left|\mathfrak{u}(\mathrm{x}) \right|>[\mathfrak{u}]_{s, p}\right\}.$ It is easy to see that $\mathcal{U}=\mathcal{U}_1 \cup \mathcal{U}_2$ and $\mathcal{U}_1 \cap \mathcal{U}_2=\emptyset.$ So, we get that
$$
\int_{\mathcal{U}}\vert  \mathfrak{u} \vert^\sigma\left|\log \frac{\mathfrak{u}(\mathrm{x})}{[\mathfrak{u}]_{s, p}}\right| \mathrm{dx}=\int_{\mathcal{U}_1}\vert  \mathfrak{u} \vert^\sigma\left|\log \frac{\mathfrak{u}(\mathrm{x})}{[\mathfrak{u}]_{s, p}}\right| \mathrm{dx}+\int_{\mathcal{U}_2}|\mathfrak{u}|^\sigma\left|\log \frac{\mathfrak{u}(\mathrm{x})}{[\mathfrak{u}]_{s, p}}\right| \mathrm{dx}.
$$

Using Lemma \ref{lemma2}, we deduce that
\begin{align}\label{equation2}
\begin{split}
\int_{\mathcal{U}_1}\vert  \mathfrak{u} \vert^\sigma\left|\log \frac{u(x)}{[\mathfrak{u}]_{s, p}}\right| \mathrm{dx} & =[\mathfrak{u}]_{s, p}^\sigma \int_{\mathcal{U}_1}\left|\frac{u(x)}{[\mathfrak{u}]_{s, p}}\right|^\sigma\left|\log \frac{u(x)}{[\mathfrak{u}]_{s, p}}\right| \mathrm{dx} \\
& \leq \frac{2^{\sigma-2}[\mathfrak{u}]_{s, p}^\sigma}{\sigma \exp (1)}
\end{split}
\end{align}

From Lemma \ref{lemma2} with $\eta=p_s^*-\sigma$ and fractional Sobolev embedding inequality, we obtain that

\begin{align}\label{equation3}
\begin{split}
\int_{U_2}|\mathfrak{u}|^\sigma\left|\log \frac{\mathfrak{u}(\mathrm{x})}{[\mathfrak{u}]_{s, p}}\right| \mathrm{dx} & \leq \frac{2^\sigma}{\exp \left(p_s^*-\sigma\right)[\mathfrak{u}]_{s, p}^{p_s^*-\sigma}} \int_{\mathcal{U}}|\mathfrak{u}|^{p_s^*} \mathrm{dx} \\
& \leq \frac{2^\sigma}{\exp \left(p_s^*-\sigma\right)[\mathfrak{u}]_{s, p}^{p_s^*-\sigma}} S_{p_x^*}^{p_x^*}[\mathfrak{u}]_{s, p}^\sigma,
\end{split}
\end{align}

with $S_{p_x^*}$ is the best constant of embedding from $W^{s, p}(\mathcal{U})$ to $L^{p_s^*}(\mathcal{U}).$ From \eqref{equation2} and \eqref{equation3} follow that
\begin{align}\label{equation4}
\begin{split}
\int_{\mathcal{U}}\vert  \mathfrak{u} \vert^\sigma\left|\log \frac{\mathfrak{u}(\mathrm{x})}{[\mathfrak{u}]_{s, p}}\right| \mathrm{dx} \leq\left(\frac{|\mathcal{U}|}{\sigma \exp (1)}+\frac{S_{p_x^*}^{p_x^*}}{\exp \left(p_s^*-\sigma\right)}\right)[\mathfrak{u}]_{s, p}^\sigma.
\end{split}
\end{align}

Using \eqref{equation4}, we deduce that

\begin{align*}
\int_{\mathcal{U}}\vert  \mathfrak{u} \vert^\sigma|\log | \mathfrak{u}(\mathrm{x})|| \mathrm{dx} & =\int_{\mathcal{U}}\vert  \mathfrak{u} \vert^\sigma\left|\log \frac{|\mathfrak{u}(\mathrm{x})|}{[\mathfrak{u}]_{s, p}}\right|+\log [\mathfrak{u}]_{s, p} \int_{\mathcal{U}}|\mathfrak{u}|^\sigma \mathrm{dx} \\
& \leq\left(\frac{|\mathcal{U}|}{\sigma \exp (1)}+\frac{S_{p_x^*}^{p_x^*}}{\exp \left(p_s^*-\sigma\right)}\right)[\mathfrak{u}]_{s, p}^\sigma+\log [\mathfrak{u}]_{s, p} \int_{\mathcal{U}}\vert  \mathfrak{u} \vert^\sigma \mathrm{dx}.
\end{align*}
Similarly, we prove that
\begin{align*}
\int_{\mathcal{U}}|\mathfrak{v}|^\sigma|\log | \mathfrak{v}(\mathrm{x})|| \mathrm{dx} \leq\left(\frac{|\mathcal{U}|}{\sigma \exp (1)}+\frac{S_{q_s^*}^{q_s^*}}{\exp \left(q_s^*-\sigma\right)}\right)[\mathfrak{v}]_{s, q}^{\sigma \sigma}+\log [\mathfrak{v}]_{s, q} \int_{\mathcal{U}}|\mathfrak{v}|^\sigma \mathrm{dx}
\end{align*}

Now, we estimate $I_3.$ It is easy to show that $\log \mathrm{x} \leq \mathrm{x}-1,$ for all $\mathrm{x}>0.$ remembering the numerical inequality
$$
a b \leq \frac{a^2+b^2}{2} \leq a^2+b^2, \text { for all }(a, b) \in \mathbb{R}^{+} \times \mathbb{R}^{+} .
$$

Using this inequality, we obtain that

\begin{align*}
\int_{\mathcal{U}}\left|\mathfrak{u}^\sigma\right| \log |\mathfrak{v}(\mathrm{x})| \mathrm{d} \mathrm{x} & \leq \int_{\mathcal{U}}|\mathfrak{u}(\mathrm{x})|^\sigma(\mathrm{|v}(\mathrm{x}) \mid-1) \mathrm{dx} \\
& \leq \int_{\mathcal{U}}|\mathfrak{u}(\mathrm{x})|^\sigma|\mathfrak{v}(\mathrm{x})|+\int_{\mathcal{U}}|\mathfrak{u}(\mathrm{x})|^\sigma \mathrm{dx} \\
& \leq \int_{\mathcal{U}}\left(|\mathfrak{u}(\mathrm{x})|^{2 \sigma}+|\mathfrak{v}(\mathrm{x})|^2+|\mathfrak{u}(\mathrm{x})|^\sigma\right) \mathrm{dx} \\
& \leq C(\sigma, \mathcal{U})\int_{\mathcal{U}}\left(|\mathfrak{u}(\mathrm{x})|^{p_x^*}+|\mathfrak{v}(\mathrm{x})|^{p_x^*}+|\mathfrak{u}(\mathrm{x})|^{p_x^*}\right) \mathrm{dx} \\
& \leq C(\sigma, \mathcal{U}, p, s) S_{p_s^*}^{p_x^*}\left([\mathfrak{u}]_{s, p}^{p_s^*}+[\mathfrak{v}]_{s, p}^{p_s^*}\right).
\end{align*}

Similarly, we prove that
\begin{align}\label{equation5}
\int_{\mathcal{U}}|\mathfrak{v}(\mathrm{x})|^{\sigma }|\log | \mathfrak{u}(\mathrm{x}) \mid \mathrm{dx} \leq C(\sigma, \mathcal{U}, q, s) S_{q_s^*}^{q_s^*}\left([\mathfrak{u}]_{s, q}^{q_s^*}+[\mathfrak{v}]_{s, q}^{q_s^*}\right) .
\end{align}

Using \eqref{equation4}- \eqref{equation5}, we deduce that
\begin{align*}
\int_{\mathcal{U}}|\mathfrak{u}|^\sigma|\mathfrak{v}|^\sigma \log |\mathrm{uv}| \mathrm{dx} \leq \log [\mathfrak{u}]_{s, p} \int_{\mathcal{U}}|\mathfrak{u}(\mathrm{x})|^\sigma \mathrm{dx}+\log [\mathfrak{v}]_{s, q} \int_{\mathcal{U}}|\mathfrak{v}(\mathrm{x})|^\sigma \mathrm{dx}+S\left([\mathfrak{u}]_{s, p}^{p_s^*}+[ \mathfrak{v}]_{s, p}^{p_s^*}+[\mathfrak{u}]_{s, q}^{q_s^*}+[\mathfrak{v}]_{s, q}^{q_s^*}+[\mathfrak{u}]_{s, p}^\sigma+[\mathfrak{u}]_{s, q}^\sigma\right),
\end{align*}
where $S:=\max \left\{C(\sigma, \mathcal{U}, q, s) S_{q_s^*}^{q_x^*}, C(\sigma, \mathcal{U}, p, s) S_{p_s^*}^{p_t^*}, \frac{|\mathcal{U}|}{\sigma \exp (1)}+\frac{s_{\psi_i^*}^{q_s^*}}{\exp \left(q_s^*-\sigma\right)}, \frac{|\mathcal{U}|}{\sigma \exp (1)}+\frac{s_{p_t^*}^{p_t^*}}{\exp \left(p_s^*-\sigma\right)}\right\}.$
\end{proof}
\begin{lemma}  \label{lemma4} Let $\left(\mathcal{H}_1\right)-\left(\mathcal{H}_2\right)$ holds and $(u, v) \in W^{s, p}(\mathcal{U}) \times W^{s . q}(\mathcal{U}) \backslash\{(0,0)\}.$ Then following results hold.
\begin{enumerate}

\item[(1)] There exists a unique constant $\varepsilon_*>0$ such that $\psi(\varepsilon \mathfrak{u}, \varepsilon \mathfrak{v})>0$ for all $0<\varepsilon<\varepsilon_*, \psi\left(\varepsilon_* \mathfrak{v}, \varepsilon_* \mathfrak{v}\right)=0,$ and $\psi(\varepsilon \mathfrak{u}, \varepsilon \mathfrak{v})<0$ for all $\varepsilon>\varepsilon_*.$
\item[(2)] The function $\varepsilon \mapsto \varphi(\varepsilon u, \varepsilon v)$ strictly decreasing on $\left(\varepsilon_*, \infty\right),$ and   it is strictly increasing on $\left(0, \varepsilon_*\right).$ Furthermore, we have $\lim _{\varepsilon \rightarrow 0^*} \varphi(\varepsilon \mathfrak{u}, \varepsilon \mathfrak{v})=0,$ and $\lim _{\varepsilon \rightarrow+\infty} \varphi(\varepsilon \mathfrak{u}, \varepsilon \mathfrak{v})=-\infty.$
\end{enumerate}
\end{lemma}
\begin{proof}
 For (1). From (2) in Lemma \ref{lemma3}, we have, for all $\varepsilon>1,$
\begin{align*}
& \psi(\varepsilon \mathfrak{u}, \varepsilon \mathfrak{v})=\mathcal{K}\left([\varepsilon \mathfrak{u}]_p^s\right) \int_{\mathcal{U} \times \mathcal{U}} \frac{\varepsilon^p|\mathfrak{u}(\mathrm{x})-\mathfrak{u}(\mathrm{y})|^p}{|\mathrm{x}-\mathrm{y}|^{N+s p}} \mathrm{dxdy}+\mathcal{K}\left([\varepsilon \mathfrak{v}]_q^s\right) \int_{\mathcal{U} \times \mathcal{U}} \frac{\varepsilon^q|\mathfrak{v}(\mathrm{x})-\mathfrak{v}(\mathrm{y})|^q}{|\mathrm{x}-\mathrm{y}|^{N+s q}} \mathrm{dxdy}-\varepsilon^{2 \sigma+2} \int_{\mathcal{U}}|\mathrm{uv}|^{\sigma +1} \log \left|\varepsilon^2 \mathrm{uv}\right| \mathrm{dx} \\
& \leq \varepsilon^{P(\beta+1)} \mathcal{K}\left([\mathfrak{u}]_p^s\right) \int_{\mathcal{U} \times U} \frac{|\mathfrak{u}(\mathrm{x})-\mathfrak{u}(\mathrm{y})|^p}{|\mathrm{x}-\mathrm{y}|^{N+s p}} \mathrm{dxdy}+\varepsilon^{q(\beta+1)} \mathcal{K}\left([\mathfrak{v}]_q^s\right) \int_{U_{\times} U} \frac{|\mathfrak{v}(\mathrm{x})-\mathfrak{v}(\mathrm{y})|^q}{|\mathrm{x}-\mathrm{y}|^{N+s q}} \mathrm{dxdy}-\varepsilon^{2 \sigma+2} \int_{\mathcal{U}}|u v|^{\sigma+1} \log \left|\varepsilon^2 \mathrm{uv}\right| \mathrm{dx} \\
& \leq \varepsilon^\eta\left[\mathcal{K}\left([\mathfrak{u}]_p^s\right) \int_{\mathcal{U} \times U} \frac{|\mathfrak{u}(\mathrm{x})-\mathfrak{u}(\mathrm{y})|^p}{|\mathrm{x}-\mathrm{y}|^{N+s p}} \mathrm{dxdy}+\mathcal{K}\left([\mathfrak{v}]_q^s\right) \int_{\mathcal{U} \times \mathcal{U}} \frac{|\mathfrak{v}(\mathrm{x})-\mathfrak{v}(\mathrm{y})|^q}{|\mathrm{x}-\mathrm{y}|^{N+s q}} \mathrm{dxdy}\right]-\varepsilon^{2 \sigma+2} \int_{\mathcal{U}}|\mathrm{uv}|^{\sigma+1} \log \left|\varepsilon^2 u v\right| \mathrm{dx},
\end{align*}

where $\eta:= \max (p(\beta + 1), q(\beta + 1)).$ Due to $\eta < 2(\sigma + 1), $ we have that
\begin{align}\label{equation8}
\lim_{\varepsilon\to \infty}\psi( \varepsilon \mathfrak{u}, \varepsilon \mathfrak{v})=-\infty.
\end{align}
Using (i) in Lemma \ref{lemma3}, for any  $\varepsilon \in (0, 1),$ we obtain that
\begin{align*}
& \psi(\varepsilon \mathfrak{u}, \varepsilon \mathfrak{v})=\mathcal{K}\left([\varepsilon \mathfrak{u}]_p^s\right) \int_{\mathcal{U} \times \mathcal{U}} \frac{\varepsilon^P|\mathfrak{u}(\mathrm{x})-\mathfrak{u}(\mathrm{y})|^p}{|\mathrm{x}-\mathrm{y}|^{N+s p}} \mathrm{dxdy}+\mathcal{K}\left([\varepsilon \mathfrak{v}]_q^s\right) \int_{\mathcal{U} \times \mathcal{U}} \frac{\varepsilon^q|\mathfrak{v}(\mathrm{x})-\mathfrak{v}(\mathrm{y})|^q}{|\mathrm{x}-\mathrm{y}|^{N+s q}} \mathrm{dxdy}-\varepsilon^{2 \sigma+2} \int_{\mathcal{U}}|\mathrm{uv}|^{\sigma +1} \log \left|\varepsilon^2 \mathrm{uv}\right| \mathrm{dx} \\
&\geq \varepsilon^{p(\beta+1)} \mathcal{K}\left([ \mathfrak{u}]_p^s\right) \int_{\mathcal{U} \times \mathcal{U}} \frac{|\mathfrak{u}(\mathrm{x})-\mathfrak{u}(\mathrm{y})|^p}{|\mathrm{x}-\mathrm{y}|^{N+s p}} \mathrm{dxdy}+ \varepsilon^{q(\beta+1)}
\mathcal{K}\left([ \mathfrak{v}]_q^s\right) \int_{\mathcal{U} \times \mathcal{U}} \frac{|\mathfrak{v}(\mathrm{x})-\mathfrak{v}(\mathrm{y})|^q}{|\mathrm{x}-\mathrm{y}|^{N+s q}} \mathrm{dxdy}-\varepsilon^{2 \sigma+2} \int_{\mathcal{U}}|\mathrm{uv}|^{\sigma +1} \log \left|\varepsilon^2 \mathrm{uv}\right| \mathrm{dx} \\
& \geq \varepsilon^{\eta} \left[ \mathcal{K}\left([ \mathfrak{u}]_p^s\right) \int_{\mathcal{U} \times \mathcal{U}} \frac{|\mathfrak{u}(\mathrm{x})-\mathfrak{u}(\mathrm{y})|^p}{|\mathrm{x}-\mathrm{y}|^{N+s p}} \mathrm{dxdy}+
\mathcal{K}\left([ \mathfrak{v}]_q^s\right) \int_{\mathcal{U} \times \mathcal{U}} \frac{|\mathfrak{v}(\mathrm{x})-\mathfrak{v}(\mathrm{y})|^q}{|\mathrm{x}-\mathrm{y}|^{N+s q}} \mathrm{dxdy}\right] - \varepsilon^{2 \sigma+2} \int_{\mathcal{U}}|\mathrm{uv}|^{\sigma +1} \log \left|\varepsilon^2 \mathrm{uv}\right| \mathrm{dx},
\end{align*}
where $\rho:=\min \left\lbrace q(\beta+1), p(\beta+1)\right\rbrace, $ which yields that $\psi(\varepsilon\mathfrak{u}, \varepsilon\mathfrak{v})>0,$ for all $0<\varepsilon<1.$ Thanks to the  intermediate value
theorem and  \eqref{equation8}, there exists a $\varepsilon_{*}>0$ such that $\psi\left(\varepsilon_{*} \mathfrak{u}, \varepsilon_{*} \mathfrak{v}\right)=0.$ For any $\varepsilon>\varepsilon_{*},$ we get that $\frac{\varepsilon}{\varepsilon_{*}}>1.$ Using (2) in Lemma  \ref{lemma3}, we obtain that

$$
\begin{aligned}
& \psi(\varepsilon \mathfrak{u}, \varepsilon \mathfrak{v})=\mathcal{K}\left([\varepsilon \mathfrak{u}]_{p}^{s}\right) \int_{\mathcal{U} \times \mathcal{U}} \frac{\varepsilon^{p}|\mathfrak{u}(\mathrm{x})-\mathfrak{u}(\mathrm{y})|^{p}}{|\mathrm{x}-\mathrm{y}|^{N+s p}} \mathrm{dxdy}+\mathcal{K}\left([\varepsilon \mathfrak{v}]_{q}^{s}\right) \int_{\mathcal{U} \times \mathcal{U}} \frac{\varepsilon^{q}|\mathfrak{v}(\mathrm{x})-\mathfrak{v}(\mathrm{y})|^{q}}{|\mathrm{x}-\mathrm{y}|^{N+s q}} d \mathrm{x} d \mathrm{y}-\varepsilon^{2 \sigma+2} \int_{\mathcal{U}}|\mathfrak{u}\mathfrak{v}|^{\sigma+1} \log \left|\varepsilon^{2} \mathfrak{u}\mathfrak{v}\right| \mathrm{dx} \\
& =\mathcal{K}\left(\frac{\left(\frac{\varepsilon}{\varepsilon_{*}}\right)^{p}}{p} \int_{\mathcal{U} \times \mathcal{U}} \frac{\varepsilon_{*}^{p}|\mathfrak{u}(\mathrm{x})-\mathfrak{u}(\mathrm{y})|^{p}}{|\mathrm{x}-\mathrm{y}|^{N+s p}} \mathrm{dxdy}\right) \int_{\mathcal{U} \times \mathcal{U}}\left(\frac{\varepsilon}{\varepsilon_{*}}\right)^{p} \frac{\varepsilon_{*}^{p}|\mathfrak{u}(\mathrm{x})-\mathfrak{u}(\mathrm{y})|^{p}}{|\mathrm{x}-\mathrm{y}|^{N+s p}} \mathrm{dxdy} \\
& +\mathcal{K}\left(\frac{\left(\frac{\varepsilon}{\varepsilon_{*}}\right)^{q}}{q} \int_{\mathcal{U} \times \mathcal{U}} \frac{\varepsilon_{*}^{q}|\mathfrak{v}(\mathrm{x})-\mathfrak{v}(\mathrm{y})|^{q}}{|\mathrm{x}-\mathrm{y}|^{N+s q}} \mathrm{dxdy}\right) \int_{\mathcal{U} \times \mathcal{U}}\left(\frac{\varepsilon}{\varepsilon_{*}}\right)^{q} \frac{\varepsilon_{*}^{q}|\mathfrak{v}(\mathrm{x})-\mathfrak{v}(\mathrm{y})|^{q}}{|\mathrm{x}-\mathrm{y}|^{N+s q}} \mathrm{dxdy} \\
& -\left(\frac{\varepsilon}{\varepsilon_{*}}\right)^{2 \sigma+2} \int_{\mathcal{U}}\left|\varepsilon_{*} \mathfrak{u}\right|^{\sigma+1}\left|\varepsilon_{*} \mathfrak{v}\right|^{\sigma+1} \log \left|\varepsilon_{*}^{2} \mathfrak{u}\mathfrak{v}\left(\frac{\varepsilon}{\varepsilon_{*}}\right)^{2}\right| \mathrm{dx} \\
& \leq\left(\frac{\varepsilon}{\varepsilon_{*}}\right)^{p(\beta+1)} \mathcal{K}\left(\left[\varepsilon_{\varepsilon} \mathfrak{u}\right]_{p}^{s}\right) \int_{\mathcal{U} \times \mathcal{U}} \frac{\varepsilon_{*}^{p}|\mathfrak{u}(\mathrm{x})-\mathfrak{u}(\mathrm{y})|^{p}}{|\mathrm{x}-\mathrm{y}|^{N+s p}} \mathrm{dxdy}+\left(\frac{\varepsilon}{\varepsilon_{*}}\right)^{q(\beta+1)} \mathcal{K}\left(\left[\varepsilon_{*} \mathfrak{v}\right]_{q}^{s}\right) \int_{\mathcal{U} \times \mathcal{U}} \frac{\varepsilon_{*}^{q}|\mathfrak{v}(\mathrm{x})-\mathfrak{v}(\mathrm{y})|^{q}}{|\mathrm{x}-\mathrm{y}|^{N+s q}} d \mathrm{x} d \mathrm{y} \\
& -\left(\frac{\varepsilon}{\varepsilon_{*}}\right)^{2 \sigma+2} \int_{\mathcal{U}}\left|\varepsilon_{*} \mathfrak{u}\right|^{\sigma+1}\left|\varepsilon_{*} \mathfrak{v}\right|^{\sigma+1} \log \left|\varepsilon_{*}^{2} \mathfrak{u}\mathfrak{v}\right| \mathrm{dx} \\
& \leq\left(\frac{\varepsilon}{\varepsilon_{*}}\right)^{\delta}\left[\mathcal{K}\left(\left[\varepsilon_{\varepsilon} \mathfrak{u}\right]_{p}^{s}\right) \int_{\mathcal{U} \times \mathcal{U}} \frac{\varepsilon_{*}^{p}|\mathfrak{u}(\mathrm{x})-\mathfrak{u}(\mathrm{y})|^{p}}{|\mathrm{x}-\mathrm{y}|^{N+s p}} \mathrm{dxdy}+\mathcal{K}\left(\left[\varepsilon_{*} \mathfrak{v}\right]_{q}^{s}\right) \int_{\mathcal{U} \times \mathcal{U}} \frac{\varepsilon_{*}^{q}|\mathfrak{u}(\mathrm{x})-\mathfrak{v}(\mathrm{y})|^{q}}{|\mathrm{x}-\mathrm{y}|^{N+s q}} d \mathrm{x} d \mathrm{y}\right] \\
& -\left(\frac{\varepsilon}{\varepsilon_{*}}\right)^{2 \sigma+2} \int_{\mathcal{U}}\left|\varepsilon_{*} \mathfrak{u}\right|^{\sigma+1}\left|\varepsilon_{*} \mathfrak{v}\right|^{\sigma+1} \log \left|\varepsilon_{*}^{2} \mathfrak{u}\mathfrak{v}\right| \mathrm{dx} \\
& \leq\left[\left(\frac{\varepsilon}{\varepsilon_{*}}\right)^{\delta}-\left(\frac{\varepsilon}{\varepsilon_{*}}\right)^{2 \sigma+2}\right] \int_{\mathcal{U}}\left|\varepsilon_{*} \mathfrak{u}\right|^{\sigma+1}\left|\varepsilon_{*} \mathfrak{v}\right|^{\sigma+1} \log \left|\varepsilon_{*}^{2} \mathfrak{u} \mathfrak{v}\right| \mathrm{dx}+\left(\frac{\varepsilon}{\varepsilon_{*}}\right)^{\delta} \psi\left(\varepsilon_{*} \mathfrak{u}, \varepsilon_{*} \mathfrak{v}\right),
\end{aligned}
$$

where $\delta:=\max \{p(\beta+1), q(\beta+1)\}.$ Since $\frac{\varepsilon}{\varepsilon_{*}}>1$ and $\delta<2(\sigma+1),$ we get that $\psi(\varepsilon \mathfrak{u}, \varepsilon \mathfrak{v})<0.$ Similarly, for any
$0<\varepsilon<\varepsilon_{*}.$ Using (i) in Lemma \ref{lemma3}, we get that

$$
\begin{aligned}
& \psi(\varepsilon \mathfrak{u}, \varepsilon \mathfrak{v})=\mathcal{K}\left([\varepsilon \mathfrak{u}]_{p}^{s}\right) \int_{\mathcal{U} \times \mathcal{U}} \frac{\varepsilon^{p}|\mathfrak{u}(\mathrm{x})-\mathfrak{u}(\mathrm{y})|^{p}}{|\mathrm{x}-\mathrm{y}|^{N+s p}} \mathrm{dxdy}+\mathcal{K}\left([\varepsilon \mathfrak{v}]_{q}^{s}\right) \int_{\mathcal{U} \times \mathcal{U}} \frac{\varepsilon^{q}|\mathfrak{v}(\mathrm{x})-\mathfrak{v}(\mathrm{y})|^{q}}{|\mathrm{x}-\mathrm{y}|^{N+s q}} d \mathrm{x} d \mathrm{y}-\varepsilon^{2 \sigma+2} \int_{\mathcal{U}}|\mathfrak{u}\mathfrak{v}|^{\sigma+1} \log \left|\varepsilon^{2} \mathfrak{u}\mathfrak{v}\right| \mathrm{dx} \\
& =\mathcal{K}\left(\frac{\left(\frac{\varepsilon}{\varepsilon_{*}}\right)^{p}}{p} \int_{\mathcal{U} \times \mathcal{U}} \frac{\varepsilon_{*}^{p}|\mathfrak{u}(\mathrm{x})-\mathfrak{u}(\mathrm{y})|^{p}}{|\mathrm{x}-\mathrm{y}|^{N+s p}} \mathrm{dxdy}\right) \int_{\mathcal{U} \times \mathcal{U}}\left(\frac{\varepsilon}{\varepsilon_{*}}\right)^{p} \frac{\varepsilon_{*}^{p}|\mathfrak{u}(\mathrm{x})-\mathfrak{u}(\mathrm{y})|^{p}}{|\mathrm{x}-\mathrm{y}|^{N+s p}} \mathrm{dxdy} \\
& +\mathcal{K}\left(\frac{\left(\frac{\varepsilon}{\varepsilon_{*}}\right)^{q}}{q} \int_{\mathcal{U} \times \mathcal{U}} \frac{\varepsilon_{*}^{q}|\mathfrak{v}(\mathrm{x})-\mathfrak{v}(\mathrm{y})|^{q}}{|\mathrm{x}-\mathrm{y}|^{N+s q}} \mathrm{dxdy}\right) \int_{\mathcal{U} \times \mathcal{U}}\left(\frac{\varepsilon}{\varepsilon_{*}}\right)^{q} \frac{\varepsilon_{*}^{q}|\mathfrak{v}(\mathrm{x})-\mathfrak{v}(\mathrm{y})|^{q}}{|\mathrm{x}-\mathrm{y}|^{N+s q}} \mathrm{dxdy} \\
& -\left(\frac{\varepsilon}{\varepsilon_{*}}\right)^{2 \sigma+2} \int_{\mathcal{U}}\left|\varepsilon_{*} \mathfrak{u}\right|^{\sigma+1}\left|\varepsilon_{*} \mathfrak{v}\right|^{\sigma+1} \log \left|\varepsilon_{*}^{2} \mathfrak{v}\left(\frac{\varepsilon}{\varepsilon_{*}}\right)^{2}\right| \mathrm{dx} \\
& \geq\left(\frac{\varepsilon}{\varepsilon_{*}}\right)^{p(\beta+1)} \mathcal{K}\left(\left[\varepsilon_{\varepsilon} \mathfrak{u}\right]_{p}^{s}\right) \int_{\mathcal{U} \times \mathcal{U}} \frac{\varepsilon_{*}^{p}|\mathfrak{u}(\mathrm{x})-\mathfrak{u}(\mathrm{y})|^{p}}{|\mathrm{x}-\mathrm{y}|^{N+s p}} \mathrm{dxdy}+\left(\frac{\varepsilon}{\varepsilon_{*}}\right)^{q(\beta+1)} \mathcal{K}\left(\left[\varepsilon_{*} \mathfrak{v}\right]_{q}^{s}\right) \int_{\mathcal{U} \times \mathcal{U}} \frac{\varepsilon_{*}^{q}|\mathfrak{v}(\mathrm{x})-\mathfrak{v}(\mathrm{y})|^{q}}{|\mathrm{x}-\mathrm{y}|^{N+s q}} d \mathrm{x} d \mathrm{y} \\
& -\left(\frac{\varepsilon}{\varepsilon_{*}}\right)^{2 \sigma+2} \int_{\mathcal{U}}\left|\varepsilon_{*} \mathfrak{u}\right|^{\sigma+1}\left|\varepsilon_{*} \mathfrak{v}\right|^{\sigma+1} \log \left|\varepsilon_{*}^{2} \mathfrak{u}\right| \mathrm{dx} \\
& \geq\left(\frac{\varepsilon}{\varepsilon_{*}}\right)^{\mu}\left[\mathcal{K}\left(\left[\varepsilon_{\varepsilon} \mathfrak{u}\right]_{p}^{s}\right) \int_{\mathcal{U} \times \mathcal{U}} \frac{\varepsilon_{*}^{p}|\mathfrak{u}(\mathrm{x})-\mathfrak{u}(\mathrm{y})|^{p}}{|\mathrm{x}-\mathrm{y}|^{N+s p}} \mathrm{dxdy}+\mathcal{K}\left(\left[\varepsilon_{*} \mathfrak{v}\right]_{q}^{s}\right) \int_{\mathcal{U} \times \mathcal{U}} \frac{\varepsilon_{*}^{q}|\mathfrak{v}(\mathrm{x})-\mathfrak{v}(\mathrm{y})|^{q}}{|\mathrm{x}-\mathrm{y}|^{N+s q}} d \mathrm{x} d \mathrm{y}\right] \\
& -\left(\frac{\varepsilon}{\varepsilon_{*}}\right)^{2 \sigma+2} \int_{\mathcal{U}}\left|\varepsilon_{*} \mathfrak{u}\right|^{\sigma+1}\left|\varepsilon_{*} \mathfrak{v}\right|^{\sigma+1} \log \left|\varepsilon_{*}^{2} \mathfrak{u}\mathfrak{v}\right| \mathrm{dx} \\
& \geq\left[\left(\frac{\varepsilon}{\varepsilon_{*}}\right)^{\mu}-\left(\frac{\varepsilon}{\varepsilon_{*}}\right)^{2 \sigma+2}\right] \int_{\mathcal{U}}\left|\varepsilon_{*} \mathfrak{u}\right|^{\sigma+1}\left|\varepsilon_{*} \mathfrak{v}\right|^{\sigma+1} \log \left|\varepsilon_{*}^{2} \mathfrak{v}\right| \mathrm{dx}+\left(\frac{\varepsilon}{\varepsilon_{*}}\right)^{\mu} \psi\left(\varepsilon_{*} \mathfrak{u}, \varepsilon_{*} \mathfrak{v}\right),
\end{aligned}
$$

with $\mu:=\min \{p(\beta+1), q(\beta+1)\}.$ Since $\frac{\varepsilon}{\varepsilon_{t}}<1$ and $\mu<2(\sigma+1),$ we get that $\psi(\varepsilon \mathfrak{u}, \varepsilon \mathfrak{v})>0.$ Now, we prove (2). According to the definition of $\varphi,$ we get that

$$
\varphi(\varepsilon \mathfrak{u}, \varepsilon \mathfrak{v})=\frac{1}{p}\widehat{\mathcal{K}}\left([\varepsilon \mathfrak{u}]_{p}^{s}\right)+\frac{1}{q}\widehat{\mathcal{K}}\left([\varepsilon \mathfrak{v}]_{q}^{s}\right)+\frac{1}{\sigma^{2}} \int_{\mathcal{U}}|\varepsilon \mathfrak{u}|^{\sigma}|\varepsilon \mathfrak{v}|^{\sigma} \mathrm{dx}-\frac{1}{\sigma} \int_{\mathcal{U}}|\varepsilon \mathfrak{u}|^{\sigma}|\varepsilon \mathfrak{v}|^{\sigma} \log |\varepsilon \mathfrak{u} \varepsilon \mathfrak{v}| \mathrm{dx}.
$$

Then one can verify that

\begin{align}\label{equation9}
\frac{d}{d \varepsilon} \varphi(\varepsilon \mathfrak{u}, \varepsilon \mathfrak{v})=\frac{1}{\varepsilon} \psi(\varepsilon \mathfrak{u}, \varepsilon \mathfrak{v}).
\end{align}

Combining \eqref{equation9} with (1) in Lemma \ref{lemma4}, we deduce that the map $\varepsilon \mapsto \varphi(\varepsilon \mathfrak{u}, \varepsilon \mathfrak{v})$ attains its maximum at $\varepsilon=\varepsilon_{*},$ is strictly decreasing on $\left(\varepsilon_{*}, \infty\right),$ and is strictly increasing on $\left(0, \varepsilon_{*}\right).$ By using a similar discussion as (1) in Lemma \ref{lemma4}, one can verify that $\lim _{\varepsilon \rightarrow 0^{+}} \varphi(\varepsilon \mathfrak{u}, \varepsilon \mathfrak{v})=0,$ and $\lim _{\varepsilon \rightarrow+\infty} \varphi(\varepsilon \mathfrak{u}, \varepsilon \mathfrak{v})=-\infty.$
\end{proof}

Considering the following elliptic system:

\begin{equation}\label{problem2}
 \left\{\begin{array}{lll}
\displaystyle
\mathcal{K}\left([\mathfrak{u}]_{p}^{s}\right) \mathscr{L}_{p}^{s} \mathfrak{u}=|\mathfrak{v}|^{\sigma}|\mathfrak{u}|^{\sigma-2} \mathfrak{u} \log |\mathfrak{v} \mathfrak{u}|, & \text { in } & \mathcal{U}, \\
\mathcal{K}\left([\mathfrak{v}]_{q}^{s}\right) \mathscr{L}_{q}^{s} \mathfrak{v}=|\mathfrak{u}|^{\sigma}|\mathfrak{v}|^{\sigma-2} \mathfrak{v} \log |\mathfrak{u} \mathfrak{v}|, & \text { in } & \mathcal{U}, \\
\mathfrak{u}(\mathrm{x})=\mathfrak{u}(\mathrm{x})=0, & \text { in } & \partial \mathcal{U},
\end{array}%
\right.
\end{equation}

By Lemma \ref{lemma4}, the problem \eqref{problem2} has a weak solution in $W^{s, p}(\mathcal{U}) \times W^{s, q}(\mathcal{U}).$

\begin{corollary}  Let us define the Nehari manifold set

$$
\mathcal{N}:=\left\lbrace (\mathfrak{u}, \mathfrak{v}) \in W^{s, p}(\mathcal{U}) \times W^{s, q}(\mathcal{U}): \psi(\mathfrak{u}, \mathfrak{v})=0\right\rbrace  \neq \emptyset, \text { and } d:=\inf _{(\mathfrak{u}, \mathfrak{v}) \in \mathcal{N}} \varphi(\mathfrak{u}, \mathfrak{v}).
$$

The number $d$ is called the depth of potential well.
\end{corollary}

In what follows, we prove an essential lemma for proving the strict positivity of depth of potential well.

\begin{lemma}\label{lemma5}

If Kirchhoff's function satisfies the conditions $\left(\mathcal{H}_{1}\right)-\left(\mathcal{H}_{2}\right)$ and the numbers $q, p, \sigma$ satisfy the condition \eqref{equation1}. Then, we obtain the following results:
\begin{enumerate}

\item[1)] $\varphi(\mathfrak{u}, \mathfrak{v})-\frac{1}{\sigma} \psi(\mathfrak{u}, \mathfrak{v}) \leq\left(\frac{1}{p}-\frac{1}{\sigma}\right)\left[\mathcal{K}\left([\mathfrak{u}]_{p}^{s}\right) \int_{\mathcal{U} \times \mathcal{U}} \frac{|\mathfrak{u}(\mathrm{x})-\mathfrak{u}(\mathrm{y})|^{p}}{|\mathrm{x}-\mathrm{y}|^{\hat{N} s p}} d \mathrm{x} d \mathrm{y}+\mathcal{K}\left([\mathfrak{v}]_{q}^{s}\right) \int_{\mathcal{U} \times \mathcal{U}} \frac{|\mathfrak{v}(\mathrm{x})-\mathfrak{v}(\mathrm{y})|^{q}}{|\mathrm{x}-\mathrm{y}|{ }^{\mid+s q}} d \mathrm{x} d \mathrm{y}\right]+\frac{1}{\sigma} \int_{\mathcal{U}}|\mathfrak{u}|^{\sigma}|\mathfrak{v}|^{\sigma} d \mathrm{x},$
\item[2)]  $\varphi(\mathfrak{u}, \mathfrak{v})-\frac{1}{\sigma} \psi(\mathfrak{u}, \mathfrak{v}) \geq\left(\frac{1}{q(\beta+1)}-\frac{1}{\sigma}\right)\left[\mathcal{K}\left([\mathfrak{u}]_{p}^{s}\right) \int_{\mathcal{U} \times \mathcal{U}} \frac{|\mathfrak{u}(\mathrm{x})-\mathfrak{u}(\mathrm{y})|^{p}}{|\mathrm{x}-\mathrm{y}|^{N+s p}} d \mathrm{x} d \mathrm{y}+\mathcal{K}\left([\mathfrak{v}]_{q}^{s}\right) \int_{\mathcal{U} \times \mathcal{U}} \frac{|\mathfrak{v}(\mathrm{x})-\mathfrak{v}(\mathrm{y})|^{q}}{|\mathrm{x}-\mathrm{y}|^{N+s q}} d \mathrm{x} d \mathrm{y}\right],$

\item[3)]  If $\psi(\mathfrak{u}, \mathfrak{v})<0 \Rightarrow \varphi(\mathfrak{u}, \mathfrak{v})-\frac{1}{\sigma} \psi(\mathfrak{u}, \mathfrak{v})>d_{*},$

\item[4)]  If $\psi(\mathfrak{u}, \mathfrak{v}) \geq 0 \Rightarrow \varepsilon_{*}^{-q(\beta+1)}-C^{\prime}\left(\varepsilon_{*}, \sigma, d, \mathfrak{u}, \mathfrak{v}\right) \leq \varphi(\mathfrak{u}, \mathfrak{v})-\frac{1}{\sigma} \psi(\mathfrak{u}, \mathfrak{v})$ and

$$
\psi(\mathfrak{u}, \mathfrak{v}) \geq\left(1-\varepsilon_{*}^{q(\beta+1)-2(\sigma+1)}\right)\left[\mathcal{K}\left([\mathfrak{u}]_{p}^{s}\right) \int_{\mathcal{U} \times \mathcal{U}} \frac{|\mathfrak{u}(\mathrm{x})-\mathfrak{u}(\mathrm{y})|^{p}}{|\mathrm{x}-\mathrm{y}|^{N+s p}} d \mathrm{x} d \mathrm{y}+\mathcal{K}\left(   [\mathfrak{v}]_{q}^{s}\right) \int_{\mathcal{U} \times \mathcal{U}} \frac{|\mathfrak{v}(\mathrm{x})-\mathfrak{v}(\mathrm{y})|^{q}}{|\mathrm{x}-\mathrm{y}|^{N+s q}} d \mathrm{x} d \mathrm{y}\right].
$$
\end{enumerate}
\end{lemma}
\begin{proof}

 According to $\varphi$ and $\psi$ 's definition, we get that
\begin{align}\label{equation22}
\begin{split}
\varphi(\mathfrak{u}, \mathfrak{v}) & = \frac{1}{p}\widehat{\mathcal{K}}\left([\mathfrak{u}]_{p}^{s}\right)+ \frac{1}{q}\widehat{\mathcal{K}}\left([\mathfrak{v}]_{q}^{s}\right)+\frac{1}{\sigma^{2}} \int_{\mathcal{U}}|\mathfrak{u}|^{\sigma}|\mathfrak{v}|^{\sigma} \mathrm{dx}-\frac{1}{\sigma} \int_{\mathcal{U}}|\mathfrak{u}|^{\sigma}|\mathfrak{v}|^{\sigma} \log |\mathfrak{u}| \mathrm{dx} \\
& \leq \mathcal{K}\left([\mathfrak{u}]_{p}^{s}\right) \int_{\mathcal{U} \times \mathcal{U}} \frac{1}{p} \frac{|\mathfrak{u}(\mathrm{x})-\mathfrak{u}(\mathrm{y})|^{p}}{|\mathrm{x}-\mathrm{y}|^{N+s p}} \mathrm{dxdy}+\mathcal{K}\left([\mathfrak{v}]_{q}^{s}\right) \int_{\mathcal{U} \times \mathcal{U}} \frac{1}{q} \frac{|\mathfrak{v}(\mathrm{x})-\mathfrak{v}(\mathrm{y})|^{q}}{|\mathrm{x}-\mathrm{y}|^{N+s q}} \mathrm{dxdy}+\frac{1}{\sigma} \int_{\mathcal{U}}|\mathfrak{u}|^{\sigma}|\mathfrak{v}|^{\sigma} \mathrm{dx} \\
& -\frac{1}{\sigma} \int_{\mathcal{U}}|\mathfrak{u}|^{\sigma}|\mathfrak{v}|^{\sigma} \log |\mathfrak{u} \mathfrak{v}| \mathrm{dx} \\
& \leq \frac{1}{p}\left[\mathcal{K}\left([\mathfrak{u}]_{p}^{s}\right) \int_{\mathcal{U} \times \mathcal{U}} \frac{|\mathfrak{u}(\mathrm{x})-\mathfrak{u}(\mathrm{y})|^{p}}{|\mathrm{x}-\mathrm{y}|^{N+s p}} \mathrm{dxdy}+\mathcal{K}\left([\mathfrak{v}]_{q}^{s}\right) \int_{\mathcal{U} \times} \frac{|\mathfrak{v}(\mathrm{x})-\mathfrak{v}(\mathrm{y})|^{q}}{|\mathrm{x}-\mathrm{y}|^{N+s q}} \mathrm{dxdy}\right]+\frac{1}{\sigma} \int_{\mathcal{U}}|\mathfrak{u}|^{\sigma}|\mathfrak{v}|^{\sigma} \mathrm{dx} \\
& -\frac{1}{\sigma} \int_{\mathcal{U}}|\mathfrak{u}|^{\sigma}|\mathfrak{v}|^{\sigma} \log |\mathfrak{u}\mathfrak{v}| \mathrm{dx} \\
& \leq \left(\frac{1}{p}-\frac{1}{\sigma}\right)\left[\mathcal{K}\left([\mathfrak{u}]_{p}^{s}\right) \int_{\mathcal{U} \times \mathcal{U}} \frac{|\mathfrak{u}(\mathrm{x})-\mathfrak{u}(\mathrm{y})|^{p}}{|\mathrm{x}-\mathrm{y}|^{N+s p}} \mathrm{dxdy}+\mathcal{K}\left([\mathfrak{v}]_{q}^{s}\right) \int_{\mathcal{U} \times \mathcal{U}} \frac{|\mathfrak{v}(\mathrm{x})-\mathfrak{v}(\mathrm{y})|^{q}}{|\mathrm{x}-\mathrm{y}|^{N+s q}} \mathrm{dxdy}\right]+\frac{1}{\sigma} \psi(\mathfrak{u}, \mathfrak{v}) \\
& +\frac{1}{\sigma} \int_{\mathcal{U}}|\mathfrak{u}|^{\sigma}|\mathfrak{v}|^{\sigma} \mathrm{dx}.
\end{split}
\end{align}

From \eqref{equation22}, we deduce that

\begin{align*}
\varphi(\mathfrak{u}, \mathfrak{v})-\frac{1}{\sigma} \psi(\mathfrak{u}, \mathfrak{v})& \leq\left(\frac{1}{p}-\frac{1}{\sigma}\right)\left[\mathcal{K}\left([\mathfrak{u}]_{p}^{s}\right) \int_{\mathcal{U} \times \mathcal{U}} \frac{|\mathfrak{u}(\mathrm{x})-\mathfrak{u}(\mathrm{y})|^{p}}{|\mathrm{x}-\mathrm{y}|^{N+s p}} \mathrm{dxdy}+\mathcal{K}\left([\mathfrak{v}]_{q}^{s}\right) \int_{\mathcal{U} \times \mathcal{U}} \frac{|\mathfrak{v}(\mathrm{x})-\mathfrak{v}(\mathrm{y})|^{q}}{|\mathrm{x}-\mathrm{y}|^{N+s q}} \mathrm{dxdy}\right]\\
&+\frac{1}{\sigma} \int_{\mathcal{U}}|\mathfrak{u}|^{\sigma}|\mathfrak{v}|^{\sigma} \mathrm{dx}.
\end{align*}

Now, we prove (2). Using (5) in Lemma  \ref{lemma3}, we have that
\begin{align}\label{equation23}
\begin{split}
\varphi(\mathfrak{u}, \mathfrak{v}) & =\frac{1}{p}\widehat{\mathcal{K}}\left([\mathfrak{u}]_{p}^{s}\right)+\frac{1}{q}\widehat{\mathcal{K}}\left([\mathfrak{v}]_{q}^{s}\right)+\frac{1}{\sigma^{2}} \int_{\mathcal{U}}|\mathfrak{u}|^{\sigma}|\mathfrak{v}|^{\sigma} \mathrm{dx}-\frac{1}{\sigma} \int_{\mathcal{U}}|\mathfrak{u}|^{\sigma}|\mathfrak{v}|^{\sigma} \log |\mathfrak{u}| \mathrm{dx} \\
& \geq \widehat{\mathcal{K}}\left([\mathfrak{u}]_{p}^{s}\right)+\widehat{\mathcal{K}}\left([\mathfrak{u}]_{q}^{s}\right)-\frac{1}{\sigma} \int_{\mathcal{U}}|\mathfrak{u}|^{\sigma}|\mathfrak{v}|^{\sigma} \log |\mathfrak{u}\mathfrak{v}| \mathrm{dx} \\
& \geq \frac{1}{\beta+1} \mathcal{K}\left([\mathfrak{u}]_{p}^{s}\right) \int_{\mathcal{U} \times \mathcal{U}} \frac{1}{p} \frac{|\mathfrak{u}(\mathrm{x})-\mathfrak{u}(\mathrm{y})|^{p}}{|\mathrm{x}-\mathrm{y}|^{N+s p}} \mathrm{dxdy}+\frac{1}{\beta+1} \mathcal{K}\left([\mathfrak{v}]_{q}^{s}\right) \int_{\mathcal{U} \times \mathcal{U}} \frac{1}{q} \frac{|\mathfrak{v}(\mathrm{x})-\mathfrak{v}(\mathrm{y})|^{q}}{|\mathrm{x}-\mathrm{y}|^{N+s q}} \mathrm{dxdy}-\frac{1}{\sigma} \int_{\mathcal{U}}|\mathfrak{u}|^{\sigma}|\mathfrak{v}|^{\sigma} \log |\mathfrak{u}| \mathrm{dx} \\
& \geq \frac{1}{q(\beta+1)}\left[\mathcal{K}\left([\mathfrak{u}]_{p}^{s}\right) \int_{\mathcal{U} \times \mathcal{U}} \frac{|\mathfrak{u}(\mathrm{x})-\mathfrak{u}(\mathrm{y})|^{p}}{|\mathrm{x}-\mathrm{y}|^{N+s p}} \mathrm{dxdy}+\mathcal{K}\left([\mathfrak{v}]_{q}^{s}\right) \int_{\mathcal{U} \times \mathcal{U}} \frac{|\mathfrak{v}(\mathrm{x})-\mathfrak{v}(\mathrm{y})|^{q}}{|\mathrm{x}-\mathrm{y}|^{N+s q}} \mathrm{dxdy}\right]-\frac{1}{\sigma} \int_{\mathcal{U}}|\mathfrak{u}|^{\sigma}|\mathfrak{v}|^{\sigma} \log |\mathfrak{u}\mathfrak{v}| \mathrm{dx} \\
& \geq\left(\frac{1}{q(\beta+1)}-\frac{1}{\sigma}\right)\left[\mathcal{K}\left([\mathfrak{u}]_{p}^{s}\right) \int_{\mathcal{U} \times \mathcal{U}} \frac{|\mathfrak{u}(\mathrm{x})-\mathfrak{u}(\mathrm{y})|^{p}}{|\mathrm{x}-\mathrm{y}|^{N+s p}} \mathrm{dxdy}+\mathcal{K}\left([\mathfrak{v}]_{q}^{s}\right) \int_{\mathcal{U} \times \mathcal{U}} \frac{|\mathfrak{v}(\mathrm{x})-\mathfrak{v}(\mathrm{y})|^{q}}{|\mathrm{x}-\mathrm{y}|^{N+s q}} \mathrm{dxdy}\right]+\frac{1}{\sigma} \psi(\mathfrak{u}, \mathfrak{v}).
\end{split}
\end{align}

From \eqref{equation23}, we deduce that

$$
\varphi(\mathfrak{u}, \mathfrak{v})-\frac{1}{\sigma} \psi(\mathfrak{u}, \mathfrak{v}) \geq\left(\frac{1}{q(\beta+1)}-\frac{1}{\sigma}\right)\left[\mathcal{K}\left([\mathfrak{u}]_{p}^{s}\right) \int_{\mathcal{U} \times \mathcal{U}} \frac{|\mathfrak{u}(\mathrm{x})-\mathfrak{u}(\mathrm{y})|^{p}}{|\mathrm{x}-\mathrm{y}|^{N+s p}} \mathrm{dxdy}+\mathcal{K}\left([\mathfrak{v}]_{q}^{s}\right) \int_{\mathcal{U} \times \mathcal{U}} \frac{|\mathfrak{v}(\mathrm{x})-\mathfrak{v}(\mathrm{y})|^{q}}{|\mathrm{x}-\mathrm{y}|^{N+s q}} \mathrm{dxdy}\right]
$$

Suppose that $\psi(\mathfrak{u}, \mathfrak{v})<0.$ Thanks to Lemma \ref{lemma4}, there exists $\varepsilon_{*} \in(0,1)$ such that $\psi\left(\varepsilon_{*} \mathfrak{u}, \varepsilon_{*} \mathfrak{v}\right)=0.$ From the definition of $d$ and (1) and (2) in Lemma  \ref{lemma5}, we obtain that

\begin{align}
\begin{split}
d & \leq \varphi\left(\varepsilon_{*} \mathfrak{u}, \varepsilon_{*} \mathfrak{v}\right) \\
& =\varphi\left(\varepsilon_{*} \mathfrak{u}, \varepsilon_{*} \mathfrak{v}\right)-\frac{1}{\sigma} \psi\left(\varepsilon_{*} \mathfrak{u}, \varepsilon_{*} \mathfrak{v}\right) \\
& \leq\left(\frac{1}{p}-\frac{1}{\sigma}\right)\left[\mathcal{K}\left(\left[\varepsilon_{*} \mathfrak{u}\right]_{p}^{s}\right) \int_{\mathcal{U} \times \mathcal{U}} \frac{\varepsilon_{*}^{p}|\mathfrak{u}(\mathrm{x})-\mathfrak{u}(\mathrm{y})|^{p}}{|\mathrm{x}-\mathrm{y}|^{N+s p}} \mathrm{dxdy}+\mathcal{K}\left(\left[\varepsilon_{*} \mathfrak{v}\right]_{q}^{s}\right) \int_{\mathcal{U} \times \mathcal{U}} \frac{\varepsilon_{*}^{q}|\mathfrak{v}(\mathrm{x})-\mathfrak{v}(\mathrm{y})|^{q}}{|\mathrm{x}-\mathrm{y}|^{N+s q}} \mathrm{dxdy}\right]+\frac{1}{\sigma} \int_{\mathcal{U}}\left|\varepsilon_{*} \mathfrak{u}\right|^{\sigma}\left|\varepsilon_{*} \mathfrak{v}\right|^{\sigma} \mathrm{dx} \\
& \leq\left(\frac{1}{p}-\frac{1}{\sigma}\right)\left[\mathcal{K}\left([\mathfrak{u}]_{p}^{s}\right) \int_{\mathcal{U} \times \mathcal{U}} \frac{|\mathfrak{u}(\mathrm{x})-\mathfrak{u}(\mathrm{y})|^{p}}{|\mathrm{x}-\mathrm{y}|^{N+s p}} \mathrm{dxdy}+\mathcal{K}\left([\mathfrak{v}]_{q}^{s}\right) \int_{\mathcal{U} \times \mathcal{U}} \frac{|\mathfrak{v}(\mathrm{x})-\mathfrak{v}(\mathrm{y})|^{q}}{|\mathrm{x}-\mathrm{y}|^{N+s q}} \mathrm{dxdy}\right]+\frac{1}{\sigma} \int_{\mathcal{U}}|\mathfrak{u}|^{\sigma}|\mathfrak{v}|^{\sigma} \mathrm{dx} \\
& \leq\left(\frac{1}{p}-\frac{1}{\sigma}\right) \frac{\varphi(\mathfrak{u}, \mathfrak{v})-\frac{1}{\sigma} \psi(\mathfrak{u}, \mathfrak{v})}{\frac{1}{q(\beta+1)}-\frac{1}{\sigma}}+\frac{1}{\sigma} \int_{\mathcal{U}}|\mathfrak{u}|^{\sigma}|\mathfrak{v}|^{\sigma} \mathrm{dx},
\end{split}
\end{align}

which implies that

$$
d_{*}:=\frac{(d-C(\sigma, \mathfrak{u}, \mathfrak{v}, \mathcal{U}))\left(\frac{1}{q(\beta+1)-\frac{1}{\sigma}}\right)}{\frac{1}{p}-\frac{1}{\sigma}} \leq \varphi(\mathfrak{u}, \mathfrak{v})-\frac{1}{\sigma} \psi(\mathfrak{u}, \mathfrak{v}),
$$

with $C(\sigma, \mathfrak{u}, \mathfrak{v}, \mathcal{U})=\frac{1}{\sigma} \int_{\mathcal{U}}|\mathfrak{u}|^{\sigma}|\mathfrak{v}|^{\sigma} \mathrm{dx}.$ Now, we prove (4). Assume that $\psi(\mathfrak{u}, \mathfrak{v}) \geq 0.$ From Lemma \ref{lemma4}, there exist $\varepsilon_{*}>1$ such that $\psi\left(\varepsilon_{*} \mathfrak{u}, \varepsilon_{*} \mathfrak{v}\right)=0.$ So, we get that

\begin{align}\label{equation15}
\begin{split}
& d \leq \varphi\left(\varepsilon_{*} \mathfrak{u}, \varepsilon_{*} \mathfrak{v}\right)=\varphi\left(\varepsilon_{*} \mathfrak{u}, \varepsilon_{*} \mathfrak{v}\right)-\frac{1}{\sigma} \psi\left(\varepsilon_{*} \mathfrak{u}, \varepsilon_{*} \mathfrak{v}\right) \\
& \leq\left(\frac{1}{p}-\frac{1}{\sigma}\right)\left[\mathcal{K}\left(\left[\varepsilon_{*} \mathfrak{u}\right]_{p}^{s}\right) \int_{\mathcal{U} \times \mathcal{U}} \frac{\varepsilon_{*}^{p}|\mathfrak{u}(\mathrm{x})-\mathfrak{u}(\mathrm{y})|^{p}}{|\mathrm{x}-\mathrm{y}|^{N+s p}} \mathrm{dxdy}+\mathcal{K}\left(\left[\varepsilon_{*} \mathfrak{v}\right]_{q}^{s}\right) \int_{\mathcal{U} \times \mathcal{U}} \frac{\varepsilon_{*}^{q}|\mathfrak{v}(\mathrm{x})-\mathfrak{v}(\mathrm{y})|^{q}}{|\mathrm{x}-\mathrm{y}|^{N+s q}} \mathrm{dxdy}\right]+\frac{1}{\sigma} \int_{\mathcal{U}}\left|\varepsilon_{*} \mathfrak{u}\right|^{\sigma}\left|\varepsilon_{*} \mathfrak{v}\right|^{\sigma} \mathrm{dx} \\
& \leq\left(\frac{1}{p}-\frac{1}{\sigma}\right)\left[\varepsilon_{*}^{p(\beta+1)} \mathcal{K}\left([\mathfrak{u}]_{p}^{s}\right) \int_{\mathcal{U} \times \mathcal{U}} \frac{|\mathfrak{u}(\mathrm{x})-\mathfrak{u}(\mathrm{y})|^{p}}{|\mathrm{x}-\mathrm{y}|^{N+s p}} \mathrm{dxdy}+\varepsilon_{*}^{q(\beta+1)} \mathcal{K}\left([\mathfrak{v}]_{q}^{s}\right) \int_{\mathcal{U} \times \mathcal{U}} \frac{|\mathfrak{u}(\mathrm{x})-\mathfrak{y}(\mathrm{y})|^{q}}{|\mathrm{x}-\mathrm{y}|^{N+s q}} \mathrm{dxdy}\right]+\frac{1}{\sigma} \int_{\mathcal{U}}\left|\varepsilon_{*} \mathfrak{u}\right|^{\sigma}\left|\varepsilon_{*} \mathfrak{v}\right|^{\sigma} \mathrm{dx} \\
& \leq\left(\frac{1}{p}-\frac{1}{\sigma}\right) \varepsilon_{*}^{q(\beta+1)}\left[\mathcal{K}\left([\mathfrak{u}]_{p}^{s}\right) \int_{\mathcal{U} \times \mathcal{U}} \frac{|\mathfrak{u}(\mathrm{x})-\mathfrak{u}(\mathrm{y})|^{p}}{|\mathrm{x}-\mathrm{y}|^{N+s p}} \mathrm{dxdy}+\mathcal{K}\left([\mathfrak{v}]_{q}^{s}\right) \int_{\mathcal{U} \times \mathcal{U}} \frac{|\mathfrak{v}(\mathrm{x})-\mathfrak{v}(\mathrm{y})|^{q}}{|\mathrm{x}-\mathrm{y}|^{N+s q}} \mathrm{dxdy}\right]+\frac{\varepsilon_{*}^{2 \sigma}}{\sigma} \int_{\mathcal{U}}|\mathfrak{u}|^{\sigma}|\mathfrak{v}|^{\sigma} \mathrm{dx} \\
& \leq \varepsilon_{*}^{q(\beta+1)}\left(\frac{1}{p}-\frac{1}{\sigma}\right) \frac{\varphi(\mathfrak{u}, \mathfrak{v})-\frac{1}{\sigma} \psi(\mathfrak{u}, \mathfrak{v})}{\frac{1}{q(\beta+1)}-\frac{1}{\sigma}}+\frac{\varepsilon_{*}^{2 \sigma}}{\sigma} \int_{\mathcal{U}}|\mathfrak{u}|^{\sigma}|\mathfrak{v}|^{\sigma} \mathrm{dx}.
\end{split}
\end{align}

From \eqref{equation15}, we deduce that

$$
d_{*} \varepsilon_{*}^{-q(\beta+1)}-C^{\prime}\left(\varepsilon_{*}, \sigma, d, \mathfrak{u}, \mathfrak{v}\right) \leq \varphi(\mathfrak{u}, \mathfrak{v})-\frac{1}{\sigma} \psi(\mathfrak{u}, \mathfrak{v}),
$$
where $C^{\prime}\left(\varepsilon_{*}, \sigma, d, \mathfrak{u}, \mathfrak{v}\right)=\frac{\varepsilon_{*}^{2 \sigma}}{\sigma} \int_{\mathcal{U}}|\mathfrak{u}|^{\sigma}|\mathfrak{v}|^{\sigma} \mathrm{dx}.$ Finally, we prove the second inequality of (4). Due to $\varepsilon_{*}>1,$ we get that

\begin{align}\label{equation20}
\begin{split}
& 0=\psi\left(\varepsilon_{*} \mathfrak{u}, \varepsilon_{*} \mathfrak{v}\right) \\
& =\mathcal{K}\left(\left[\varepsilon_{*} \mathfrak{u}\right]_{p}^{s}\right) \int_{\mathcal{U} \times \mathcal{U}} \frac{\varepsilon_{*}^{p}|\mathfrak{u}(\mathrm{x})-\mathfrak{u}(\mathrm{y})|^{p}}{|\mathrm{x}-\mathrm{y}|^{N+s p}} \mathrm{dxdy}+\mathcal{K}\left(\left[\varepsilon_{*} \mathfrak{v}\right]_{q}^{s}\right) \int_{\mathcal{U} \times \mathcal{U}} \frac{\varepsilon_{*}^{q}|\mathfrak{v}(\mathrm{x})-\mathfrak{v}(\mathrm{y})|^{q}}{|\mathrm{x}-\mathrm{y}|^{N+s q}} \mathrm{dxdy}-2 \int_{\mathcal{U}}\left|\varepsilon_{*} \mathfrak{u}\right|^{\sigma+1}\left|\varepsilon_{*} \mathfrak{v}\right|^{\sigma+1} \log \left|\varepsilon_{*}^{2} \mathfrak{u}\mathfrak{v}\right| \mathrm{dx} \\
& \leq \varepsilon_{*}^{p(\beta+1)} \mathcal{K}\left([\mathfrak{u}]_{p}^{s}\right) \int_{\mathcal{U} \times \mathcal{U}} \frac{|\mathfrak{u}(\mathrm{x})-\mathfrak{u}(\mathrm{y})|^{p}}{|\mathrm{x}-\mathrm{y}|^{N+s p}} \mathrm{dxdy}+\varepsilon_{*}^{q(\beta+1)} \mathcal{K}\left([\mathfrak{v}]_{q}^{s}\right) \int_{\mathcal{U} \times \mathcal{U}} \frac{|\mathfrak{v}(\mathrm{x})-\mathfrak{v}(\mathrm{y})|^{q}}{|\mathrm{x}-\mathrm{y}|^{N+s q}} \mathrm{dxdy}\\
&-2 \int_{\mathcal{U}}\left|\varepsilon_{*} \mathfrak{l}\right|^{\sigma+1}\left|\varepsilon_{*} \mathfrak{v}\right|^{\sigma+1} \log \left|\varepsilon_{*}^{2}  \mathfrak{v}\mathfrak{u}\right| \mathrm{dx} \\
& \leq \varepsilon_{*}^{q(\beta+1)}\left[\mathcal{K}\left([\mathfrak{u}]_{p}^{s}\right) \int_{\mathcal{U} \times \mathcal{U}} \frac{|\mathfrak{u}(\mathrm{x})-\mathfrak{u}(\mathrm{y})|^{p}}{|\mathrm{x}-\mathrm{y}|^{N+s p}} \mathrm{dxdy}+\mathcal{K}\left([\mathfrak{v}]_{q}^{s}\right) \int_{\mathcal{U} \times \mathcal{U}} \frac{|\mathfrak{v}(\mathrm{x})-\mathfrak{v}(\mathrm{y})|^{q}}{|\mathrm{x}-\mathrm{y}|^{N+s q}} \mathrm{dxdy}\right]-2 \int_{\mathcal{U}}\left|\varepsilon_{*} \mathfrak{u}\right|^{\sigma+1}\left|\varepsilon_{*} \mathfrak{v}\right|^{\sigma+1} \log \left|\varepsilon_{*}^{2} \mathfrak{u}\mathfrak{v}\right| \mathrm{dx} \\
& \leq \varepsilon_{*}^{q(\beta+1)}\left[\mathcal{K}\left([\mathfrak{u}]_{p}^{s}\right) \int_{\mathcal{U} \times \mathcal{U}} \frac{|\mathfrak{u}(\mathrm{x})-\mathfrak{u}(\mathrm{y})|^{p}}{|\mathrm{x}-\mathrm{y}|^{N+s p}} \mathrm{dxdy}+\mathcal{K}\left([\mathfrak{v}]_{q}^{s}\right) \int_{\mathcal{U} \times \mathcal{U}} \frac{|\mathfrak{v}(\mathrm{x})-\mathfrak{v}(\mathrm{y})|^{q}}{|\mathrm{x}-\mathrm{y}|^{N+s q}} \mathrm{dxdy}\right]-\varepsilon_{*}^{2(\sigma+1)} \int_{\mathcal{U}}|\mathfrak{u}|^{\sigma+1}|\mathfrak{v}|^{\sigma+1} \log |\mathfrak{u} \mathfrak{v}| \mathrm{dx} \\
& \leq\left(\varepsilon_{*}^{q(\beta+1)}-\varepsilon_{*}^{2(\sigma+1)}\right)\left[\mathcal{K}\left([\mathfrak{u}]_{p}^{s}\right) \int_{\mathcal{U} \times \mathcal{U}} \frac{|\mathfrak{u}(\mathrm{x})-\mathfrak{u}(\mathrm{y})|^{p}}{|\mathrm{x}-\mathrm{y}|^{N+s p}} \mathrm{dxdy}+\mathcal{K}\left([\mathfrak{v}]_{q}^{s}\right) \int_{\mathcal{U} \times \mathcal{U}} \frac{|\mathfrak{v}(\mathrm{x})-\mathfrak{v}(\mathrm{y})|^{q}}{|\mathrm{x}-\mathrm{y}|^{N+s q}} \mathrm{dxdy}\right]+\varepsilon_{*}^{2(\sigma+1)} \psi(\mathfrak{u}, \mathfrak{v}).
\end{split}
\end{align}

From \eqref{equation20}, we deduce that $$\psi(\mathfrak{u}, \mathfrak{v}) \geq\left(1-\varepsilon_{*}^{q(\beta+1)-2(\sigma+1)}\right)\left[\mathcal{K}\left([\mathfrak{u}]_{p}^{s}\right) \int_{\mathcal{U} \times \mathcal{U}} \frac{\mid \mathfrak{u}(\mathrm{x})-\mathfrak{u}(\mathrm{y}) p^{p}}{|\mathrm{x}-\mathrm{y}|^{(+s p}} \mathrm{dxdy}+\mathcal{K}\left([\mathfrak{v}]_{q}^{s}\right) \int_{\mathcal{U} \times \mathcal{U}} \frac{\mid \mathfrak{v}(\mathrm{x})-\mathfrak{v}(\mathrm{y}))^{q}}{|\mathrm{x}-\mathrm{y}|^{(y+q q}} \mathrm{dxdy}\right].$$
\end{proof}
\begin{theorem}

 If Kirchhoff's function satisfies the conditions $\left(\mathcal{H}_{1}\right)-\left(\mathcal{H}_{2}\right)$ and the numbers $q, p, \sigma$ satisfy the condition \eqref{equation1}. Then, the potential well depth $d$ is positive.
\end{theorem}

\begin{proof} We assume that $\mathcal{K}(1)>0.$ Using (1) and (3) in Lemma  \ref{lemma3}, we get that

\begin{align}\label{equ1}
\begin{split}
& \mathcal{K}\left([\mathfrak{u}]_{p}^{s}\right) \int_{\mathcal{U} \times \mathcal{U}} \frac{|\mathfrak{u}(\mathrm{x})-\mathfrak{u}(\mathrm{y})|^{p}}{|\mathrm{x}-\mathrm{y}|^{N+s p}} \mathrm{dxdy}+\mathcal{K}\left([\mathfrak{v}]_{q}^{s}\right) \int_{\mathcal{U} \times \mathcal{U}} \frac{|\mathfrak{v}(\mathrm{x})-\mathfrak{v}(\mathrm{y})|^{q}}{|\mathrm{x}-\mathrm{y}|^{N+s q}} \mathrm{dxdy} \\
& \geq \frac{1}{p^{\beta}} \mathcal{K}\left(\int_{\mathcal{U} \times \mathcal{U}} \frac{|\mathfrak{v}(\mathrm{x})-\mathfrak{v}(\mathrm{y})|^{q}}{|\mathrm{x}-\mathrm{y}|^{N+s q}} \mathrm{dxdy}\right) \int_{\mathcal{U} \times \mathcal{U}} \frac{|\mathfrak{u}(\mathrm{x})-\mathfrak{u}(\mathrm{y})|^{p}}{|\mathrm{x}-\mathrm{y}|^{N+s p}} \mathrm{dxdy}\\
&+\frac{1}{q^{\beta}} \mathcal{K}\left(\int_{\mathcal{U} \times \mathcal{U}} \frac{|\mathfrak{v}(\mathrm{x})-\mathfrak{v}(\mathrm{y})|^{q}}{|\mathrm{x}-\mathrm{y}|^{N+s q}} \mathrm{dxdy}\right) \int_{\mathcal{U} \times \mathcal{U}} \frac{|\mathfrak{v}(\mathrm{x})-\mathfrak{v}(\mathrm{y})|^{q}}{|\mathrm{x}-\mathrm{y}|^{N+s q}} \mathrm{dxdy} \\
& \geq \frac{\mathcal{K}(1)}{p^{\beta}} \min \left\{\|\mathfrak{u}\|_{W^{s, p}(\mathcal{U})}^{p},\|\mathfrak{u}\|_{W^{s, p}(\mathcal{U})}^{p(\beta+1)}\right\}+\frac{\mathcal{K}(1)}{q^{\beta}} \min \left\{\|\mathfrak{v}\|_{W^{s, q}(\mathcal{U})}^{q},\|\mathfrak{v}\|_{W^{s, q}(\mathcal{U})}^{q(\beta+1)}\right\} \\
& \geq \frac{\mathcal{K}(1)}{q^{\beta}} \min \left\{\|\mathfrak{u}\|_{W^{s, p}(\mathcal{U})}^{p} + \|\mathfrak{v}\|_{W^{s, q}(\mathcal{U})}^{q},\|\mathfrak{u}\|_{W^{s, p}(\mathcal{U})}^{p(\beta+1)}+\|\mathfrak{v}\|_{W^{s, q}(\mathcal{U})}^{q(\beta+1)}\right\}\\
&\geq  \frac{\mathcal{K}(1)}{q^{\beta}} \min \left\lbrace \Vert (\mathfrak{u}, \mathfrak{v})\Vert^{q}_{Y},
\Vert (\mathfrak{u}, \mathfrak{v})\Vert^{q(\beta+1)}_{Y}
 \right\rbrace.
\end{split}
\end{align}
Let $\left( \mathfrak{u}, \mathfrak{v} \right)\in  \mathcal{N}.$ Using Lemma \ref{injection} and embedding compact theorem, we have that
\begin{align}\label{equ2}
\begin{split}
& \mathcal{K}\left([\mathfrak{u}]_{p}^{s}\right) \int_{\mathcal{U} \times \mathcal{U}} \frac{|\mathfrak{u}(\mathrm{x})-\mathfrak{u}(\mathrm{y})|^{p}}{|\mathrm{x}-\mathrm{y}|^{N+s p}} \mathrm{dxdy}+\mathcal{K}\left([\mathfrak{v}]_{q}^{s}\right) \int_{\mathcal{U} \times \mathcal{U}} \frac{|\mathfrak{v}(\mathrm{x})-\mathfrak{v}(\mathrm{y})|^{q}}{|\mathrm{x}-\mathrm{y}|^{N+s q}} \mathrm{dxdy} \\
&=2  \int_{\mathcal{U}}\left|\mathfrak{u}\right|^{\sigma+1}\left|\mathfrak{v}\right|^{\sigma+1} \log \left|\mathfrak{u} \mathfrak{v}\right| \mathrm{dx}\\
&\leq 2\left[  \log [\mathfrak{u}]_{s, p} \int_{\mathcal{U}}|\mathfrak{u}|^{\sigma+1} \mathrm{dx}+\log [\mathfrak{v}]_{s, q} \int_{\mathcal{U}}|\mathfrak{v}|^{\sigma+1} \mathrm{dx}+S\left([\mathfrak{u}]_{s, p}^{p_s^*}+[ \mathfrak{v}]_{s, p}^{p_s^*}+[\mathfrak{u}]_{s, q}^{q_s^*}+[\mathfrak{v}]_{s, q}^{q_s^*}+[\mathfrak{u}]_{s, p}^{\sigma}+[\mathfrak{u}]_{s, q}^\sigma\right)\right]\\
&\leq 2\left[c_{1}(\mathfrak{u}, s, p, \mathcal{U})\int_{\mathcal{U}}|\mathfrak{u}|^{\sigma+1} \mathrm{dx}+
c_{2}(\mathfrak{v}, s, q, \mathcal{U})\int_{\mathcal{U}}|\mathfrak{v}|^{\sigma+1} \mathrm{dx}+ S^{'}_{1} \right] \\
&\leq  S^{''}_{1}\left(C_{1} \Vert \mathfrak{u}\Vert^{\sigma+1}_{L ^{\sigma+1}(\mathcal{U})}+
C_{2} \Vert \mathfrak{v}\Vert^{\sigma+1}_{L ^{\sigma+1}(\mathcal{U})}+1
 \right) \\
 &\leq S^{''}_{1} C^{'}_{1} \Vert\left(\mathfrak{u}, \mathfrak{v} \right)  \Vert ^{\sigma+1}_{Y},
\end{split}
\end{align}
where $C^{'}_{1} :=\max (c_{1}, c_{2}),$ $c_{1},$  and  $ c_{2}$
are embedding constant of  $W^{s, p}(\mathcal{U}),  W^{s, q}(\mathcal{U})\hookrightarrow  L ^{\sigma+1}(\mathcal{U}).$
 If $\Vert(\mathfrak{u}, \mathfrak{v}) \Vert_{Y}\leq 1.  $  Combining \eqref{equ1} with  \eqref{equ2}, we obtain that

 \begin{equation*}
 \frac{\mathcal{K}(1)}{q^{\beta}} \Vert\left(\mathfrak{u}, \mathfrak{v} \right)  \Vert ^{\sigma+1}_{Y}\leq S^{''}_{1} C^{'}_{1} \Vert\left(\mathfrak{u}, \mathfrak{v} \right)  \Vert ^{\sigma+1}_{Y},
 \end{equation*}
which yields that
\begin{equation*}
\Vert\left(\mathfrak{u}, \mathfrak{v} \right)  \Vert_{Y}\geq  \left( \frac{\mathcal{K}(1)}{q^{\beta} S^{''}_{1} C^{'}_{1}}   \right) ^{\frac{1}{\sigma+1- q(\beta+1)}}.
 \end{equation*}
Using the fact $\Vert\left(\mathfrak{u}, \mathfrak{v} \right)  \Vert_{Y}>0$ and $\sigma+1>q(\beta+1).$ We get that, for any $\left(\mathfrak{u}, \mathfrak{v} \right) \in \mathcal{N}$
\begin{equation*}
\Vert\left(\mathfrak{u}, \mathfrak{v} \right)  \Vert_{Y}\geq   \min \left\lbrace \left( \frac{\mathcal{K}(1)}{q^{\beta} S^{''}_{1} C^{'}_{1}}   \right) ^{\frac{1}{\sigma+1- q(\beta+1)}}, 1\right\rbrace:=  \varsigma \in(0, 1].
 \end{equation*}
 Thanks to Lemma \ref{lemma5}, we have that
 \begin{align}\label{equ5}
 \begin{split}
\varphi(\mathfrak{u}, \mathfrak{v})&= \varphi(\mathfrak{u}, \mathfrak{v})-\frac{1}{\sigma} \psi(\mathfrak{u}, \mathfrak{v})\\
& \geq\left(\frac{1}{q(\beta+1)}-\frac{1}{\sigma}\right)\left[\mathcal{K}\left([\mathfrak{u}]_{p}^{s}\right) \int_{\mathcal{U} \times \mathcal{U}} \frac{|\mathfrak{u}(\mathrm{x})-\mathfrak{u}(\mathrm{y})|^{p}}{|\mathrm{x}-\mathrm{y}|^{N+s p}} d \mathrm{x} d \mathrm{y}+\mathcal{K}\left([\mathfrak{v}]_{q}^{s}\right) \int_{\mathcal{U} \times \mathcal{U}} \frac{|\mathfrak{v}(\mathrm{x})-\mathfrak{v}(\mathrm{y})|^{q}}{|\mathrm{x}-\mathrm{y}|^{N+s q}} d \mathrm{x} d \mathrm{y}\right]  \\
&\geq  \left(\frac{1}{q(\beta+1)}-\frac{1}{\sigma}\right) \varsigma ^{q(\beta+1)}.
 \end{split}
 \end{align}
 Using the definition of d and \eqref{equ5}, we have that
  \begin{align*}
  d= \inf_{\left(\mathfrak{u}, \mathfrak{v} \right) \in \mathcal{N} }\varphi \left(\mathfrak{u}, \mathfrak{v} \right) \geq \left(\frac{1}{q(\beta+1)}-\frac{1}{\sigma}\right) \varsigma ^{q(\beta+1)}>0.
  \end{align*}
\end{proof}

\begin{theorem}

  Assume that $\left(\mathfrak{u}_{0}, \mathfrak{v}_{0}\right) \in Y.$ Under assumptions $\left(\mathcal{H}_{1}\right)-\left(\mathcal{H}_{2}\right).$ Then, there exists $T>0$ such that the system \eqref{problem} has a weak solution in Y. Moreover, there holds

\begin{align}\label{eq7}
\int_{0}^{t}\left\|\mathfrak{u}_{s}(s)\right\|_{2}^{2} d s+\int_{0}^{t}\left\|\mathfrak{v}_{s}(s)\right\|_{2}^{2} d s+\varphi(\mathfrak{u}(t), \mathfrak{v}(t)) \leq \varphi\left(\mathfrak{u}_{0}, \mathfrak{v}_{0}\right), \quad \text { for all } t \in[0, T].
\end{align}
\end{theorem}
\begin{proof} Since $Y$ is separable space and is dense in $L^{2}(\mathcal{U}) \times L^{2}(\mathcal{U}),$ there exist a base $\overline{\mathcal{Z}}^{\mid 1 \cdot \|_{Y}}=Y$ and $\left(w_{i}, w_{j}\right)=\left(z_{i}, z_{j}\right)=$ $\delta_{i, j},$ for all $j, i=1,2, \ldots$ where $\delta_{i, j}$ denotes the Kronecker symbol. For any $n \in \mathbb{N},$ we construct the approximation solutions $\mathfrak{u}_{m}(\mathrm{x}, t)=\mathfrak{u}_{m}(t)=\sum_{i=1}^{m} g_{i m}(t) w_{i}(\mathrm{x})$ and $\mathfrak{v}_{m}(\mathrm{x}, t)=\mathfrak{v}_{m}(t)=\sum_{i=1}^{m} h_{i m}(t) z_{i}(\mathrm{x})$ satisfying

\begin{equation}\label{equation26}
\left\{\begin{array}{l}
\int_{\mathcal{U}} \mathfrak{u}_{t m} w_{i} \mathrm{dx}+\left\langle\mathfrak{u}_{m}, w_{i}\right\rangle_{W^{s, p}(\mathcal{U})}=\int_{\mathcal{U}}\left|\mathfrak{v}_{m}\right|^{\sigma}\left|\mathfrak{u}_{m}\right|^{\sigma-2} \mathfrak{u}_{m} w_{i} \log \left|\mathfrak{u}_{m} \mathfrak{v}_{m}\right|\mathrm{dx} \\
\int_{\mathcal{U}} \mathfrak{v}_{t m} z_{i} \mathrm{dx}+\left\langle\mathfrak{v}_{m}, z_{i}\right\rangle_{W^{s, q}(\mathcal{U})}=\int_{\mathcal{U}}\left|\mathfrak{u}_{m}\right|^{\sigma}\left|\mathfrak{v}_{m}\right|^{\sigma-2} \mathfrak{v}_{m} z_{i} \log \left|\mathfrak{u}_{m} \mathfrak{v}_{m}\right|\mathrm{dx} \\
\mathfrak{u}_{m}(0)=\mathfrak{u}_{0 m}, \quad \mathfrak{v}_{m}(0)=\mathfrak{v}_{0 m},
\end{array}\right.
\end{equation}

where $\left\langle\mathfrak{u}_{m}, w_{i}\right\rangle_{W^{s, p}(\mathcal{U})}=\mathcal{K}\left(\left[\mathfrak{u}_{m}\right]_{s}^{p}\right) \int_{\mathcal{U} \times \mathcal{U}} \frac{\left|\mathfrak{u}_{m}(\mathrm{x}, t)-\mathfrak{u}_{m}(\mathrm{y}, t)\right|^{q-2}\left(\mathfrak{u}_{m}(\mathrm{x}, t)-\mathfrak{u}_{m}(\mathrm{y}, t)\right)\left(w_{i}(\mathrm{x})-w_{i}(\mathrm{y})\right)}{|\mathrm{x}-\mathrm{y}|^{N+s q}} \mathrm{dx} d \mathrm{y}, \mathfrak{v}_{t m}=\frac{\partial}{\partial t} \mathfrak{v}_{m}(t),$

$\left\langle\mathfrak{v}_{m}, z_{i}\right\rangle_{W^{s, q}(\mathcal{U})}=\mathcal{K}\left(\left[\mathfrak{v}_{m}\right]_{s}^{q}\right) \int_{\mathcal{U} \times \mathcal{U}} \frac{\left|\mathfrak{v}_{m}(\mathrm{x}, t)-\mathfrak{v}_{m}(\mathrm{y}, t)\right|^{q-2}\left(\mathfrak{v}_{m}(\mathrm{x}, t)-\mathfrak{v}_{m}(\mathrm{y}, t)\right)\left(z_{i}(\mathrm{x})-z_{i}(\mathrm{y})\right)}{|\mathrm{x}-\mathrm{y}|^{N+s q}} \mathrm{dx} d \mathrm{y}$ and $\mathfrak{u}_{t m}=\frac{\partial}{\partial t} \mathfrak{u}_{m}(t).$ Since $\left(\mathfrak{u}_{0}, \mathfrak{v}_{0}\right) \in Y,$ there exists $\left\{\left(\xi_{m i}, \zeta_{m i}\right), i=1, \ldots, m\right\}$ such that
\begin{align}\label{equation42}
\mathfrak{u}_{0 m}=\sum_{i=1}^{m} \xi_{m i}(t) w_{i}(\mathrm{x}) \rightarrow \mathfrak{u}_{0} \text { in } W^{s, p}(\mathcal{U}) \text { and } \mathfrak{v}_{0 m}=\sum_{i=1}^{m} \zeta_{m i} z_{i}(\mathrm{x}) \rightarrow \mathfrak{v}_{0} \text { in } W^{s, q}(\mathcal{U}).
\end{align}

Then \eqref{equation26}  is equivalent to the following initial value problem for a system of nonlinear ordinary differential system:

$$
\left\{\begin{array}{l}
g_{i m}^{\prime}(t)=G_{i}(g, h), i=1, \ldots, m, \\
h_{i m}^{\prime}(t)=H_{i}(g, h), i=1, \ldots, m, \\
g_{i m}(0)=\xi_{i m}, \quad h_{i m}(0)=\zeta_{i m}, i=1, \ldots, m,
\end{array}\right.
$$

with $G_{i}\left( g, h\right)  =-\left\langle\mathfrak{u}_{m}, w_{i}\right\rangle_{W^{s, p}(\mathcal{U})}+\int_{\mathcal{U}}\vert \mathfrak{v}_{m}\vert ^{\sigma}\vert |\mathfrak{u}_{m}\vert ^{\sigma-2} \mathfrak{u}_{m} w_{i} \log \vert \mathfrak{u}_{m} \mathfrak{v}_{m}\vert \mathrm{dx}$
 and
  $$H_{i}\left( g, h\right)  =-\left\langle\mathfrak{v}_{m}, z_{i}\right\rangle_{W^{s, p}(\mathcal{U})}+\int_{\mathcal{U}}\vert \mathfrak{u}_{m}\vert ^{\sigma}\vert |\mathfrak{v}_{m}\vert ^{\sigma-2} \mathfrak{u}_{m} z_{i} \log \vert \mathfrak{u}_{m} \mathfrak{v}_{m}\vert \mathrm{dx}.$$

From the Picardo iteration method, there exists $t_{0, m}>0$ depending on $|\xi|$ and $|\zeta|$ such that the system \eqref{equation26}  admits a local solution $\left(g_{i m}, g_{i m}\right) \in C^{1}\left(\left[0, t_{0, m}\right]\right) \times C^{1}\left(\left[0, t_{0, m}\right]\right).$ The following is an attempt to obtain some a priori estimations of the approximate solutions. We multiply the first equation in \eqref{equation26}  by $g_{i m}$ and the second equation by $h_{i m},$ after that summing with respect to $i,$ we get that

\begin{align}\label{equation34}
\begin{split}
& \frac{1}{2} \frac{d}{d t}\left\|\mathfrak{u}_{m}\right\|_{2}^{2}+\frac{1}{2} \frac{d}{d t}\left\|\mathfrak{v}_{m}\right\|_{2}^{2}+\mathcal{K}\left(\left[\mathfrak{u}_{m}\right]_{s}^{p}\right) \int_{\mathcal{U} \times \mathcal{U}} \frac{\left|\mathfrak{u}_{m}(\mathrm{x}, t)-\mathfrak{u}_{m}(\mathrm{y}, t)\right|^{p}}{|\mathrm{x}-\mathrm{y}|^{N+s p}} \mathrm{dx} d \mathrm{y}+\mathcal{K}\left(\left[\mathfrak{v}_{m}\right]_{s}^{q}\right) \int_{\mathcal{U} \times \mathcal{U}} \frac{\left|\mathfrak{v}_{m}(\mathrm{x}, t)-\mathfrak{v}_{m}(\mathrm{y}, t)\right|^{q}}{|\mathrm{x}-\mathrm{y}|^{N+s q}} \mathrm{dx} d \mathrm{y} \\
& =2 \int_{\mathcal{U}}\left|\mathfrak{u}_{m}\right|^{\sigma}\left|\mathfrak{v}_{m}\right|^{\sigma} \log \left|\mathfrak{u}_{m} \mathfrak{v}_{m}\right| \mathrm{dx}.
\end{split}
\end{align}

Considering set

$$
\mathcal{U}_{1}:=\left\{\mathrm{x} \in \mathcal{U}:\left|\mathfrak{u}_{m} \mathfrak{v}_{m}(\mathrm{x})\right| \leq 1\right\} \text { and } \mathcal{U}_{2}:=\left\{\mathrm{x} \in \mathcal{U}:\left|\mathfrak{u}_{m} \mathfrak{v}_{m}(\mathrm{x})\right|>1\right\}.
$$

It is easy to see that $\mathcal{U}=\mathcal{U}_{1} \cup \mathcal{U}_{2}$ and $\mathcal{U}_{1} \cap \mathcal{U}_{2}=\emptyset.$ So, we get that

$$
\begin{aligned}
\int_{\mathcal{U}}\left|\mathfrak{u}_{m}\right|^{\sigma}\left|\mathfrak{v}_{m}\right|^{\sigma} \log \left|\mathfrak{u}_{m} \mathfrak{v}_{m}\right|\mathrm{dx} & =\int_{\mathcal{U}_{1}}\left|\mathfrak{u}_{m}\right|^{\sigma}\left|\mathfrak{v}_{m}\right|^{\sigma} \log \left|\mathfrak{u}_{m} \mathfrak{v}_{m}\right| \mathrm{dx}+\int_{\mathcal{U}_{2}}\left|\mathfrak{u}_{m}\right|^{\sigma}\left|\mathfrak{v}_{m}\right|^{\sigma} \log \left|\mathfrak{u}_{m} \mathfrak{v}_{m}\right|\mathrm{dx} \\
& \leq \int_{\mathcal{U}_{1}}\left|\mathfrak{u}_{m}\right|^{\sigma}\left|\mathfrak{v}_{m}\right|^{\sigma} \frac{1}{\exp (1)} \mathrm{dx}+\int_{\mathcal{U}_{2}} \frac{1}{\sigma} \frac{1}{\exp (1)} \mathrm{dx} \\
& \leq \frac{|\mathcal{U}|}{\sigma \exp (1)}+\frac{1}{\sigma \exp (1)} \int_{\mathcal{U}}\left|\mathfrak{u}_{m} \mathfrak{v}_{m}\right|^{\sigma+1} \mathrm{dx} .
\end{aligned}
$$

Using H\"{o}lder's inequality and the classical interpolation inequality, we get that

\begin{equation}\label{equation27}
\int_{\mathcal{U}}\left|\mathfrak{u}_{m}\right|^{\sigma}\left|\mathfrak{v}_{m}\right|^{\sigma} \log \left|\mathfrak{u}_{m} \mathfrak{v}_{m}\right|\mathrm{dx} \leq\left\|\mathfrak{u}_{m}\right\|_{2}^{((\sigma+1)}\left\|\mathfrak{v}_{m}\right\|_{2}^{\kappa(\sigma+r)}\left\|\mathfrak{v}_{m}\right\|_{q}^{(1-\kappa)(\sigma+r)}\left\|\mathfrak{u}_{m}\right\|_{p}^{(1-\mu)(\sigma+r)},
\end{equation}

where $\mu, \kappa \in(0,1),$ and $r>0$ satisfy $\frac{1}{2(r+\mu)}=\frac{\mu}{2}+\frac{1-\mu}{p}, \frac{1}{2(r+\mu)}=\frac{\kappa}{2}+\frac{1-\kappa}{q},$ and $2<2(\sigma+r)<\min \{q, p\}.$ From \eqref{equation27} , and the compacts embeddings $W^{s, p}(\mathcal{U}) \hookrightarrow L^{p}(\mathcal{U}), W^{s, q}(\mathcal{U}) \hookrightarrow L^{q}(\mathcal{U}),$ we have that

$$
\int_{\mathcal{U}}\left|\mathfrak{u}_{m}\right|^{\sigma}\left|\mathfrak{v}_{m}\right|^{\sigma} \log \left|\mathfrak{u}_{m} \mathfrak{v}_{m}\right|\mathrm{dx} \leq C\left\|\mathfrak{u}_{m}\right\|_{2}^{(\sigma+1)}\left\|\mathfrak{v}_{m}\right\|_{2}^{\kappa(\sigma+r)}\left\|\mathfrak{v}_{m}\right\|_{W^{s, p}(\mathcal{U})}^{(1-\kappa)(\sigma+r)}\left\|\mathfrak{u}_{m}\right\|_{W^{s, q}(\mathcal{U})}^{(1-\mu)(\sigma+r)}.
$$

If $\left\|\mathfrak{u}_{m}\right\|_{W^{s q}(\mathcal{U})},\left\|\mathfrak{v}_{m}\right\|_{W^{s, q}(\mathcal{U})} \leq 1.$ Once that is done, there is nothing more to prove. From now, we distinguish three cases:\\
 Case 1: $\left\|\mathfrak{u}_{m}\right\|_{W^{s, q}(\mathcal{U})}>1$ and $\left\|\mathfrak{v}_{m}\right\|_{W^{s, q}(\mathcal{U})} \leq 1.$

$$
\int_{\mathcal{U}}\left|\mathfrak{u}_{m}\right|^{\sigma}\left|\mathfrak{v}_{m}\right|^{\sigma} \log \left|\mathfrak{u}_{m} \mathfrak{v}_{m}\right|\mathrm{dx} \leq \int_{\mathcal{U}}\left|\mathfrak{u}_{m}\right|^{\sigma} \log \left(\left|\mathfrak{u}_{m}\right|\right) \mathrm{dx}.
$$

Considering the following sets

$$
\mathcal{U}_{1}^{\prime}:=\left\{\mathrm{x} \in \mathcal{U}:\left|\mathfrak{u}_{m}(\mathrm{x})\right| \leq 1\right\} \text { and } \mathcal{U}_{2}^{\prime}:=\left\{\mathrm{x} \in \mathcal{U}:\left|\mathfrak{u}_{m}(\mathrm{x})\right|>1\right\}.
$$

Lemma \ref{lemma2} follows that

\begin{align}\label{equation36}
\begin{split}
\int_{\mathcal{U}}\left|\mathfrak{u}_{m}\right|^{\sigma} \log \left(\left|\mathfrak{u}_{m}\right|\right) \mathrm{dx} & =\int_{\mathcal{U}_{1}^{\prime}}\left|\mathfrak{u}_{m}\right|^{\sigma} \log \left(\left|\mathfrak{u}_{m}\right|\right) \mathrm{dx}+\int_{\mathcal{U}_{2}^{\prime}}\left|\mathfrak{u}_{m}\right|^{\sigma} \log \left(\left|\mathfrak{u}_{m}\right|\right) \mathrm{dx} \\
& \leq \frac{1}{\rho} \int_{\mathcal{U}_{1}^{\prime}}\left|\mathfrak{u}_{m}\right|^{\sigma+\rho} \mathrm{dx},
\end{split}
\end{align}

where $\rho$ is chosen sufficiently such that $0<\rho<\frac{s \sigma^{2}}{N}$ and $\sigma<\sigma+\rho<p_{s}^{*}.$ Combining \eqref{equation36} with continuous embedding $W^{s, p}(\mathcal{U}) \hookrightarrow L^{p_{s}^{*}}(\mathcal{U}),$ we get that

$$
\begin{aligned}
\int_{\mathcal{U}}\left|\mathfrak{u}_{m}\right|^{\sigma} \log \left(\left|\mathfrak{u}_{m}\right|\right) \mathrm{dx} & \leq C\left\|\mathfrak{u}_{m}\right\|_{L^{p_{s}^{*}}}^{\mu(\sigma+\rho)}\left\|\mathfrak{u}_{m}\right\|_{L^{\sigma}}^{(1-\mu)(\sigma+\rho)} \\
& \leq C_{1}\left\|\mathfrak{u}_{m}\right\|_{W^{s, p}(\mathcal{U})}^{\mu(\sigma+\rho)}\left\|\mathfrak{u}_{m}\right\|_{L^{\sigma}}^{(1-\mu)(\sigma+\rho)},
\end{aligned}
$$

where $\mu=\frac{s \rho}{N \sigma(\sigma+\rho)} \in(0,1).$ Since $0<\rho<\frac{s \sigma^{2}}{N},$ it follows that $\mu(\sigma+\rho)<\sigma.$ Thanks to Young inequality for any $\varepsilon \in(0,1),$ we get that

\begin{align}\label{equation31}
\int_{\mathcal{U}}\left|\mathfrak{u}_{m}\right|^{\sigma}\left|\mathfrak{v}_{m}\right|^{\sigma} \log \left|\mathfrak{u}_{m} \mathfrak{v}_{m}\right|\mathrm{dx} \leq \int_{\mathcal{U}}\left|\mathfrak{u}_{m}\right|^{\sigma} \log \left(\left|\mathfrak{u}_{m}\right|\right) \mathrm{dx} \leq \varepsilon\left\|\mathfrak{u}_{m}\right\|_{W^{s, p}(\mathcal{U})}^{\sigma}+C(\varepsilon)\left(\left\|\mathfrak{u}_{m}\right\|_{L^{2}}^{2}\right)^{\gamma},
\end{align}

with $\gamma=\frac{(1-\mu)(\sigma+\rho)}{\sigma-\mu(\sigma+\rho)}>1.$\\
- Case 2: $\left\|\mathfrak{u}_{m}\right\|_{W^{s, q}(\mathcal{U})}<1$ and $\left\|\mathfrak{v}_{m}\right\|_{W^{s, q}(\mathcal{U})} \geq 1.$

Similarly, we show that

\begin{align}\label{equation32}
\int_{\mathcal{U}}\left|\mathfrak{u}_{m}\right|^{\sigma}\left|\mathfrak{v}_{m}\right|^{\sigma} \log \left|\mathfrak{u}_{m} \mathfrak{v}_{m}\right|\mathrm{dx} \leq \int_{\mathcal{U}}\left|\mathfrak{v}_{m}\right|^{\sigma} \log \left(\left|\mathfrak{v}_{m}\right|\right) \mathrm{dx} \leq \varepsilon \|\left.\mathfrak{v}_{m}\right|_{W^{s, q}(\mathcal{U})} ^{\sigma}+C(\varepsilon)\left(\left\|\mathfrak{v}_{m}\right\|_{L^{2}}^{2}\right)^{\gamma}.
\end{align}
- Case 3: $\left\|\mathfrak{u}_{m}\right\|_{W^{s, q}(\mathcal{U})}>1$ and $\left\|\mathfrak{v}_{m}\right\|_{W^{s, q}(\mathcal{U})}>1.$

From \eqref{equation27} , we get that

\begin{align}\label{equation28}
\int_{\mathcal{U}}\left|\mathfrak{u}_{m}\right|^{\sigma}\left|\mathfrak{v}_{m}\right|^{\sigma} \log \left|\mathfrak{u}_{m} \mathfrak{v}_{m}\right|\mathrm{dx} \leq C\left\|\mathfrak{u}_{m}\right\|_{2}^{\mu(\sigma+1)}\left\|\mathfrak{v}_{m}\right\|_{2}^{\kappa(\sigma+r)}\left\|\mathfrak{v}_{m}\right\|_{W^{s, p}(\mathcal{U})}^{(1-\kappa)}\left\|\mathfrak{u}_{m}\right\|_{W^{s, q}(\mathcal{U})}^{(1-\mu)(\sigma+r)}.
\end{align}

Recalling the numeral inequality

\begin{align}\label{equation29}
a b c \leq \frac{a^{3}+ b^{3} +c^{3}}{3} \leq a^{3}+ b^{3}+ c^{3}, \text { for all }a, b, c \in \mathbb{R}^{+}.
\end{align}
Using \eqref{equation28} and \eqref{equation29}, we obtaint that

\begin{align}\label{equation30}
\int_{\mathcal{U}}\left|\mathfrak{u}_{m}\right|^{\sigma}\left|\mathfrak{v}_{m}\right|^{\sigma} \log \left|\mathfrak{u}_{m} \mathfrak{v}_{m}\right|\mathrm{dx} \leq C\left\|\mathfrak{u}_{m}\right\|_{2}^{3 \mu(\sigma+1)}\left\|\mathfrak{v}_{m}\right\|_{2}^{3 \kappa(\sigma+r)}+\left\|\mathfrak{v}_{m}\right\|_{W^{s, p}(\mathcal{U})}^{3(\sigma+r)}+\left\|\mathfrak{u}_{m}\right\|_{W^{s, q}(\mathcal{U})}^{3(1-\mu)(\sigma+r)}.
\end{align}

Using the Young inequality and \eqref{equation30}, we get that

\begin{align}\label{equation33}
\int_{\mathcal{U}}\left|\mathfrak{u}_{m}\right|^{\sigma}\left|\mathfrak{v}_{m}\right|^{\sigma} \log \left|\mathfrak{u}_{m} \mathfrak{v}_{m}\right|\mathrm{dx} \leq C\left\|\mathfrak{u}_{m}\right\|_{2}^{6 \mu(\sigma+1)}+C\left\|\mathfrak{v}_{m}\right\|_{2}^{6 \kappa(\sigma+r)}+\varepsilon\left\|\mathfrak{v}_{m}\right\|_{W^{s, p}(\mathcal{U})}^{p}+\varepsilon\left\|\mathfrak{u}_{m}\right\|_{W^{s, q}(\mathcal{U})}^{q}, \quad \text { for all } \varepsilon>0 \text {. }
\end{align}

From \eqref{equation31}, \eqref{equation32}, and \eqref{equation33}, there exist $\lambda>0$ such that

\begin{align}\label{equation35}
\begin{split}
\int_{\mathcal{U}}\left|\mathfrak{u}_{m}\right|^{\sigma}\left|\mathfrak{v}_{m}\right|^{\sigma} \log \left|\mathfrak{u}_{m} \mathfrak{v}_{m}\right|\mathrm{dx} & \leq \varepsilon\left(\left\|\mathfrak{v}_{m}\right\|_{W^{s, p}(\mathcal{U})}^{p}+\varepsilon\left\|\mathfrak{u}_{m}\right\|_{W^{s, q}(\mathcal{U})}^{q}\right)+C_{\varepsilon}\left(\left\|\mathfrak{u}_{m}\right\|_{2}^{2}+\left\|\mathfrak{v}_{m}\right\|_{2}^{2}\right)^{\lambda} \\
& \leq \varepsilon\left(1+\min \left\{\left\|\mathfrak{u}_{m}\right\|_{W^{s, p}(\mathcal{U})}^{p},\left\|\mathfrak{u}_{m}\right\|_{W^{s, p}(\mathcal{U})}^{p(\beta+1)}\right\}+\min \left\{\left\|\mathfrak{v}_{m}\right\|_{W^{s, q}(\mathcal{U})}^{q},\left\|\mathfrak{v}_{m}\right\|_{W^{s, q}(\mathcal{U})}^{q(\beta+1)}\right\}\right)\\
&+C_{\varepsilon}\left(\left\|\mathfrak{u}_{m}\right\|_{2}^{2}+\left\|\mathfrak{v}_{m}\right\|_{2}^{2}\right)^{\lambda} \\
& \leq \varepsilon \frac{p^{\beta}}{\mathcal{K}(1)}\left(\mathcal{K}\left(\left[\mathfrak{u}_{m}\right]_{s}^{p}\right) \int_{\mathcal{U} \times \mathcal{U}} \frac{\left|\mathfrak{u}_{m}(\mathrm{x}, t)-\mathfrak{u}_{m}(\mathrm{y}, t)\right|^{p}}{|\mathrm{x}-\mathrm{y}|^{N+s p}} \mathrm{dx} d \mathrm{y}\right) \\
& +\varepsilon \frac{q^{\beta}}{\mathcal{K}(1)}\left(\mathcal{K}\left(\left[\mathfrak{v}_{m}\right]_{s}^{q}\right) \int_{\mathcal{U} \times \mathcal{U}} \frac{\left|\mathfrak{v}_{m}(\mathrm{x}, t)-\mathfrak{v}_{m}(\mathrm{y}, t)\right|^{q}}{|\mathrm{x}-\mathrm{y}|^{N+s q}} \mathrm{dx} d \mathrm{y}\right)+C_{\varepsilon}\left(\left\|\mathfrak{u}_{m}\right\|_{2}^{2}+\left\|\mathfrak{v}_{m}\right\|_{2}^{2}\right)^{\lambda}.
\end{split}
\end{align}

From \eqref{equation34}, \eqref{equation35}, we have that

$$
\begin{aligned}
\frac{1}{2}\left(\frac{d}{d t}\Vert\mathfrak{u}_{m}\Vert_{2}^{2}+\frac{d}{d t}\Vert \mathfrak{v}_{m} \Vert_{2}^{2}\right) & \leq\left(\frac{\varepsilon q^{\beta}}{\mathcal{K}(1)}-1\right)\mathcal{K}\left(\left[\mathfrak{u}_{m}\right]_{s}^{p}\right) \int_{\mathcal{U} \times \mathcal{U}} \frac{\left|\mathfrak{u}_{m}(\mathrm{x}, t)-\mathfrak{u}_{m}(\mathrm{y}, t)\right|^{p}}{|\mathrm{x}-\mathrm{y}|^{N+s p}} \mathrm{dx} d \mathrm{y}\\
&+ \left(\frac{\varepsilon q^{\beta}}{\mathcal{K}(1)}-1\right)\mathcal{K}\left(\left[\mathfrak{v}_{m}\right]_{s}^{q}\right) \int_{\mathcal{U} \times \mathcal{U}} \frac{\left|\mathfrak{v}_{m}(\mathrm{x}, t)-\mathfrak{v}_{m}(\mathrm{y}, t)\right|^{q}}{|\mathrm{x}-\mathrm{y}|^{N+s q}} \mathrm{dxdy} \\
& +C_{\varepsilon}\left(\left\|\mathfrak{u}_{m}\right\|_{2}^{2}+\left\|\mathfrak{v}_{m}\right\|_{2}^{2}\right)^{\lambda}.
\end{aligned}
$$

By choosing $\varepsilon<\frac{\mathcal{K}(1)}{q^{\beta}},$ it follows from the above inequality that

\begin{align}\label{equation37}
\frac{1}{2}\left(\frac{d}{d t}\Vert\mathfrak{u}_{m}\Vert_{2}^{2}+\frac{d}{d t}\Vert \mathfrak{v}_{m} \Vert_{2}^{2}\right) \leq 2 C_{\varepsilon}\left(\left\|\mathfrak{u}_{m}\right\|_{2}^{2}+\left\|\mathfrak{v}_{m}\right\|_{2}^{2}+1\right)^{\lambda}.
\end{align}

Therefore, we arrive at the following differential inequality

\begin{align}\label{equation3996}
Z_{n}^{\prime}(t) \leq 4 C\left(Z_{n}(t)+1\right)^{\lambda}.
\end{align}

with $Z_{n}(t)=\left\|\mathfrak{u}_{m}\right\|_{2}^{2}+\left\|\mathfrak{v}_{m}\right\|_{2}^{2}.$ We distinguish three cases:\\
- If $\lambda \in(0,1),$ then there exists $C_{1}^{\prime}=\sup _{n \mathbb{N}} Z_{n}(0) \in(0, \infty)$ such that

$$
Z_{n}(t) \leq\left(\left(C_{1}^{\prime}+1\right)^{1-\lambda}+2(1-\lambda) t C\right)^{\frac{1}{1-\lambda}},  \qquad \text { for all } t \in\left(0, t_{0, n}\right).
$$
- If $\lambda=1,$ using \eqref{equation3996}, we arrive at the following inequality

$$
Z_{n}(t) \leq\left(1+c_{1}\right) \exp (2 C t), \qquad \text { for all } t \in\left(0, t_{0, n}\right).
$$
- If $\lambda \geq 1.$ Then, by integrating the inequality \eqref{equation37} with respect to time from 0 to $t,$ we have

\begin{align}\label{equation40}
\left\|\mathfrak{u}_{m}\right\|_{2}^{2}+\left\|\mathfrak{v}_{m}\right\|_{2}^{2} \leq 2 C^{\prime}(\lambda) \int_{0}^{t}\left(\left\|\mathfrak{u}_{m}(z)\right\|_{2}^{2 \lambda} d z+\left\|\mathfrak{v}_{m}(z)\right\|_{2}^{2 \lambda} d z\right)+\lambda \int_{0}^{t} d z.
\end{align}

By using the Gronwall- Bellman inequality, it follows from \eqref{equation40} that there exists $T_{0}>0$ such that

$$
\left\|\mathfrak{u}_{m}\right\|_{2}^{2}+\left\|\mathfrak{v}_{m}\right\|_{2}^{2} \leq C_{6}\left(T_{0}\right) \text {, for all } t \in\left(0, t_{0, n}\right).
$$

We next multiply the first equation in \eqref{equation26}  by $g_{\text {in }}^{\prime}(t)$ and the second equation by $h_{i n}^{\prime}(t),$ summing over $\mathrm{i}$ and integrate with respect to time from 0 to $t,$ we obtain that

\begin{align}\label{equation41}
\int_{0}^{t}\left\|\mathfrak{u}_{m s}(s)\right\|_{2}^{2} d s+\left\|\mathfrak{v}_{m s}\right\|_{2}^{2} d s+\varphi\left(\mathfrak{u}_{m}(t), \mathfrak{v}_{m}(t)\right)=\varphi\left(\mathfrak{u}_{m}(0), \mathfrak{v}_{m}(0)\right), \text { for all } t \in\left(0, T_{0}\right).
\end{align}

Using \eqref{equation41}, we have that

\begin{align}\label{equation43}
\begin{split}
\varphi\left(\mathfrak{u}_{m}(0), \mathfrak{v}_{m}(0)\right) & \geq \varphi\left(\mathfrak{u}_{m}(t), \mathfrak{v}_{m}(t)\right) \\
& =\frac{1}{p}\widehat{\mathcal{K}}\left(\left[\mathfrak{u}_{m}\right]_{p}^{s}\right)+\frac{1}{q}\widehat{\mathcal{K}}\left(\left[\mathfrak{v}_{m}\right]_{q}^{s}\right)+\frac{1}{\sigma^{2}} \int_{\mathcal{U}}\left|\mathfrak{u}_{m}\right|^{\sigma}\left|\mathfrak{v}_{m}\right|^{\sigma} \mathrm{dx}-\frac{1}{\sigma} \int_{\mathcal{U}}\left|\mathfrak{u}_{m}\right|^{\sigma}\left|\mathfrak{v}_{m}\right|^{\sigma} \log \left|\mathfrak{u}_{m} \mathfrak{v}_{m}\right| \mathrm{dx} \\
& \geq \widehat{\mathcal{K}}\left(\left[\mathfrak{u}_{m}\right]_{p}^{s}\right)+\widehat{\mathcal{K}}\left(\left[\mathfrak{v}_{m}\right]_{q}^{s}\right)-\frac{1}{\sigma} \int_{\mathcal{U}}\left|\mathfrak{u}_{m}\right|^{\sigma}\left|\mathfrak{v}_{m}\right|^{\sigma} \log \left|\mathfrak{u}_{m} \mathfrak{v}_{m}\right| \mathrm{dx} \\
& \geq \frac{1}{p(\beta+1)} \mathcal{K}\left(\left[\mathfrak{u}_{m}\right]_{s}^{p}\right) \int_{\mathcal{U} \times \mathcal{U}} \frac{\left|\mathfrak{u}_{m}(\mathrm{x}, t)-\mathfrak{u}_{m}(\mathrm{y}, t)\right|^{p}}{|\mathrm{x}-\mathrm{y}|^{N+s p}} \mathrm{dxdy}\\
&+\frac{1}{q(\beta+1)} \mathcal{K}\left(\left[\mathfrak{v}_{m}\right]_{s}^{p}\right) \int_{\mathcal{U} \times \mathcal{U}} \frac{\left|\mathfrak{v}_{m}(\mathrm{x}, t)-\mathfrak{v}_{m}(\mathrm{y}, t)\right|^{p}}{|\mathrm{x}-\mathrm{y}|^{N+s p}} \mathrm{dxdy}
 -\frac{1}{\sigma} \int_{\mathcal{U}}\left|\mathfrak{u}_{m}\right|^{\sigma}\left|\mathfrak{v}_{m}\right|^{\sigma} \log \left|\mathfrak{u}_{m} \mathfrak{v}_{m}\right| \mathrm{dx} \\
& \geq \frac{1}{q(\beta+1)}\left[\mathcal{K}\left(\left[\mathfrak{u}_{m}\right]_{s}^{p}\right) \int_{\mathcal{U} \times \mathcal{U}} \frac{\left|\mathfrak{u}_{m}(\mathrm{x}, t)-\mathfrak{u}_{m}(\mathrm{y}, t)\right|^{p}}{|\mathrm{x}-\mathrm{y}|^{N+s p}} \mathrm{dxdy}+\mathcal{K}\left(\left[\mathfrak{v}_{m}\right]_{s}^{q}\right) \int_{\mathcal{U} \times \mathcal{U}} \frac{\left|\mathfrak{v}_{m}(\mathrm{x}, t)-\mathfrak{v}_{m}(\mathrm{y}, t)\right|^{p}}{|\mathrm{x}-\mathrm{y}|^{N+s p}} \mathrm{dxdy}\right] \\
& -\frac{1}{\sigma} \int_{\mathcal{U}}\left|\mathfrak{u}_{m}\right|^{\sigma}\left|\mathfrak{v}_{m}\right|^{\sigma} \log \left|\mathfrak{u}_{m} \mathfrak{v}_{m}\right| \mathrm{dx} \\
& \geq\left(\frac{1}{q(\beta+1)}-\frac{\varepsilon q^{\beta}}{\mathcal{K}(1) \sigma}\right) \frac{\mathcal{K}(1)}{q^{\beta}} \min \left\{\left\|\mathfrak{u}_{m}\right\|_{W^{s, p}(\mathcal{U})}^{p}+\left.\left\|\left.\mathfrak{v}_{m}\right|_{W^{s, q}(\mathcal{U})} ^{q},\right\| \mathfrak{u}_{m}\left\|_{W^{s, p}(\mathcal{U})}^{p(\beta+1)}+\right\| \mathfrak{v}_{m}\right|_{W^{s, q}(\mathcal{U})} ^{q(\beta+1)}\right\}-\frac{C^{\gamma}}{\sigma} \\
& \geq\left(\frac{1}{q(\beta+1)}-\frac{\varepsilon q^{\beta}}{\mathcal{K}(1) \sigma}\right) \frac{\mathcal{K}(1)}{q^{\beta}} \min \left\{\left\|\mathfrak{u}_{m}\right\|_{W^{s, p}(\mathcal{U})}^{p}+\|\left.\mathfrak{v}_{m}\right|_{W^{s, q}(\mathcal{U})} ^{q}-1\right\}-\frac{C^{\gamma}}{\sigma} \\
& \geq-\left(\frac{1}{q(\beta+1)}-\frac{\varepsilon q^{\beta}}{\mathcal{K}(1) \sigma}\right) \frac{\mathcal{K}(1)}{q^{\beta}}-\frac{C^{\gamma}}{\sigma}.
\end{split}
\end{align}

On the other hand, we deduce from \eqref{equation42} that $\varphi\left(\mathfrak{u}_{m}(0), \mathfrak{v}_{m}(0)\right) \rightarrow \varphi\left(\mathfrak{u}_{0}(\mathrm{x}), \mathfrak{v}_{0}(\mathrm{x})\right)$ as $m \rightarrow \infty.$ So, we have that

\begin{align}\label{equation44}
\varphi\left(\mathfrak{u}_{m}(0), \mathfrak{v}_{m}(0)\right) \leq C_{7}.
\end{align}

From \eqref{equation41} and \eqref{equation44}, it follows that

\begin{align}\label{equation46}
\begin{split}
\left\|\mathfrak{u}_{m}\right\|_{W^{s, p}(\mathcal{U})}^{p}+\left\|\mathfrak{v}_{m}\right\|_{W^{s, q}(\mathcal{U})}^{q} \leq C_{8}.
\end{split}
\end{align}

From continuous embedding and Lemma \ref{lemma2}, we deduce that

\begin{align}\label{equation47}
\begin{split}
\int_{\mathcal{U}}\left|\mathfrak{u}_{m}\right|^{\sigma}\left|\mathfrak{v}_{m}\right|^{\sigma} \log \left|\mathfrak{u}_{m} \mathfrak{v}_{m}\right|\mathrm{dx} & =\int_{\mathcal{U}_{1}}\left|\mathfrak{u}_{m}\right|^{\sigma}\left|\mathfrak{v}_{m}\right|^{\sigma} \log \left|\mathfrak{u}_{m} \mathfrak{v}_{m}\right|\mathrm{dx}+\int_{\mathcal{U}_{2}}\left|\mathfrak{u}_{m}\right|^{\sigma}\left|\mathfrak{v}_{m}\right|^{\sigma} \log \left|\mathfrak{u}_{m} \mathfrak{v}_{m}\right|\mathrm{dx} \\
& \leq \frac{|\mathcal{U}|}{\sigma \exp (1)}+\frac{1}{\exp (1)} \int_{\mathcal{U}}\left|\mathfrak{u}_{m}\right|^{\sigma}\left|\mathfrak{v}_{m}\right|^{\sigma} \mathrm{dx} \\
& \leq \frac{|\mathcal{U}|}{\sigma \exp (1)}+\frac{1}{\exp (1)} \int_{\mathcal{U}}\left|\mathfrak{u}_{m}\right|^{\sigma+1}\left|\mathfrak{v}_{m}\right|^{\sigma+1} \mathrm{dx} \\
& \leq \frac{|\mathcal{U}|}{\sigma \exp (1)}+\frac{1}{\exp (1)} \int_{\mathcal{U}}\left(\left.\left|\mathfrak{u}_{m}\right|^{2(\sigma+1)}|+| \mathfrak{v}_{m}\right|^{2(\sigma+1)} \mid\right) \mathrm{dx} \\
& \leq \frac{|\mathcal{U}|}{\sigma \exp (1)}+\frac{1}{\exp (1)}\left(M_{1} C_{1}^{2(\sigma+1)}+M_{2} C^{2(\sigma+1)}\right),
\end{split}
\end{align}

with $M_{1}=\sup _{m>1}\left\|\mathfrak{u}_{m}\right\|_{L^{2(\sigma+1)}(\mathcal{U})}^{2(\sigma+1)}, M_{2}=\sup _{m>1}\left\|\mathfrak{v}_{m}\right\|_{L^{2(\sigma+1)}(\mathcal{U})}^{2(\sigma+1)}, C,$ and $C_{1}$ are the embeddings constants of $W^{s, p}(\mathcal{U}) \hookrightarrow$ $L^{2(\sigma+1)}(\mathcal{U})$ and $W^{s, q}(\mathcal{U}) \hookrightarrow L^{2(\sigma+1)}(\mathcal{U}).$ It follows from \eqref{equation43}, \eqref{equation46}, and \eqref{equation47} that there exist a functions $(\mathfrak{u}, \mathfrak{v}) \in$ $L^{\infty}\left(0, T ; W^{s, p}(\mathcal{U})\right) \times L^{\infty}\left(0, T ; W^{s, q}(\mathcal{U})\right)$ and a sub-sequence of $\left(\mathfrak{u}_{m}, \mathfrak{v}_{m}\right)_{m \in \mathbb{N}}$ still denote by itself such that
\begin{equation}\label{eq50}
 \left\{\begin{array}{lll}
\displaystyle
\left(\mathfrak{u}_{m}, \mathfrak{v}_{m}\right) \rightarrow(\mathfrak{u}, \mathfrak{v}), \, \,  \mbox{weakly star in } L^{\infty}\left(0, T ; W^{s, p}(\mathcal{U})\right) \times L^{\infty}\left(0, T ; W^{s, q}(\mathcal{U})\right),\\
\left(\mathfrak{u}_{m t}, \mathfrak{v}_{m t}\right) \rightarrow\left(\mathfrak{u}_{t}, \mathfrak{v}_{t}\right),  \text {  weakly in }  L^{2}\left(0, T ; L^{2}(\mathcal{U})\right) \times L^{2}\left(0, T ; L^{2}(\mathcal{U})\right), \\
\mathcal{K}\left(\left[\mathfrak{u}_{m}\right]_{s}^{p}\right) \int_{\mathcal{U} \times \mathcal{U}} \frac{\left|\mathfrak{u}_{m}(\mathrm{x}, t)-\mathfrak{u}_{m}(\mathrm{y}, t)\right|^{p-2}\left(\mathfrak{u}_{m}(\mathrm{x}, t)-\mathfrak{u}_{m}(\mathrm{y}, t)\right)}{|\mathrm{x}-\mathrm{y}|^{N+s p}} \mathrm{dxdy} \rightarrow \chi_{1},  \text {  weakly in }  L^{\infty}\left(0, T ; L^{\frac{p}{p-1}}(\mathcal{U})\right), \\
\mathcal{K}\left(\left[\mathfrak{v}_{m}\right]_{s}^{q}\right) \int_{\mathcal{U} \times \mathcal{U}} \frac{\left|\mathfrak{v}_{m}(\mathrm{x}, t)-\mathfrak{v}_{m}(\mathrm{y}, t)\right|^{q-2}\left(\mathfrak{v}_{m}(\mathrm{x}, t)-\mathfrak{v}_{m}(\mathrm{y}, t)\right)}{|\mathrm{x}-\mathrm{y}|^{N+s q}} \mathrm{dxdy} \rightarrow \chi_{2},  \text {weakly  in } L^{\infty}\left(0, T ; L^{\frac{q}{q-1}}(\mathcal{U})\right),
\end{array}%
\right.
\end{equation}

Then by using Aubin-Lions-Simon lemma and \eqref{eq50}, we have that

\begin{align}\label{eq1}
\left(\mathfrak{u}_{m}, \mathfrak{v}_{m}\right) \rightarrow(\mathfrak{u}, \mathfrak{v}) \text { in } C\left(\left[0, T_{0}\right] ; L^{r_{1}}(\mathcal{U})\right) \times C\left(\left[0, T_{0}\right] ; L^{r_{2}}(\mathcal{U})\right), \quad \text { for all } 2<r_{1}, r_{2}<\sigma \text { and } \mathrm{x} \in \mathcal{U}.
\end{align}

We followed the same method appearing in Lemma 8 \cite{aberq}, we get that

\begin{align}\label{eq2}
\left|\mathfrak{v}_{m}\right|^{\sigma}\left|\mathfrak{u}_{m}\right|^{\sigma-2} \mathfrak{u}_{m} \log \left|\mathfrak{u}_{m} \mathfrak{v}_{m}\right| \rightarrow|\mathfrak{v}|^{\sigma}|\mathfrak{u}|^{\sigma-2} \mathfrak{u} \log |\mathfrak{u} \mathfrak{v}| \quad \text { weakly in } L^{\frac{\sigma}{\sigma-1}}\left(\mathcal{U} \times\left(0, T_{0}\right)\right),  and
\end{align}
\begin{align}\label{eq3}
\quad\left|\mathfrak{u}_{m}\right|^{\sigma}\left|\mathfrak{v}_{m}\right|^{\sigma-2} \mathfrak{v}_{m} \log \left|\mathfrak{u}_{m} \mathfrak{v}_{m}\right| \rightarrow|\mathfrak{u}|^{\sigma}|\mathfrak{v}|^{\sigma-2} \mathfrak{v} \log |\mathfrak{u}\mathfrak{v}| \quad \text { weakly in } L^{\frac{\sigma}{\sigma-1}}\left(\mathcal{U} \times\left(0, T_{0}\right)\right).
\end{align}

Using \eqref{eq2}, \eqref{eq3}, and \eqref{eq50}, and letting $m \rightarrow \infty$ in \eqref{equation26}  and \eqref{equation42} we have that

\begin{equation}\label{eq8}
\left\{\begin{array}{l}
\left\langle\mathfrak{u}_{t}, w_{1}\right\rangle+\left\langle\chi_{1}, w_{1}\right\rangle=\int_{\mathcal{U}}|\mathfrak{v}|^{\sigma}|\mathfrak{u}|^{\sigma-2} \mathfrak{u} w_{1} \log (|\mathfrak{u}|) \mathrm{d} x \\
\left\langle\mathfrak{v}_{t}, w_{2}\right\rangle+\left\langle\chi_{2}, w_{2}\right\rangle=\int_{\mathcal{U}}|\mathfrak{u}|^{\sigma}|\mathfrak{v}|^{\sigma-2} \mathfrak{v} w_{2} \log (|\mathfrak{u w}|) \mathrm{dx} \\
(\mathfrak{u}(0), \mathfrak{v}(0))=\left(\mathfrak{u}_{0}, \mathfrak{v}_{0}\right)
\end{array}\right.
\end{equation}

for all $\left(w_{1}, w_{2}\right) \in W^{s, p}(\mathcal{U}) \times W^{s, q}(\mathcal{U})$ and $t \in\left(0, T_{0}\right).$ We followed the same method appearing in Lemma 4.3 \cite{aberq}, we get that

$$
\begin{gathered}
\chi_{1}=\mathcal{K}\left([\mathfrak{u}]_{s}^{p}\right) \int_{\mathcal{U} \times \mathcal{U}} \frac{|\mathfrak{u}(\mathrm{x}, t)-\mathfrak{u}(\mathrm{y}, t)|^{p-2}(\mathfrak{u}(\mathrm{x}, t)-\mathfrak{u}(\mathrm{y}, t))}{|\mathrm{x}-\mathrm{y}|^{N+s p}} \mathrm{dxdy}\quad \text { and } \\
\chi_{2}=\mathcal{K}\left([\mathfrak{v}]_{s}^{q}\right) \int_{\mathcal{U} \times \mathcal{U}} \frac{|\mathfrak{v}(\mathrm{x}, t)-\mathfrak{v}(\mathrm{y}, t)|^{q-2}(\mathfrak{v}(\mathrm{x}, t)-\mathfrak{v}(\mathrm{y}, t))}{|\mathrm{x}-\mathrm{y}|^{N+s q}} \mathrm{dx} d \mathrm{y}.
\end{gathered}
$$

Hence, our system \eqref{problem} admits solution $(\mathfrak{u}, \mathfrak{v}) \in Y.$ Now we turn to show inequality \eqref{eq7}. In view of \eqref{equation42}, \eqref{eq2}, and \eqref{eq3}, we have

$$
\varphi\left(\mathfrak{u}_{0 m}, \mathfrak{v}_{0 m}\right) \rightarrow \varphi\left(\mathfrak{u}_{0}, \mathfrak{v}_{0}\right) \text { as } m \rightarrow \infty .
$$

Furthermore, due to the weak lower semicontinuity of the $Y$-norm \eqref{eq2}, and \eqref{eq3}, we conclude

$$
\varphi(\mathfrak{u}, \mathfrak{v}) \leq \lim _{m \rightarrow \infty} \varphi\left(\mathfrak{u}_{m}, \mathfrak{v}_{m}\right)
$$

Letting $m \rightarrow \infty$ in \eqref{equation41}, we obtain

$$
\begin{aligned}
\int_{0}^{t}\left(\left\|\mathfrak{u}_{s}(s)\right\|_{2}^{2}+\left\|\mathfrak{v}_{s}(s)\right\|_{2}^{2}\right) d s+\varphi(\mathfrak{u}(t), \mathfrak{v}(t)) & \leq \liminf _{m \rightarrow \infty}\left[\int_{0}^{t}\left(\left\|\mathfrak{u}_{m s}(s)\right\|_{2}^{2}+\left\|\mathfrak{v}_{m s}(s)\right\|_{2}^{2}\right) d s+\varphi\left(\mathfrak{u}_{m}(t), \mathfrak{v}_{m}(t)\right)\right] \\
& \leq \varphi\left(\mathfrak{u}_{0}, \mathfrak{v}_{0}\right), \quad \text { for all } t \in[0, T].
\end{aligned}
$$
\end{proof}
\section{ Global Existence and Asymptotic Behavior}\label{Global Existence and Asymptotic Behavior}

\begin{theorem} Let $\left(\mathcal{H}_{1}\right)-\left(\mathcal{H}_{2}\right)$ hold. if $\varphi\left(\mathfrak{u}_{0}, \mathfrak{v}_{0}\right)<d_{*}^{\prime}$ and $\psi\left(\mathfrak{u}_{0}, \mathfrak{v}_{0}\right) \geq 0.$ Then, we have
\begin{enumerate}
    \item [1)] $T_{\max }=+\infty$

\item [2)]  the energy functional $\varphi$ satisfies the following estimates

$$
\varphi(\mathfrak{u}(t), \mathfrak{v}(t)) \leq \begin{cases}\varphi\left(\mathfrak{u}_{0}, \mathfrak{v}_{0}\right)\left(\frac{p(\beta+1)}{2+(p(\beta+1)-2) \alpha t}\right)^{\frac{p(\beta+1)}{p(\beta+1)-2}}, & \text { if } q>\frac{2}{(\beta+1)} \\ \varphi\left(\mathfrak{u}_{0}, \mathfrak{v}_{0}\right) e^{1-\mu t}, & \text { if } q \leq \frac{2}{(\beta+1)}\end{cases}
$$
\end{enumerate}
with $\alpha$ and $\mu$ are some positive constants.
\end{theorem}

\begin{proof} Firstly, we prove that $\psi(\mathfrak{u}(t), \mathfrak{v}(t)) \geq 0,$ for all $t \in\left[0, T_{\max }\right).$ Considering the following sets

$$
\mathcal{W}:=\{(\mathfrak{u}, \mathfrak{v}) \in Y: \varphi(\mathfrak{u}, \mathfrak{v})<d, \psi(\mathfrak{u},(\mathfrak{v})>0\} \cup\{0\}, \quad \text { and } \mathcal{M}:=\{(\mathfrak{u}, \mathfrak{v}) \in Y: \varphi(\mathfrak{u},(\mathfrak{v})<d, \psi(\mathfrak{u},(\mathfrak{v})<0\}
$$

Due to $\varphi\left(\mathfrak{u}_{0}, \mathfrak{v}_{0}\right)<d_{*}^{\prime}<d$ and $\psi\left(\mathfrak{u}_{0}, \mathfrak{v}_{0}\right)>0,$ we have $\left(\mathfrak{u}_{0}, \mathfrak{v}_{0}\right) \in \mathcal{W}.$ Using a contradiction for proving $(\mathfrak{u}(t), \mathfrak{v}(t)) \in \mathcal{W},$ for all $t \in\left[0, T_{\max }\right).$ We suppose that there exist a $t_{0} \in\left[0, T_{\max }\right)$ such that $\left(\mathfrak{u}\left(t_{0}\right), \mathfrak{v}\left(t_{0}\right)\right) \in \partial \mathcal{W}.$ Noticing that 0 is an interior of $\mathcal{W}.$ So, we get that $\varphi\left(\mathfrak{u}\left(t_{0}\right), \mathfrak{v}\left(t_{0}\right)\right)=d$ or $\psi\left(\mathfrak{u}\left(t_{0}\right), \mathfrak{v}\left(t_{0}\right)\right)=0$ and $\left(\mathfrak{u}\left(t_{0}\right), \mathfrak{v}\left(t_{0}\right)\right) \neq(0,0).$ Using \eqref{eq7}, we obtain that $\varphi\left(\mathfrak{u}\left(t_{0}\right), \mathfrak{v}\left(t_{0}\right)\right) \leq\left(\mathfrak{u}_{0}, \mathfrak{v}_{0}\right)<d,$ and hence $\psi\left(\mathfrak{u}\left(t_{0}\right), \mathfrak{v}\left(t_{0}\right)\right)=0.$ So, we deduce that $\left(\mathfrak{u}\left(t_{0}\right), \mathfrak{v}\left(t_{0}\right)\right) \in \mathcal{N}.$ Using the definition of $\mathrm{d},$ we have $\varphi\left(\mathfrak{u}\left(t_{0}\right), \mathfrak{v}\left(t_{0}\right)\right) \geq d.$ We arrive that contradiction. Finally, we deduce that $\psi\left(\mathfrak{u}\left(t_{0}\right), \mathfrak{v}\left(t_{0}\right)\right) \geq 0$ for all $t \in\left[0, T_{\max }\right).$ Now, we show that $T_{\max }=\infty.$
 Using Lemma \ref{lemma5}, we get that
 \begin{align}\label{eq25}
 \begin{split}
 d_{*}&>\varphi(\mathfrak{u}_{0}, \mathfrak{u}_{0})\geq
 \varphi\left( \mathfrak{u}(t), \mathfrak{u}(t)\right) \\
 & \geq \varphi\left( \mathfrak{u}(t), \mathfrak{u}(t)\right) -\frac{1}{\sigma}
  \psi\left( \mathfrak{u}(t), \mathfrak{u}(t)\right) \\
  &\geq\left(\frac{1}{q(\beta+1)}-\frac{1}{\sigma}\right)\left[\mathcal{K}\left([\mathfrak{u}(t)]_{p}^{s}\right) \int_{\mathcal{U} \times \mathcal{U}} \frac{|\mathfrak{u}(\mathrm{x},t)-\mathfrak{u}(\mathrm{y},t)|^{p}}{|\mathrm{x}-\mathrm{y}|^{N+s p}} d \mathrm{x} d \mathrm{y}+\mathcal{K}\left([\mathfrak{v}(t)]_{q}^{s}\right) \int_{\mathcal{U} \times \mathcal{U}} \frac{|\mathfrak{v}(\mathrm{x}, t)-\mathfrak{v}(\mathrm{y},t)|^{q}}{|\mathrm{x}-\mathrm{y}|^{N+s q}} d \mathrm{x} d \mathrm{y}\right]\\
   &\geq\left(\frac{1}{q(\beta+1)}-\frac{1}{\sigma}\right)
 \left[ \frac{\mathcal{K}(1)}{p^{\beta}} \min \left\lbrace \Vert \mathfrak{u}\Vert^{p}_{W^{s, p}(\mathcal{U})}, \Vert \mathfrak{u}\Vert^{p(\beta+1)}_{W^{s, p}(\mathcal{U})} \right\rbrace +  \frac{\mathcal{K}(1)}{q^{\beta}} \min \left\lbrace \Vert \mathfrak{v}\Vert^{q}_{W^{s, q}(\mathcal{U})}, \Vert \mathfrak{u}\Vert^{q(\beta+1)}_{W^{s, q}(\mathcal{U})} \right\rbrace \right]\\
 &  \geq\left(\frac{1}{q(\beta+1)}-\frac{1}{\sigma}\right)
 \frac{\mathcal{K}(1)}{q^{\beta}}
 \min \left\lbrace\Vert \mathfrak{u}\Vert^{q}_{W^{s, p}(\mathcal{U})} +\Vert \mathfrak{v}\Vert^{q}_{W^{s, q}(\mathcal{U})} , \Vert \mathfrak{u}\Vert^{q(\beta+1)}_{W^{s, q}(\mathcal{U})}+ \Vert \mathfrak{v}\Vert^{q(\beta+1)}_{W^{s, q}(\mathcal{U})}\right\rbrace\\
   &  \geq\left(\frac{1}{q(\beta+1)}-\frac{1}{\sigma}\right)
 \frac{\mathcal{K}(1)}{q^{\beta}}\min\left\lbrace  \Vert (\mathfrak{u}, \mathfrak{v})\Vert^{q}_{Y},   \Vert (\mathfrak{u}, \mathfrak{v})\Vert^{q(\beta+1)}_{Y}\right\rbrace,
 \end{split}
 \end{align}
 which yields that
 \begin{align}\label{eq26}
 \Vert (\mathfrak{u}, \mathfrak{v})\Vert_{Y} \leq \max \left\lbrace 1, \left(   \frac{d_{*} q^{\beta+1}(\beta+1)\sigma}{\left(\sigma-q(\beta+1) \right)\mathcal{K}(1) }\right)  ^{\frac{1}{q}}\right\rbrace:= C_{1}.
 \end{align}
 In the other words, $\left( \mathfrak{u}(t),  \mathfrak{v}(t)\right) $ is uniformly bounded in time in $W^{s, p}(\mathcal{U})\times W^{s, q}(\mathcal{U}).$ So, we have $T_{\max}=+\infty.$  Using (4) of Lemma \ref{lemma5}, we get that
 \begin{align}\label{eq20}
 \varphi(\mathfrak{u}(t), \mathfrak{v}(t))-\frac{1}{\sigma}
 \psi(\mathfrak{u}(t), \mathfrak{v}(t))\geq  \varepsilon_{*}^{-q(\beta+1)}-C.
 \end{align}

 However, we also have
 \begin{align}\label{eq21}
 \varphi\left( \mathfrak{u}(t), \mathfrak{u}(t)\right) -\frac{1}{\sigma}
  \psi\left( \mathfrak{u}(t), \mathfrak{u}(t)\right)  \leq
 \varphi\left( \mathfrak{u}(t), \mathfrak{u}(t)\right) \leq \varphi(\mathfrak{u}_{0}, \mathfrak{u}_{0}).
 \end{align}
 It follows from \eqref{eq20} and  \eqref{eq21} that

 \begin{align}
 \varepsilon_{*}\geq  \left( \frac{d_{*}}{\varphi(\mathfrak{u}_{0}, \mathfrak{v}_{0})+C}\right) ^{\frac{1}{q(\beta+1)}}.
 \end{align}
 Applying (4) of Lemma \ref{lemma5}, we get that

 \begin{align}\label{eq22}
 \begin{split}
 &\psi(\mathfrak{u}(t), \mathfrak{v}(t))\geq \left(1- \varepsilon_{*}^{-q(\beta+1)-2(\sigma+1)} \right) \left[\mathcal{K}\left([\mathfrak{u}(t)]_{p}^{s}\right) \int_{\mathcal{U} \times \mathcal{U}} \frac{|\mathfrak{u}(\mathrm{x},t)-\mathfrak{u}(\mathrm{y},t)|^{p}}{|\mathrm{x}-\mathrm{y}|^{N+s p}} d \mathrm{x} d \mathrm{y}+\mathcal{K}\left([\mathfrak{v}(t)]_{q}^{s}\right) \int_{\mathcal{U} \times \mathcal{U}} \frac{|\mathfrak{v}(\mathrm{x},t)-\mathfrak{v}(\mathrm{y},t)|^{q}}{|\mathrm{x}-\mathrm{y}|^{N+s q}} d \mathrm{x} d \mathrm{y}\right]\\
 & \geq \left( 1- \left( \frac{d_{*}}{\varphi(\mathfrak{u}_{0}, \mathfrak{v}_{0})}\right) ^{1- \frac{2(\sigma+1)}{q(\beta+1)}} \right) \left[\mathcal{K}\left([\mathfrak{u}(t)]_{p}^{s}\right) \int_{\mathcal{U} \times \mathcal{U}} \frac{|\mathfrak{u}(\mathrm{x},t)-\mathfrak{u}(\mathrm{y},t)|^{p}}{|\mathrm{x}-\mathrm{y}|^{N+s p}} d \mathrm{x} d \mathrm{y}+\mathcal{K}\left([\mathfrak{v}(t)]_{q}^{s}\right) \int_{\mathcal{U} \times \mathcal{U}} \frac{|\mathfrak{v}(\mathrm{x},t)-\mathfrak{v}(\mathrm{y},t)|^{q}}{|\mathrm{x}-\mathrm{y}|^{N+s q}} d \mathrm{x} d \mathrm{y}\right].
 \end{split}
  \end{align}
  Combining  \eqref{eq22} with Lemma \ref{lemma5}, one has
  \begin{align}\label{eq23}
 \begin{split}
  &\varphi(\mathfrak{u}(t), \mathfrak{u}(t))-\frac{1}{\sigma}
  \psi\left( \mathfrak{u}(t), \mathfrak{u}(t)\right) \leq\left(\frac{1}{p(\beta+1)}-\frac{1}{\sigma}\right)\mathcal{K}\left([\mathfrak{u}(t)]_{p}^{s}\right) \int_{\mathcal{U} \times \mathcal{U}} \frac{|\mathfrak{u}(\mathrm{x},t)-\mathfrak{u}(\mathrm{y},t)|^{p}}{|\mathrm{x}-\mathrm{y}|^{N+s p}} d \mathrm{x} d \mathrm{y}\\
  &+  \left(\frac{1}{p(\beta+1)}-\frac{1}{\sigma}\right)\mathcal{K}\left([\mathfrak{v}(t)]_{q}^{s}\right) \int_{\mathcal{U} \times \mathcal{U}} \frac{|\mathfrak{v}(\mathrm{x}, t)-\mathfrak{v}(\mathrm{y},t)|^{q}}{|\mathrm{x}-\mathrm{y}|^{N+s q}} d \mathrm{x} d \mathrm{y}\\
 &\leq \left(\frac{1}{p(\beta+1)}-\frac{1}{\sigma}\right)\left( 1- \left( \frac{d_{*}}{\varphi(\mathfrak{u}_{0}, \mathfrak{v}_{0})}\right) ^{1- \frac{2(\sigma+1)}{q(\beta+1)}} \right)^{-1}  \psi\left( \mathfrak{u}(t), \mathfrak{u}(t)\right).
   \end{split}
  \end{align}
  So, we get that

  \begin{align}\label{eq30}
  \varphi\left( \mathfrak{u}(t), \mathfrak{u}(t)\right) \leq C_{2} \psi\left( \mathfrak{u}(t), \mathfrak{u}(t)\right) ,
  \end{align}
  with $C_{2}:=\left(\frac{1}{p(\beta+1)}-\frac{1}{\sigma}\right)\left( 1- \left( \frac{d_{*}}{\varphi(\mathfrak{u}_{0}, \mathfrak{v}_{0})}\right) ^{1- \frac{2(\sigma+1)}{q(\beta+1)}} \right)^{-1}>0.$
  Now, our contribution is the decay estimates of the functional  $ \varphi\left( \mathfrak{u}(t), \mathfrak{u}(t)\right).$ Using \eqref{eq25} and \eqref{eq26}, we have that
   \begin{align}\label{eq27}
   \varphi\left( \mathfrak{u}(t), \mathfrak{u}(t)\right) &\geq\left(\frac{1}{q(\beta+1)}-\frac{1}{\sigma}\right)
 \frac{\mathcal{K}(1)}{q^{\beta}}\min\left\lbrace  \Vert (\mathfrak{u}, \mathfrak{v})\Vert^{q}_{Y},   \Vert (\mathfrak{u}, \mathfrak{v})\Vert^{q(\beta+1)}_{Y}\right\rbrace\\
 &\geq  C_{3 }\Vert (\mathfrak{u}, \mathfrak{v})\Vert^{q(\beta+1)}_{Y},
     \end{align}
 where $C_{3}:= \left(\frac{1}{q(\beta+1)}-\frac{1}{\sigma}\right)\frac{\mathcal{K}(1)}{q^{\beta}} C_{1}^{q\beta}.$
 Conversely, however. By using  the embedding $W^{s, p}(\mathcal{U})\times W^{s, q}(\mathcal{U})\hookrightarrow L^{2}(\mathcal{U}) \times L^{2}(\mathcal{U}), $ we have that
 \begin{align}\label{eq28}
  \Vert\left( \mathfrak{u}, \mathfrak{v} \right)\Vert_{L^{2}(\mathcal{U})}
\leq C_{4}  \Vert\left( \mathfrak{u}, \mathfrak{v} \right) \Vert_{Y}.
 \end{align}
 Combining \eqref{eq27} with \eqref{eq28}, we get that
  \begin{align}\label{eq29}
   \Vert\left( \mathfrak{u}, \mathfrak{v} \right)\Vert_{L^{2}(\mathcal{U})}\leq C_{5}\left[   \varphi\left( \mathfrak{u}(t), \mathfrak{u}(t)\right)\right]^{\frac{2}{q(\beta+1)}},
   \end{align}

where $ C_{5}=C_{3}^{-\frac{2}{p^{+}(s+1)}} C_{4}^{2}>0.
$
Let be z a number such that $t<z.$ Combining \eqref{eq30} with \eqref{eq29}, we obtain that
\begin{align}
\begin{split}
\int_{t}^{z}  \varphi\left( \mathfrak{u}(h), \mathfrak{v}(h)\right)\mathrm{d} h &\leq C_{2} \int_{t}^{z}   \psi\left( \mathfrak{u}(h), \mathfrak{v}(h)\right)\mathrm{d} h\\
&= \frac{C_{2}}{2}\left[   \left(\left\|\mathfrak{u}(t)\right\|_{2}^{2}+\left\|\mathfrak{v}(t)\right\|_{2}^{2}\right)-    \left(\left\|\mathfrak{u}(z)\right\|_{2}^{2}+\left\|\mathfrak{v}(z)\right\|_{2}^{2}\right) \right]\\
&\leq \frac{C_{2}C_{5}}{2}
\left[ \varphi\left( \mathfrak{u}(t), \mathfrak{u}(t)\right)\right] ^{\frac{2}{q(\beta+1)}}
\end{split}
\end{align}

Letting $z \rightarrow \infty,$ we obtain

\begin{equation}\label{a}
\int_{t}^{\infty} \varphi\left( \mathfrak{u}(h), \mathfrak{v}(h)\right)\mathrm{d} h\leq C_{6} \varphi\left( \mathfrak{u}(t), \mathfrak{u}(t)\right)^{\frac{2}{q(\beta+1)}},
\end{equation}

where $C_{6}=\frac{C_{2} C_{5}}{2}>0.$

It is easy yo show that  $ t \mapsto\varphi\left( \mathfrak{u}(t), \mathfrak{u}(t)\right)$ is non-increasing with respect to $t.$ Next, we separate the following two cases.

 \textbf{Case 1: If $q\leq \frac{2}{\beta+1}.$}\\

Using \eqref{eq7} and  \eqref{a}, we obtain that

$$
\int_{t}^{\infty} \varphi\left( \mathfrak{u}(h), \mathfrak{v}(h)\right)\mathrm{d} h\leq C_{5} \left[\varphi\left( \mathfrak{u}_{0}, \mathfrak{v}_{0}\right)\right]^{\frac{2}{q(\beta+1)}-1} \varphi\left( \mathfrak{u}(t), \mathfrak{u}(t)\right)=\frac{1}{\mu} \varphi\left( \mathfrak{u}(h), \mathfrak{v}(h)\right),
$$

where

$$
\mu=C_{6}^{-1} \left[\varphi\left( \mathfrak{u}_{0}, \mathfrak{v}_{0}\right)\right]^{1-\frac{2}{q(\beta+1)}}>0 .
$$

Then using  Lemma \ref{lem1} to $G(t)=\varphi\left( \mathfrak{u}(t), \mathfrak{u}(t)\right)$ and $\gamma=0,$ we obtain

$$
\varphi\left( \mathfrak{u}(t), \mathfrak{u}(t)\right) \leq \varphi\left(\mathfrak{u}_{0}, \mathfrak{v}_{0}\right) e^{1-\mu t}.
$$

\textbf{Case2: If $q>\frac{2}{\beta+1}.$}

Considering  the function $G(t)=\left[ \varphi\left( \mathfrak{u}(t), \mathfrak{u}(t)\right)\right]  ^{\frac{2}{q(\beta+1)}}.$ Using  \eqref{a}, one has

$$
\int_{t}^{\infty} G^{\frac{q(\beta+1)}{2}}(h) \mathrm{d} h \leq \frac{1}{\eta} G^{\frac{q(\beta+1)}{2}-1}(0) G(t),
$$

where

\begin{align*}
\eta=\frac{\left[\varphi\left( \mathfrak{u}_{0}, \mathfrak{v}_{0}\right)\right]^{1-\frac{2}{q(\beta+1)}}}{C_{6}}>0.
\end{align*}

Then applying  Lemma \ref{lem1} to $F(t)=G(t)$ and $\gamma=\frac{q(\beta+1)}{2}-1>0,$ we obtain

\begin{align*}
G(t) \leq G(0)\left(\frac{q(\beta+1)}{2+\left(q(\beta+1)-2\right) \gamma t}\right)^{\frac{2}{q(\beta+1)-2}}.
\end{align*}

Then, we deduce that

\begin{align*}
\varphi\left( \mathfrak{u}(t), \mathfrak{u}(t)\right) \leq \varphi\left(\mathfrak{u}_{0}, \mathfrak{v}_{0}\right)\left(\frac{q(\beta+1)}{2+\left(q(\beta+1)-2\right) \gamma t}\right)^{\frac{q(\beta+1)}{q(\beta+1)-2}}.
\end{align*}

\end{proof}

\section{ Blow-up Phenomena}\label{Blow-up Phenomena}

In this section, applying Nehari manifold and  the potential well theory combined , we show that the local weak solution of system \eqref{problem} exist Blow-up phenomena.
\begin{theorem}
 Assume that the assumptions $\left(\mathcal{H}_{1}\right)-\left(\mathcal{H}_{2}\right)$ hold. If $\varphi\left(\mathfrak{u}_{0}, \mathfrak{v}_{0}\right)<d_{*}$ and $\psi\left(\mathfrak{u}_{0}, \mathfrak{v}_{0}\right)<0.$ Then, we have $T_{\max }<+\infty.$ Moreover, we also give an upper bound for the maximal existence time as follows

$$
T_{\max } \leq \frac{4(\sigma-1)\left(\left\|\mathfrak{u}_{0}\right\|_{2}^{2}+\left\|\mathfrak{v}_{0}\right\|_{2}^{2}\right)}{\sigma\left(d_{*}-\varphi\left(\mathfrak{u}_{0}, \mathfrak{v}_{0}\right)\right)(\sigma-2)^{2}}
$$
\end{theorem}
\begin{proof} We show that $\psi(\mathfrak{u}(t), \mathfrak{v}(t))<0$ for all $t \in\left[0, T_{\max }\right).$ Indeed, suppose by contradiction that, there exist $t_{0} \in$ $\left(0, T_{\max }\right)$ such that $\psi(\mathfrak{u}(t), \mathfrak{v}(t))<0$ for all $t \in\left[0, t_{0}\right)$ and $\psi\left(\mathfrak{u}\left(t_{0}\right), \mathfrak{v}\left(t_{0}\right)\right)=0.$ Testing by $(\mathfrak{u}, \mathfrak{v})$ in \eqref{eq8}, we get that

$$
\begin{aligned}
\frac{d}{d t}\|\mathfrak{u}(t)\|_{2}^{2}+\frac{d}{d t}\|\mathfrak{v}(t)\|_{2}^{2} & =2\left(\left\langle\mathfrak{u}_{t}, \mathfrak{u}\right\rangle+\left\langle\mathfrak{v}_{t}, \mathfrak{v}\right\rangle\right) \\
&=-2 \psi\left(\mathfrak{u}(t), \mathfrak{v} (t)\right)> 0, \quad \text { for all } t \in\left[0, t_{0}\right).
\end{aligned}
$$

So, the map $t \rightarrow\|\mathfrak{u}(t)\|_{2}^{2}+\|\mathfrak{v}(t)\|_{2}^{2}$ is strictly increasing on $\left[0, t_{0}\right].$ Therefore, $\left\|\mathfrak{u}\left(t_{0}\right)\right\|_{2}^{2}+\left\|\mathfrak{v}\left(t_{0}\right)\right\|_{2}^{2}.$ So, we get that $\left(\mathfrak{u}\left(t_{0}\right), \mathfrak{v}\left(t_{0}\right)\right) \in \mathcal{N}.$ Using the definition of $\mathrm{d},$ we obtain

\begin{align}\label{eq9}
\varphi\left(\mathfrak{u}\left(t_{0}\right), \mathfrak{v}\left(t_{0}\right)\right) \geq d.
\end{align}

However, it follows from \eqref{eq7}

\begin{align}\label{eq10}
\varphi\left(\mathfrak{u}\left(t_{0}\right), \mathfrak{v}\left(t_{0}\right)\right) \leq \varphi\left(\mathfrak{u}_{0}, \mathfrak{v}_{0}\right) < d_{*} \leq d.
\end{align}

Using \eqref{eq9} and \eqref{eq10}, we arrive at a contradiction. So, we have $\psi(\mathfrak{u}(t), \mathfrak{v}(t))<0$ for all $t \in\left[0, T_{\max }\right).$ Now, we give an upper bound of $T_{\max }.$ Fixing $T \in\left(0, T_{\max }\right)$ and considering the function $\mathrm{L}$ defined by

$$
L(t)=\int_{0}^{t}\left(\left\|\mathfrak{u}_{s}(s)\right\|_{2}^{2}+\left\|\mathfrak{v}_{s}(s)\right\|_{2}^{2}\right) d s+(T-t)\left(\left\|\mathfrak{u}_{0}\right\|_{2}^{2}+\left\|\mathfrak{v}_{0}\right\|_{2}^{2}\right)+(a t+b)^{2}, \text { for all } t \in[0, T]
$$

with $\mathrm{a}, \mathrm{b}$ are positive constents. By a direct computation, we have

$$
\begin{aligned}
L^{\prime}(t) & =\left\|\mathfrak{u}_{s}(s)\right\|_{2}^{2}+\left\|\mathfrak{v}_{s}(s)\right\|_{2}^{2}-\left(\left\|\mathfrak{u}_{0}\right\|_{2}^{2}+\left\|\mathfrak{v}_{0}\right\|_{2}^{2}\right)+2 a(a t+b) \\
& =2\left[\int_{0}^{t}\left(\left\langle\mathfrak{u}_{s}(s), \mathfrak{u}(s)\right\rangle+\left\langle\mathfrak{v}_{s}(s), \mathfrak{v}(s)\right\rangle\right) d s+a(a t+b)\right] \\
& >0.
\end{aligned}
$$

and

\begin{align}\label{eq11}
\begin{split}
L^{\prime \prime}(t) & =2\left(\left\langle\mathfrak{u}_{s}(s), \mathfrak{u}(s)\right\rangle+\left\langle\mathfrak{v}_{s}(s), \mathfrak{v}(s)\right\rangle+a^{2}\right) \\
& =2(-\psi(\mathfrak{u}(t), \mathfrak{v}(t)))+a^{2}>0.
\end{split}
\end{align}

Thanks to Cauchy-Schwarz inequality, we get that

\begin{align}\label{eq12}
\begin{split}
\left(L^{\prime}(t)\right)^{2} & =4\left[\int_{0}^{t}\left(\left\langle\mathfrak{u}_{s}(s), \mathfrak{u}(s)\right\rangle+\left\langle\mathfrak{v}_{s}(s), \mathfrak{v}(s)\right\rangle\right) d s+a(a t+b)\right]^{2} \\
& \leq 4\left[\sqrt{\int_{0}^{t}\left\|\mathfrak{u}_{s}(s)\right\|_{2}^{2} d s} \sqrt{\int_{0}^{t}\left\|\mathfrak{u}(s)_{2}^{2}\right\| d s}+\sqrt{\int_{0}^{t}\left\|\mathfrak{v}_{s}(s)\right\|_{2}^{2} d s} \sqrt{\int_{0}^{t}\|\mathfrak{v}(s)\|_{2}^{2} d s}+a(a t+b)\right]^{2} \\
& \leq 4\left[\int_{0}^{t}\left(\|\mathfrak{u}(s)\|_{2}^{2}+\|\mathfrak{v}(s)\|_{2}^{2}\right) d s+(a t+b)^{2}\right]\left[\int_{0}^{t}\left(\left\|\mathfrak{u}_{s}(s)\right\|_{2}^{2}+\left\|\mathfrak{v}_{s}(s)\right\|_{2}^{2}\right) d s+a^{2}\right] \\
& =4 L(t)\left[\int_{0}^{t}\left(\left\|\mathfrak{u}_{s}(s)\right\|_{2}^{2}+\left\|\mathfrak{v}_{s}(s)\right\|_{2}^{2}\right) d s+a^{2}\right].
\end{split}
\end{align}

On the other hand, using (3) of Lemma  \ref{lemma5}, we have

\begin{align}\label{eq13}
\begin{split}
\varphi(\mathfrak{u}(t), \mathfrak{v}(t))-\frac{1}{\sigma} \psi(\mathfrak{u}(t), \mathfrak{v}(t))>d_{*}.
\end{split}
\end{align}

Using \eqref{eq7},  \eqref{eq11},  \eqref{eq12}, and  \eqref{eq13} we obtain that

\begin{align}\label{eq14}
\begin{split}
L^{\prime \prime}(t)= & 2\left(-2 \psi(\mathfrak{u}(t), \mathfrak{v}(t))+a^{2}\right) \\
& \geq 2\left[\sigma\left(d_{*}-\varphi(\mathfrak{u}(t), \mathfrak{v}(t))\right)+a^{2}\right] \\
& \geq 2\left[\sigma \int_{0}^{t}\left(\|\mathfrak{u}(s)\|_{2}^{2}+\|\mathfrak{v}(s)\|_{2}^{2}\right) d s+\sigma\left(d_{*}-\varphi\left(\mathfrak{u}_{0}, \mathfrak{v}_{0}\right)\right)+a^{2}\right] \\
& \geq\left[\sigma\left(\left(\frac{\left.L^{\prime}(t)\right)^{2}}{4 L(t)}-a^{2}\right)+\sigma\left(d_{*}-\varphi\left(\mathfrak{u}_{0}, \mathfrak{v}_{0}\right)\right)+a^{2}\right] .\right.
\end{split}
\end{align}

Usin \eqref{eq14}, we have that

\begin{align}\label{eq15}
\begin{split}
L^{\prime \prime}(t) L(t)-\frac{\sigma}{2}\left(L^{\prime}(t)\right)^{2} \geq 2 L(t)\left[\sigma\left(d_{*}-\varphi\left(\mathfrak{u}_{0}, \mathfrak{v}_{0}\right)+a^{2}\right]\right.
\end{split}
\end{align}

Putting

\begin{align}\label{eq16}
\begin{split}
0<a<\sqrt{\frac{\sigma\left(d_{*}-\varphi\left(\mathfrak{u}_{0}, \mathfrak{v}_{0}\right)\right)}{\sigma}}.
\end{split}
\end{align}

Combining \eqref{eq15}) with \eqref{eq16}, we have the above inequality

$$
L^{\prime \prime}(t)-\frac{\sigma}{2}\left(L^{\prime}(t)\right)^{2} \geq 0
$$

By a direct computation, we have $L(0), L^{\prime}(0)>0.$ Using Lemma \ref{lemma2.2}, we have that $T \leq \frac{L(0)}{\left(\frac{\sigma}{2}-1\right) L^{\prime}(0)}=\frac{T\left(\left\|u_{0}\right\|_{2}^{2}+\left\|\mathfrak{v}_{0}\right\|_{2}^{2}\right)+b^{2}}{a b\left(\frac{\sigma}{2}-1\right)},$ which implies that $T \leq \frac{b^{2}}{a b\left(\frac{\sigma}{2}-1\right)-\left(\left\|\mathfrak{u}_{0}\right\|_{2}^{2}+\left\|\mathfrak{u}_{0}\right\|_{2}^{2}\right)}:=\xi(b),$ for all $b \in\left(\frac{\left\|\mathfrak{u}_{0}\right\|_{2}^{2}+\left\|\mathfrak{v}_{0}\right\|_{2}^{2}}{a\left(\frac{\sigma}{2}-1\right)}, \infty\right).$ Minimizing $\xi(b)$ for all $b \in\left(\frac{\left\|\mathfrak{u}_{0}\right\|_{2}^{2}+\left\|\mathfrak{v}_{0}\right\|_{2}^{2}}{a\left(\frac{\sigma}{2}-1\right)}, \infty\right).$ So, we get that

$$
T \leq \frac{4\left(\left\|\mathfrak{u}_{0}\right\|_{2}^{2}+\left\|\mathfrak{v}_{0}\right\|_{2}^{2}\right)}{\left(a\left(\frac{\sigma}{2}-1\right)\right)^{2}}:=\varsigma(a), \text { for all } 0<a<\sqrt{\frac{\sigma\left(d_{*}-\varphi\left(\mathfrak{u}_{0}, \mathfrak{v}_{0}\right)\right)}{\sigma}}.
$$

Minimizing $\varsigma(a)$ and letting $T \rightarrow T_{\max },$ we deduce that

$$
T_{\max } \leq \frac{4(\sigma-1)\left(\left\|\mathfrak{u}_{0}\right\|_{2}^{2}+\left\|\mathfrak{v}_{0}\right\|_{2}^{2}\right)}{\sigma\left(d_{*}-\varphi\left(\mathfrak{u}_{0}, \mathfrak{v}_{0}\right)\right)(\sigma-2)^{2}}.
$$
\end{proof}

\subsection*{Aknoledgements}

The third author is a member of INdAM (Istituto Nazionale di Alta Matematica \textquotedblleft Francesco Severi\textquotedblright ) Research group GNAMPA (Gruppo Nazionale
per l'Analisi Matematica, la Probabilit`a e le loro Applicazioni). The third author likes to thank Faculty of Fundamental Science, Industrial University of Ho Chi Minh City, Vietnam, for the opportunity to work in it.

\section*{Declarations}

\subsection*{Conflict of Interests}
The authors declare to have no conflict of interests.

\subsection*{Authors contributions}
The authors declare that their contributions are equal.

 \subsection*{Data availabilty}
This manuscript contains no associated data.

\subsection*{Funding}
Not Applicable

\subsection*{Orcid Data}

\noindent
ORCID ID of Ahmed Aberqi~~~~~~~ https://orcid.org/0000-0003-3599-9099

 \noindent
ORCID ID of Abdesslam Ouaziz~ https://orcid.org/0000-0002-5523-5654

\noindent
 ORCID ID of Maria Alessandra Ragusa\, https://orcid.org/0000-0001-6611-6370

%\bibliographystyle{...}
%\bibliography{...}

\begin{thebibliography}{9}
%Please refer to the journal's website for the corresponding reference style.

\bibitem{aberq}
  Aberqi A., Benslimane O., Ouaziz A., Repovs D.D.:
Fractional Sobolev spaces with kernel function on compact Riemannian manifolds. Mediterranean Journal of Mathematics. \textbf{21} (1), 6  (2024)

\bibitem{Elmassoudi}
  Aberqi A., Elmassoudi M., Hammoumi M.: Descret solution for a
nonlinear parabolic equations with diffusion terms in museilak-spaces. Mathematical Modeling and Computing. \textbf{8} (4), 584-600 (2021)

\bibitem{aberq}
Aberqi A., Ouaziz A.:
 Morse's theory and local linking for a fractional $\left(p _{1} (x.,), p _{2} (x.,)\right):$ Laplacian problems on compact manifolds,
  Journal of Pseudo-Differential Operators and Applications.  \textbf{14} (3), 41 (2023)

\bibitem{Aboulaich}
Aboulaich R., Meskine D., Souissi A.: New diffusion models in image processing. Computers  Mathematics with Applications. \textbf{56} (4), 874-882 (2008)

\bibitem{Applebaum}
Applebaum D.: L{\'e}vy processes-from probability to finance and quantum groups. Notices of the AMS. \textbf{51} (11), 1336-1347 (2004)
\bibitem{Antontsev}
  Antontsev, S.N., Rodrigues, J.F.: On stationary thermo-rheological viscous flows. Annali- Universita Di Ferrara. \textbf{752} (1), 19 (2006)

  \bibitem{autuori}
Autuori G., Pucci P., Salvatori M. C.:
Global nonexistence for nonlinear Kirchhoff systems,  Archive for Rational Mechanics and Analysis.  \textbf{196},  489--516 (2010)

  \bibitem{Ball}
Ball J.:
 Remarks on blow-up and nonexistence theorems for nonlinear evolution equations, Quart. J. Math. Oxford Ser. \textbf{28}, 473-486 (1977)

 \bibitem{bucur}
Bucur C., Valdinoci E.: Nonlocal diffusion and applications, 20. Ed. Springer. (2016)

\bibitem{adams}
 Brezis H.: Functional Analysis, Sobolev Spaces and Partial Differential Equations, Ed. Springer New York, (2010).
\bibitem{Caffarelli}
Caffarelli, L.: Non-local equations, drifts and games. In: Nonlinear Partial Differential Equations. Abel Symposia.  \textbf{7},  37-52 (2012)

\bibitem{Chipot}
Chipot M., Valente V., Caffarelli G. V.:  Remarks on a nonlocal problems involving the Dirichlet energy, Rend. Sem. Math. Univ. Padova.  \textbf{110}, 199-220  (2003)

%\bibitem{Bollt}
 % title={Graduated adaptive image denoising: local compromise between total variation and isotropic diffusion},
 % author={Bollt, Erik M and Chartrand, Rick and Esedo-lu, Selim and Schultz, Pete and Vixie, Kevin R},
  %journal={Advances in Computational Mathematics},
  %volume={31},
  %pages={61--85},
  %year={2009},
  %publisher={Springer}
%}

\bibitem{Chen}
Chen Y., Levine S., Rao M.: Variable exponent, linear growth functionals in image restoration. SIAM Journal on Applied Mathematics. \textbf{66} (4), 1383-1406 (2006)

\bibitem{iezza}
Di Nezza, E,. Palatucci G., Valdinoci E.: Hitchhker's guide to the fractional Sobolev spaces. Bull. Sci. Math. \textbf{136} (5), 519-527  (2012)

\bibitem{dAvenia}
D'Avenia P., Squassina M., Zenari M., Fractional logarithmic Schrodinger equations. Math. Methods Appl. Sci. \textbf{38}, 5207-5216  (2015)

\bibitem{Evans}
 Evans L.C.: Partial Differential Equations,  Graduate Studies in Mathematics. American Mathematical Society, Providence. (1998)

\bibitem{Guo}
 Guo L., Zhang B.,  Zhang Y.: Fractional p-Laplacian equations on Riemannian manifolds. Electron. J. Differ. Equations.  \textbf{156}, 1-17  (2018)

\bibitem{Ferreira}
 Ferreira L. C., Queiroz O. S. A.: singular parabolic equation with logarithmic nonlinearity and $L^{p}$-initial data, J. Differ. Equations. \textbf{249} (2), 349-365 (2010)

\bibitem{Ghisi}
 Ghisi M.,  Gobbino  M.: Hyperbolic-parabolic singular perturbation for middly degenerate Kirchhoff equations: time-decay estimates, J. Differ. Equaqions. \textbf{245}, 2979-3007 (2008)

\bibitem{khaldi}
 Khaldi A., Ouaoua A., Maouni M.: Global existence and stability of solution for a nonlinear Kirchhoff type reaction-diffusion equation with variable exponents. Mathematica Bohemica. \textbf{147} (4), 471-484 (2022)

\bibitem{kirchhoff}
 Kirchhoff G., Hensel K.: Vorlesungen uber mathematische physik, Druck und Verlag von BG Teubner. (1883)

\bibitem{Komornik}
 Komornik V.: Exact Controllability and Stabilization, RAM: Research in Applied Mathematics. Masson, Paris; John Wiley, Ltd., Chichester. (1994)

\bibitem{Lee}
Lee J., Kim J.-M., Kim Y.-H.: Existence and multiplicity of solutions for kirchhoff-schrodinger type equations involving $p(x)-$laplacian on the entire space RN. Nonlinear Analysis: Real World Applications. \textbf{45}, 620-649 (2019)

\bibitem{levine}
Levine H. A,: Instability and nonexistence of global solutions to nonlinear wave equations of the form  $ {P}u_{tt}=-{A}u+ {F} (u).$ Trans Amer. Math. Soc. \textbf{192}, 1-21 (1974)

\bibitem{lions}
Lions J.-L.: Quelques m{\'e}thodes de r{\'e}solution des problemes aux limites non lin{\'e}aires, Dunod. 88-233 (1969)

\bibitem{Liu}
 Liu D.: On a $p-$ Kirchhoff equation via fountain theorem and dual fountain theorem. Nonlinear Analysis: Theory, Methods  Applications. \textbf{72} (1), 302-308 (2010)

\bibitem{Liu1}
 Liu Y., Zhao J.: On the potential wells and applications to semilinear hyperbolic equations and parabolic equations, Nonlinear Anal. \textbf{64} (12), 2665-2687
(2006)

\bibitem{Han}
Han Y., Li Q.: Threshold results for the existence of global and blow-up solutions to kirchhoff equations with arbitrary initial energy. Computers  Mathematics with Applications. \textbf{75} (9), 3283-3297 (2018)

\bibitem{Komornik}
Keyantuo V., Tebou L., Warma M.: A gevrey class semigroup for a thermoelastic plate model with a fractional laplacian: Between the
Euler -Bernoulli and Kirchhoff models. Discrete  Continuous Dynamical Systems: Series A. \textbf{40} (5) (2020)

\bibitem{Le}
Le C.N., Le X. T..: Global solution and blow-up for a class of p-Laplacian evolution equation with logarithmic
nonlinearity, Acta. Appl Math. \textbf{151}, 149-169 (2017)

\bibitem{Pan}
Pan N., Zhang B., Cao J.: Degenerate Kirchhoff-type diffusion problems involving the fractional p-Laplacian, Nonlinear Anal. Real World Appl. \textbf{37}, 56-70  (2017)

\bibitem{Aberqi}
Ouaziz A., Aberqi A.: Infinitely many solutions to a Kirchhoff-type equation involving logarithmic nonlinearity via Morse's theory. Bolet{\'i}n de la Sociedad Matem{\'a}tica Mexicana. \textbf{30} (1), 10 (2024)

\bibitem{aberqi2022}
Ouaziz A., Aberqi A., Benslimane O.: On some weighted fractional $ p(., .)-$Laplacian problems. Palestine Journal of Mathematics.  \textbf{12} (4), (2023)

\bibitem{Payne}
Payne L.E., Sattinger D.H.:  Saddle points and instability of nonlinear hyperbolic equations, Israel J. Math. \textbf{22}, 273-303  (1975)

\bibitem{Rajagopal}
Rajagopal K.R., R\u uzick\v a M.: Mathematical modeling of electrorheological materials. Continuum mechanics and thermodynamics. \textbf{13} (1), 59-78 (2001)

\bibitem{Samko}
Samko S.: On progress in the theory of Lebesgue spaces with variable exponents: maximal and singular operators. Integral Transforms and
Special Functions.  \textbf{16} (5-6), 461-482 (2005)

\bibitem{Sattinger}
Sattinger  D.H.: On global solution of nonlinear hyperbolic equations, Arch. Rat. Mech. Anal. \textbf{30}, 148-172 (1968)

\bibitem{Tan}
Tan Z.: Global solution and blow-up of semilinear heat equation with critical Sobolev exponent, Commun. Partial Differ. Equ. \textbf{26}, 717-741 (2001)

\bibitem{Tian}
Tian  S.: Multiple solutions for the semilinear elliptic equations with the sign-changing logarithmic nonlinearity, J. Math. Anal. Appl. \textbf{454}, 816-828 (2017)

\bibitem{Zloshchastiev}
Zloshchastiev K. G.: Logarithmic nonlinearity in the theories of quantum gravity: origin of time and
observational consequences. Grav. Cosmol. \textbf{16}, 288-297 (2017)



%%%%%%%%%%%%%%%%%%%%%%%%%%%%%%%%%%%%



\end{thebibliography}

\end{document}